\pdfoutput=1 
\documentclass[12pt,reqno]{amsart}
\usepackage{amssymb,latexsym,graphicx,amscd}
\usepackage{lscape}
\usepackage[all]{xy}
\usepackage{xcolor}       
\usepackage{epic,eepic}
\usepackage{xspace}
\usepackage{mathrsfs}
\usepackage{setspace}
\usepackage{float}

\usepackage{tikz}
\usepackage{graphicx}
\usetikzlibrary{calc}
\usetikzlibrary{matrix,arrows}
\usepackage{mathtools} 
\tikzset{desc/.style=auto}
\tikzset{diagmat/.style={matrix of math nodes, row sep=3em, column sep=3em, text height=1.5ex, text depth=0.25ex}}
\tikzset{diagpath/.style={->, font=\scriptsize}}
\usepackage{pgfplots}
\pgfplotsset{my style/.append style={axis x line=middle, axis y line=
middle, axis equal }}
\def\centerarc[#1](#2)(#3:#4:#5){ \draw[#1] ($(#2)+({#5*cos(#3)},{#5*sin(#3)})$) arc (#3:#4:#5); }
\usepackage{mdframed}

\newcommand{\greenx}[1]{\textcolor{green}{#1}}

\usepackage{caption}
\DeclareCaptionLabelFormat{cont}{#1~#2\alph{ContinuedFloat}}
\usepackage{hyperref}
\hypersetup{
    colorlinks,
    linkcolor={blue},
    citecolor={blue},
    urlcolor={blue!80!black}
}

\newtheorem{Thm}{Theorem}[section]
\newtheorem{Lem}[Thm]{Lemma}

\newtheorem{Cor}[Thm]{Corollary}
\newtheorem{Prop}[Thm]{Proposition}

\theoremstyle{definition}           
\newtheorem{Rem}[Thm]{Remark}
\newtheorem{Def}[Thm]{Definition}   
\numberwithin{equation}{section}    

\newcommand{\Z}{\mathbb{Z}}

\newcommand{\N}{\mathbb{N}}
\newcommand{\C}{\mathbb{C}}

\newcommand{\df}{\colon}

\newcommand{\cA}{{\mathcal A}}
\newcommand{\cB}{{\mathcal B}}
\newcommand{\cC}{{\mathcal C}}

\newcommand{\cF}{{\mathcal F}}
\newcommand{\cG}{{\mathcal G}}

\newcommand{\cM}{{\mathcal M}}
\newcommand{\cN}{{\mathcal N}}
\newcommand{\cO}{{\mathcal O}}
\newcommand{\cP}{{\mathcal P}}

\newcommand{\cS}{{\mathcal S}}

\newcommand{\cU}{{\mathcal U}}

\newcommand{\cX}{{\mathcal X}}

\newcommand{\bbS}{{\mathbb S}}
\newcommand{\bbM}{{\mathbb M}}

\newcommand{\m}{\mathfrak{m}}

\newcommand{\bd}{\mathbf{d}}
\newcommand{\be}{\mathbf{e}}
\newcommand{\bg}{\mathbf{g}}   

\newcommand{\bs}{{\mathbf s}}

\newcommand{\bv}{{\mathbf v}}

\newcommand{\Del}{\Delta}      

\newcommand{\gam}{\gamma}      
\newcommand{\eps}{\varepsilon} 
\newcommand{\si}{\sigma}

\newcommand{\la}{\lambda}

\newcommand{\sfA}{\mathsf{A}}

\newcommand{\ux}{\underline{x}}
\newcommand{\hy}{\hat{y}}
\newcommand{\uhy}{\underline{\hat{y}}}

\newcommand{\tQ}{\widetilde{Q}}
\newcommand{\wtA}{\widetilde{A}}
\newcommand{\tA}{\widetilde{\mathsf{A}}}

\newcommand{\rk}{\operatorname{rank}}

\newcommand{\Lam}{\operatorname{Lam}}

\newcommand{\SM}{(\bbS,\bbM)}

\newcommand{\cross}{\operatorname{cross}}    
\newcommand{\md}{\operatorname{mod}}
\newcommand{\ind}{\operatorname{ind}}
\newcommand{\rad}{\operatorname{rad}}
\newcommand{\smooth}{\operatorname{smooth}}

\newcommand{\CaCh}{\operatorname{CC}}    

\newcommand{\nil}{\operatorname{nil}}
\newcommand{\rep}{\operatorname{rep}}
\newcommand{\decrep}{\operatorname{decmod}}

\newcommand{\pdim}{\operatorname{proj.dim}}
\newcommand{\idim}{\operatorname{inj.dim}}

\newcommand{\dimv}{\underline{\dim}}
\newcommand{\rd}{{\rm red}}

\newcommand{\proj}{{\rm proj}}
\newcommand{\Hom}{\operatorname{Hom}}
\newcommand{\sHom}{\operatorname{\underline{Hom}}}

\newcommand{\Ext}{\operatorname{Ext}}
\newcommand{\ext}{\operatorname{ext}}
\newcommand{\End}{\operatorname{End}}

\newcommand{\length}{\operatorname{length}}

\newcommand{\Ima}{\operatorname{Im}}
\newcommand{\Ker}{\operatorname{Ker}}

\newcommand{\op}{{\operatorname{op}}}    

\newcommand{\Int}{\operatorname{Int}}

\newcommand{\irr}{\operatorname{Irr}}
\newcommand{\decirr}{\operatorname{decIrr}}

\newcommand{\bsm}{\begin{smallmatrix}}
\newcommand{\esm}{\end{smallmatrix}}

\newcommand{\bbm}{\begin{matrix}}
\newcommand{\ebm}{\end{matrix}}

\newcommand{\bbsm}{\left(\begin{smallmatrix}}
\newcommand{\besm}{\end{smallmatrix}\right)}

\newcommand{\bbbm}{\left(\begin{matrix}}
\newcommand{\bebm}{\end{matrix}\right)}

\newcommand{\GL}{\operatorname{GL}}
\newcommand{\Gr}{\operatorname{Gr}}
\newcommand{\Spec}{\operatorname{Spec}}

\newcommand{\ra}{\rightarrow}

\newcommand{\red}{\color{red}}
\newcommand{\blue}{\color{blue}}
\newcommand{\green}{\color{green}}

\newcommand{\redx}[1]{\textcolor{red}{#1}}   

\newcommand{\sr}{\tau}

\begin{document}

\today

\title[Schemes of modules over gentle algebras]{Schemes of modules over gentle algebras and 
laminations of surfaces}


\author{Christof Gei{\ss}}
\address{Christof Gei{\ss}\newline
Instituto de Matem\'aticas\newline
Universidad Nacional Aut{\'o}noma de M{\'e}xico\newline
Ciudad Universitaria\newline
04510 M{\'e}xico D.F.\newline
M{\'e}xico}
\email{christof.geiss@im.unam.mx}

\author{Daniel Labardini-Fragoso}
\address{Daniel Labardini-Fragoso\newline
Instituto de Matem\'aticas\newline
Universidad Nacional Aut{\'o}noma de M{\'e}xico\newline
Ciudad Universitaria\newline
04510 M{\'e}xico D.F.\newline
M{\'e}xico}
\email{labardini@im.unam.mx}

\author{Jan Schr\"oer}
\address{Jan Schr\"oer\newline
Mathematisches Institut\newline
Universit\"at Bonn\newline
Endenicher Allee 60\newline
53115 Bonn\newline
Germany}
\email{schroer@math.uni-bonn.de}



\begin{abstract}
We study the affine schemes of modules over gentle algebras.
We describe the smooth points of
these schemes, and we also analyze their irreducible
components in detail.
Several of our results generalize formerly known results, e.g. by dropping
acyclicity, and by incorporating band modules.
A special class of gentle algebras are
Jacobian algebras arising from 
triangulations of unpunctured marked surfaces.
For these we obtain a bijection between the set of
generically $\tau$-reduced decorated irreducible components and the
set of laminations of the surface.
As an application, we get that the set of bangle functions
(defined by Musiker-Schiffler-Williams)
in the upper cluster algebra associated with the surface 
coincides with the set of generic Caldero-Chapoton functions 
(defined by Gei{\ss}-Leclerc-Schr\"oer). 
\end{abstract}

\maketitle

\setcounter{tocdepth}{1}

\tableofcontents

\parskip2mm



\section{Introduction and main results}\label{sec:intro}


\subsection{Overview}
We study some geometric aspects of the representation theory
of gentle algebras.
This class of finite-dimensional algebras was defined by Assem and Skowro\'nski \cite{AS}, who were classifying the iterated tilted algebras of path algebras of
extended Dynkin type $\widetilde{A}$.
Gentle algebras are special biserial, which implies that their 
module categories can be described combinatorially,
see \cite{WW} and also \cite{BR}.

The irreducible components of the affine schemes of modules over gentle algebras are easy to classify 
(see Proposition~\ref{prop:gentlecomponents}).
As a first main result, we describe all smooth points of
these schemes, and we show that
most components are generically reduced.

A special class of gentle algebras are
Jacobian algebras arising from 
triangulations of unpunctured marked surfaces $\SM$.
For these
we obtain a bijection between the set of
generically $\tau$-reduced decorated irreducible components and the
set of laminations of the surface.
This bijection is compatible with the parametrization of these two sets via $g$-vectors and 
shear coordinates.
This bijection has some application to cluster algebras,
a class of combinatorially defined commutative algebras discovered by
Fomin and Zelevinsky \cite{FZ}.
Initially meant as a tool to describe parts of Lusztig's dual
canonical basis of quantum groups in a combinatorial way, cluster algebras turned out to appear at
numerous different places of mathematics and mathematical
physics.
The generically $\tau$-reduced decorated components parametrize the generic Caldero-Chapoton functions,
which belong to the coefficient-free upper cluster algebra $\cU_{\SM}$ associated with $\SM$.
In many cases, these generic Caldero-Chapoton functions are known to form a basis, called the \emph{generic basis}, of $\cU_{\SM}$, see
for example \cite{GLS} and \cite{Q}.
We use the bijection mentioned above to show that the generic basis coincides with 
Musiker-Schiffler-Williams' bangle basis (see \cite[Corollary~1.3]{MSW2}) of the 
coefficient-free cluster algebra $\cA_{\SM}$ associated with $(\bbS,\bbM)$.
It is known in most cases (for example, if $|\bbM| \ge 2$) that
$\cA_{\SM} = \cU_{\SM}$, see \cite{Mu1,Mu2}.

In the following subsections we describe our results in more detail.

\subsection{Gentle algebras}
Let $Q = (Q_0,Q_1,s,t)$ be a quiver.
Thus by definition, $Q_0$ and $Q_1$ are finite sets, where the elements of $Q_0$ and $Q_1$ are the \emph{vertices} and \emph{arrows} of $Q$, respectively.
Furthermore, $s$ and $t$ are maps $s,t\df Q_1 \to Q_0$,
where $s(a)$ and $t(a)$ are the \emph{starting vertex} and
\emph{terminal vertex} of an arrow $a \in Q_1$, respectively. 
A \emph{loop} in $Q$ is an arrow $a \in Q_1$ with $s(a) = t(a)$.
 
A basic algebra $A = KQ/I$ is a \emph{gentle algebra}
provided the following hold:
\begin{itemize}

\item[(i)]
For each $i \in Q_0$ we have
$|\{ a \in Q_1 \mid s(a) = i \}| \le 2$ and
$|\{ a \in Q_1 \mid t(a) = i \}| \le 2$.

\item[(ii)]
The ideal $I$ is generated by a set $\rho$ of paths of length 2.

\item[(iii)]
Let $a,b,c \in Q_1$ such that $a \not= b$ and $t(a) = t(b) = s(c)$.
Then exactly one of the paths $ca$ and $cb$ is in $I$.

\item[(iv)]
Let $a,b,c \in Q_1$ such that $a \not= b$ and $s(a) = s(b) = t(c)$.
Then exactly one of the paths $ac$ and $bc$ is in $I$.

\end{itemize}

A gentle algebra $A = KQ/I$ is a Jacobian algebra in the sense of
\cite{DWZ1}  
if and only if the following hold:
\begin{itemize}

\item[(v)]
$Q$ is connected.

\item[(vi)]
$Q$ does not have any loops.

\item[(vii)]
Let $a,b \in Q_1$ such that $s(a) = t(b)$ and $ab \in I$.
Then there exists an arrow $c \in Q_1$ with
$s(c) = t(a)$ and $t(c) = s(b)$ such that 
$bc,ca \in I$.

\end{itemize}
The gentle Jacobian algebras are exactly the Jacobian algebras 
associated to triangulations of unpunctured marked
surfaces.
This follows from \cite[Section~2]{ABCP}.

\subsection{Smooth locus and generic reducedness of 
module schemes}
Let $Q$ be a quiver with $Q_0 = \{ 1,\ldots,n \}$, and let
$A = KQ/I$ be a basic algebra.
For $\bd \in \N^n$ let
$\irr(A,\bd)$
be the set of irreducible components
of the affine scheme $\md(A,\bd)$ of $A$-modules
with dimension vector $\bd$.
For $Z \in \irr(A,\bd)$ we write $\dimv(Z) := \bd$.
Let 
$$
\irr(A) := \bigcup_{\bd \in \N^n} \irr(A,\bd).
$$
The group 
$$
\GL_\bd(K) := \prod_{i=1}^n \GL_{d_i}(K)
$$ 
acts on 
the $K$-rational points of
$\md(A,\bd)$ by conjugation, where $\bd = (d_1,\ldots,d_n)$.
The orbit of $M \in \md(A,\bd)$ is denoted by $\cO_M$.
The orbits in $\md(A,\bd)$ correspond bijectively to the isomorphism
classes of $A$-modules with dimension vector $\bd$.

For $Z \in \irr(A,\bd)$ let
$Z^\circ$
be the \emph{interior} of $Z$.
These are all $M \in Z$ such that $M$ is not contained in any other
irreducible component of $\md(A,\bd)$.
Obviously $Z^\circ$ is a non-empty, open, irreducible subset of $\md(A,\bd)$.

A module $M \in \md(A,\bd)$ is
\emph{smooth}, if 
$$
\dim T_M = \max \{ \dim(Z) \mid Z \in \irr(A,\bd),\; M \in Z \},
$$
where $T_M$ is the tangent space of $M$ at the affine scheme
$\md(A,\bd)$.
Otherwise, $M$ is \emph{singular}.
Let $\smooth(A,\bd)$ denote the set of smooth points of $\md(A,\bd)$.

For each gentle algebra $A$ we obtain a complete description of smooth points of $\md(A,\bd)$ for all $\bd$, see
Theorem~\ref{thm:main1b}.
As a consequence we get the following
neat characterization for the case of gentle Jacobian algebras.

\begin{Thm}[Smooth points]\label{thm:main1}
Let $A$ be a gentle Jacobian algebra.
For each dimension vector $\bd$ we have
$$
\smooth(A,\bd) = \bigcup_{Z \in \irr(A,\bd)} Z^\circ.
$$ 
\end{Thm}

Note that the inclusion $\subseteq$ in Theorem~\ref{thm:main1}
is true for arbitrary basic algebras $A$.
The other inclusion $\supseteq$ is wrong in general.
For example, it fails for most gentle algebras which are not
Jacobian algebras.

A module $M \in \md(A,\bd)$ is \emph{reduced} if 
$$
\dim T_M = \dim T_M^\rd,
$$
where $T_M^\rd$ is the tangent space of $M$ at the reduced affine scheme $\md(A,\bd)^\rd$ associated with $\md(A,\bd)$.
We call $\md(A,\bd)$ \emph{reduced} if 
$\md(A,\bd) = \md(A,\bd)^\rd$.
This is the case if and only if 
$M$ is reduced for all $M \in \md(A,\bd)$.

An irreducible component $Z \in \irr(A)$ is 
\emph{generically reduced}
provided $Z$ contains a dense open subset $U$ such that each
$M \in U$ is reduced.

\begin{Thm}[Generic reducedness]\label{thm:main2}
Let $A$ be a gentle algebra without loops.
Then each $Z \in \irr(A)$ is generically reduced.
\end{Thm}

We prove a slightly more general version of Theorem~\ref{thm:main2}
where we characterize all generically reduced components
for arbitrary gentle algebras, see Theorem~\ref{thm:main2b}.

For acyclic gentle algebras, Theorem~\ref{thm:main2} is a
consequence of \cite{DS}.

\subsection{Generically $\tau$-reduced components}
For $M \in \md(A,\bd)$ let
\begin{align*}
c_A(M) &:= \max \{ \dim(Z) \mid 
Z \in \irr(A,\bd),\; M \in Z \} - \dim \cO_M,\\
e_A(M) &:= \dim \Ext_A^1(M,M),\\
h_A(M) &:= \dim \Hom_A(M,\tau_A(M)).
\end{align*}
Here $\tau_A$ denotes the Auslander-Reiten translation of $A$.

For each $Z \in \irr(A)$ there is a dense open subset $U \subseteq Z$
such that the maps $c_A$, $e_A$ and $h_A$ are constant on
$U$.
These generic values are denoted by $c_A(Z)$, $e_A(Z)$ and
$h_A(Z)$.

It follows that
\begin{align*}
c_A(Z) &= \min \{ \dim(Z) - \dim \cO_M \mid M \in Z \},\\
e_A(Z) &= \min \{ \dim \Ext_A^1(M,M) \mid M \in Z \}.
\end{align*}

Voigt's Lemma \ref{lem:Voigt} and the Auslander-Reiten formulas
(see Theorem~\ref{ARformula1}) imply that
$$
c_A(Z) \le e_A(Z) \le h_A(Z).
$$

Clearly, an irreducible component $Z$ is generically reduced if
and only if
$c_A(Z) = e_A(Z)$. 
We say that $Z$ is \emph{generically $\tau$-reduced} provided
$$
c_A(Z) = e_A(Z) = h_A(Z).
$$
Such irreducible components were first defined and studied in \cite{GLS}, where they ran under the name \emph{strongly reduced components}.

Let $\irr^\sr(A)$ be the subset of $\irr(A)$ consisting of the generically $\tau$-reduced components.

Recall that an $A$-module $M$ is \emph{rigid} (resp. $\tau$-\emph{rigid}) if
$\Ext_A^1(M,M) = 0$ (resp. $\Hom_A(M,\tau_A(M)) = 0$).
By the Auslander-Reiten formulas, any $\tau$-rigid module is
rigid, wheras the converse is wrong in general.
Each rigid $A$-module $M$ yields a generically reduced component
$Z = \overline{\cO_M}$.
If $M$ is $\tau$-rigid, then this $Z$ is generically $\tau$-reduced.

The next result says that for gentle Jacobian algebras, the 
generically $\tau$-reduced components
are determined by their dimension vectors.

\begin{Thm}\label{thm:main3}
Let $A$ be a gentle Jacobian algebra.
For $Z_1,Z_2 \in \irr^\tau(A)$ the following are
equivalent:
\begin{itemize}

\item[(i)]
$\dimv(Z_1) = \dimv(Z_2)$;

\item[(ii)]
$Z_1 = Z_2$.

\end{itemize}
\end{Thm}

Let $A = KQ/I$ be a gentle Jacobian algebra with 
$Q_0 = \{ 1,\ldots,n \}$.
Recall that the ideal $I$ is generated by a set $\rho$ of paths
of length $2$.
We denote the standard idempotents of $A$ by $e_1,\ldots,e_n$.
Let $a \in Q_1$.
Then we are in one of the following two cases:
\begin{itemize}

\item[(i)]
There is no arrow $b \in Q_1$ with $s(a) = t(b)$ such that
$ab \in I$.
In this case, the 3-dimensional subalgebra of $A$ spanned by $e_{s(a)}$, $e_{t(a)}$ and $a$
is called a \emph{$2$-block} of $A$.

\item[(ii)]
There are arrows $b,c \in Q_1$ with
$s(a) = t(b)$, $s(c) = t(a)$ and $s(b) = t(c)$
such that
$ab,ca,bc \in I$.
In this case, the 6-dimensional subalgebra of $A$ spanned by
$e_{s(a)}$, $e_{s(b)}$, $e_{s(c)}$, $a$, $b$ and $c$ is called a \emph{$3$-block}
of $A$.

\end{itemize}
In the special case where the quiver $Q$ consists just of a single vertex, we call $A$ itself a \emph{$1$-block}.
A \emph{$\rho$-block} of $A$ is a subalgebra which is either a 
$1$-block, $2$-block or $3$-block.
(Note that the $\rho$-blocks are not necessarily unital subalgebras, i.e.
the unit of a $\rho$-block of $A$ does in
general not coincide with the unit of $A$.)


We say that 
a vertex $j \in Q_0$ or an arrow $a \in Q_1$ \emph{belongs to a $\rho$-block} $A_i$ of $A$ if $e_j \in A_i$ or $a \in A_i$, respectively.
Note that each arrow of $Q$ belongs to exactly one $\rho$-block of $A$,
and each vertex of $Q$ belongs to at most two $\rho$-blocks.

The restriction of representations of a gentle Jacobian
algebra $A$ to its $\rho$-blocks $A_1,\ldots,A_t$ yields a bijection
\begin{align*}
\irr(A) &\to \irr(A_1) \times \cdots \times \irr(A_t)
\\
Z &\mapsto (\pi_1(Z),\ldots,\pi_t(Z)).
\end{align*}
In Section~\ref{sec:blocks} we extend this observation to arbitrary basic algebras $A = KQ/I$.
This reduces the study of schemes of modules over gentle algebras to schemes of complexes.

Our next result characterizes the generically $\tau$-reduced
components of a gentle Jacobian algebra in terms of
the generically $\tau$-reduced components of its $\rho$-blocks.

The fact that
the generic reducedness or the smooth locus of a component $Z$ relate
to the generic reducedness or the smooth locus of the components $\pi_i(Z)$
does not come as a surprise.
The following result however is somewhat unexpected, since the
Auslander-Reiten translation for $A$ is quite different from the
Auslander-Reiten translations for the $\rho$-blocks of $A$.

\begin{Thm}\label{thm:main4}
Let $A = KQ/I$ be a gentle Jacobian algebra, and
let $A_1,\ldots,A_t$ be its $\rho$-blocks.
For an irreducible component $Z \in \irr(A)$ the following
are equivalent:
\begin{itemize}

\item[(i)]
$Z \in \irr^\tau(A)$;

\item[(ii)]
$\pi_i(Z) \in \irr^\tau(A_i)$ for all $1 \le i \le t$.

\end{itemize}
\end{Thm}

One might ask if Theorem~\ref{thm:main4} holds for arbitrary 
finite-dimensional $K$-algebras using of course a generalized
definition for \emph{$\rho$-blocks}.

\subsection{Band components}
The indecomposable modules over a gentle algebra $A$
(or more generally, over a string algebra)
are either \emph{string modules} or \emph{band modules},
see \cite{BR,WW} for details.
The band modules occur naturally in $K^*$-parameter families.
An irreducible component $Z \in \irr(A)$ is a \emph{string component} if
it contains a string module whose orbit is dense in $Z$, and $Z$
is a \emph{band component} if it contains a $K^*$-parameter
family of band modules whose union of orbits is dense in $Z$.

An irreducible component $Z \in \irr(A)$ is a \emph{brick component}
if it contains a \emph{brick}, i.e. an $A$-module $M$ with
$\dim \End_A(M) = 1$.
In this case, by upper semicontinuity the bricks in $Z$ form a dense open subset of $Z$.

\begin{Thm}\label{thm:main6}
Let $A$ be a gentle algebra.
Then each band component is a brick component.
\end{Thm}

Using the terminology of \cite{CBS}, each irreducible component
$Z \in \irr(A)$ is a direct sum of uniquely determined indecomposable
irreducible components.
The string and band components are the only indecomposable
components for string algebras.

The generically $\tau$-reduced string components are exactly the components
containing an indecomposable $\tau$-rigid module, which 
is then automatically a string module.

\begin{Thm}\label{thm:main5}
Let $A$ be a gentle algebra.
For $Z \in \irr(A,\bd)$ the following are equivalent:
\begin{itemize}

\item[(i)]
$Z$ is a direct sum of band components. 

\item[(ii)]
$\dim(Z) = \dim(\GL_\bd(K))$.

\end{itemize}
In this case, $Z$ is generically $\tau$-reduced.
\end{Thm}

Theorem~\ref{thm:main5} is closely related to the seemingly
different
\cite[Proposition~4.3]{CCKW}.
The proofs follow the same line of arguments.
We thank Ryan Kinser for pointing this out to us.

For acyclic gentle algebras,
Theorems~\ref{thm:main6} and
\ref{thm:main5} can be extracted from Carroll and
Chindris
\cite[Corollary~10]{CC} and \cite[Proposition~11]{CC}, see
also \cite[Theorem~2]{C}.
As a consequence of Theorem~\ref{thm:main6}, one gets the
known result that
a gentle algebra $A$ is representation-finite if and only if 
$\md(A)$ contains just finitely many bricks, compare 
\cite[Theorem~1.1]{P2}.

\subsection{Laminations of marked surfaces and generically $\tau$-reduced components}
A \emph{lamination} of an unpunctured marked surface $(\bbS,\bbM)$ is a set of homotopy classes of curves and loops in $(\bbS,\bbM)$, which do not intersect each other, together
with a positive integer attached to each class.
Let $\Lam\SM$ be the set of such laminations.
(For more precise definitions, we refer to Section~\ref{sec:laminations}.)

Let $T$ be a triangulation of $(\bbS,\bbM)$, and let
$A_T$ be the associated gentle Jacobian algebra.
A \emph{decorated irreducible component} is roughly speaking 
an irreducible component of $\md(A_T,\bd)$ equipped with a
certain integer datum.
Similarly as before, one defines generically $\tau$-reduced decorated
irreducible components.
Let $\decirr^\tau(A_T)$ be the set of all generically $\tau$-reduced decorated components of $\decrep(A_T,(\bd,\bv))$, where $(\bd,\bv)$ runs over all dimension vectors. 
A precise definition can be found in Section~\ref{sec:decrep}.

\begin{Thm}\label{thm:main7}
Let $(\bbS,\bbM)$ be an unpunctured marked surface, and
let $T$ be a triangulation of $(\bbS,\bbM)$.
Let $A = A_T$ be the associated Jacobian algebra.
Then there is a natural bijection
$$
\eta_T\df \Lam\SM \to \decirr^\tau(A).
$$
\end{Thm}

In their ground breaking work, 
Fomin, Shapiro and Thurston \cite{FST} proved that 
the laminations of $\SM$ consisting of curves 
are in bijection with
the cluster monomials of a cluster algebra $\cA_{\SM}$ 
associated with $\SM$.
Note that Fomin, Shapiro and Thurston work with cluster algebras with arbitrary coefficient
systems, whereas we always assume that
$\cA_{\SM}$ is a coefficient-free 
cluster algebra.

Musiker, Schiffler and Williams \cite{MSW2}  
defined a set 
$$
\cB_T := \{ \psi_L \mid L \in \Lam\SM \} 
$$
of \emph{bangle functions}, whose elements are
parametrized by $\Lam\SM$, and which 
(by results in \cite{MSW1}) contains all cluster monomials.
They show that $\cB_T$ forms a basis of $\cA_{\SM}$ provided
$|\bbM| \ge 2$, see \cite[Corollary~1.3]{MSW2}.

A result by W. Thurston (see \cite[Theorem~12.3]{FT})
says that there is a bijection 
$$
\bs_T\df \Lam\SM \to \Z^n
$$ 
sending a lamination to its shear coordinate.
Combining a theorem by Br\"ustle and Zhang \cite[Theorem~1.6]{BZ}
with a result by Adachi, Iyama and
Reiten \cite[Theorem~4.1]{AIR}, one gets a bijection between the laminations in
$\Lam\SM$ which consist only of curves, and the
set of generically $\tau$-reduced decorated components 
in $\decirr^\tau(A_T)$, which contain a dense orbit.
On the other hand, Plamondon \cite{P1} proved that there is a bijection
$$
\bg_T\df \decirr^\tau(A_T) \to \Z^n
$$ 
sending a component to its
$g$-vector.
Theorem~\ref{thm:main7} extends Br\"ustle-Zhang's bijection mentioned above
to a bijection 
$$
\eta_T\df \Lam\SM \to \decirr^\tau(A_T)
$$ 
such that
$\bg_T \circ \eta_T = \bs_T$.

Let 
$$
\cG_T := \{ \phi_Z \mid Z \in \decirr^\tau(A_T) \}
$$ 
be the set of \emph{generic Caldero-Chapoton functions} as defined in \cite{GLS}.
As a consequence of more general results in \cite{DWZ2}, the set $\cG_T$ is contained in the upper cluster algebra $\cU_{\SM}$ and contains all cluster monomials.
Furthermore,
by \cite[Theorem~1.3]{P1}, the set $\cG_T$ is (in a certain sense) independent of the choice of the triangulation $T$ of $\SM$.

The proof of the next theorem is based on the bijection from
Theorem~\ref{thm:main7}.

\begin{Thm}\label{thm:main8}
$\cB_T = \cG_T$.
\end{Thm}

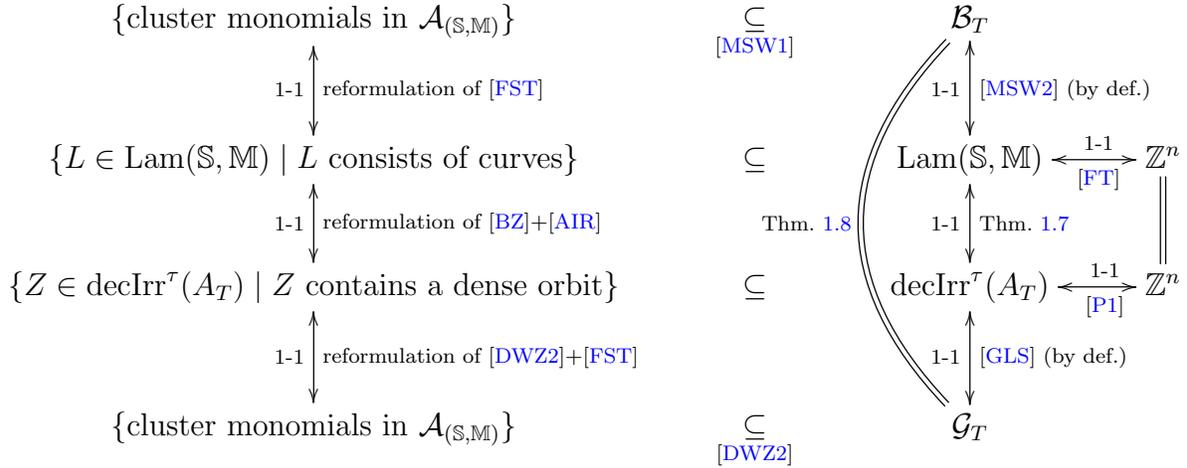
\begin{figure}[H]
$$
\xymatrix@+1ex{
\{ \text{cluster monomials in } \cA_{\SM} \}
\ar@{<->}[d]^{\text{reformulation of \cite{FST}}}_{\text{1-1}}
&\underset{\text{\cite{MSW1}}}{\subseteq}
& \cB_T  
\ar@{<->}[d]^{\text{\cite{MSW2} (by def.)}}_{\text{1-1}}
\ar@/_8ex/@{=}[ddd]_{\text{Thm.~\ref{thm:main8}}}
\\
\{ L \in \Lam\SM \mid L \text{ consists of curves} \} 
\ar@{<->}[d]^{\text{reformulation of \cite{BZ}+\cite{AIR}}}_{\text{1-1}}
&\subseteq
&\Lam\SM 
\ar@{<->}[d]^{\text{Thm.~\ref{thm:main7}}}_{\text{1-1}}
\ar@{<->}[r]_<<<<<<{\text{\cite{FT}}}^<<<<<<{\text{1-1}}
& \Z^n \ar@{=}[d]
\\
\{ Z \in \decirr^\tau(A_T) \mid Z \text{ contains a dense orbit} \} 
&\subseteq
&\decirr^\tau(A_T)
\ar@{<->}[d]^{\text{\cite{GLS} (by def.)}}_{\text{1-1}}
\ar@{<->}[r]_<<<<<<{\text{\cite{P1}}}^<<<<<<{\text{1-1}}
& \Z^n
\\
\{ \text{cluster monomials in } \cA_{\SM} \}
\ar@{<->}[u]_{\text{reformulation of \cite{DWZ2}+\cite{FST}}}^{\text{1-1}}
&\underset{\text{\cite{DWZ2}}}{\subseteq} & \cG_T 
}
$$
\caption{Bangle functions $\cB_T$ and generic Caldero-Chapoton functions $\cG_T$ for the coefficient-free cluster algebra $\cA_{\SM}$ associated with an unpunctured marked surface $\SM$.}\label{fig:intro}
\end{figure}

The diagram in Figure~\ref{fig:intro} summarizes the situation.

\subsection{Overall structure of the article}
The article is organized as follows.
After the introduction (Section~\ref{sec:intro}), we recall
in Section~\ref{sec:moduleschemes}
some fundamentals on schemes of modules over basic algebras.
Section~\ref{sec:tame} contains a characterization of 
generically $\tau$-reduced components for tame algebras.
In Section~\ref{sec:blocks} we introduce $\rho$-block
decompositions of schemes of modules and derive
some consequences on tangent spaces.
Section~\ref{sec:modules} contains a few facts on the
representation theory of gentle algebras.
We also recall the definition of rank functions of 
modules over gentle algebras.
Section~\ref{sec:complex} consists of a detailed study of 
schemes of complexes.
We determine their smooth points, and we describe all
rigid and $\tau$-rigid modules over the associated basic algebras.
In Section~\ref{sec:schemesgentle} we apply the results obtained
in Section~\ref{sec:complex} and prove 
Theorems~\ref{thm:main1}, \ref{thm:main2}, \ref{thm:main5} and \ref{thm:main6}.
The proofs of
Theorems~\ref{thm:main3} and \ref{thm:main4} can be found in
Section~\ref{sec:taureduced}.
In Section~\ref{sec:decrep} we recall some basics
on decorated modules and schemes of decorated modules
over finite-dimensional algebras.
Section~\ref{sec:laminations} contains the proof of Theorem~\ref{thm:main7}, and also the proof that
under the bijection in Theorem~\ref{thm:main7},
shear coordinates and $g$-vectors are compatible.
Theorem~\ref{thm:main8} is proved in Section~\ref{sec:bangle}. 
In Section~\ref{sec:example} we illustrate the combinatorics used in Section~\ref{sec:bangle} by an example.


\section{Scheme of modules}\label{sec:moduleschemes}


In this section, we recall some definitions and elementary
facts on the representation theory of 
basic algebras and on schemes of modules over such algebras.
Throughout, let $K$ be an algebraically closed field.

\subsection{Orbits, tangent spaces and Voigt's Lemma} 
Let $Q = (Q_0,Q_1,s,t)$ be a quiver.
If not mentioned otherwise, we always assume that 
$Q_0 = \{ 1,\ldots,n \}$.

A \emph{path} in $Q$ 
is a tuple $p = (a_1,\ldots,a_m)$ of
arrows $a_i \in Q_1$ such that $s(a_i) = t(a_{i+1})$ for all
$1 \le i \le m-1$.
Then $\length(p) := m$ is the \emph{length} of $p$, and we set
$s(p) := s(a_m)$ and $t(p) := t(a_1)$.
Additionally, for each vertex $i \in Q_0$ there is a path
$e_i$ of length $0$, and let $s(e_i) = t(e_i) = i$.
We often just write $a_1 \cdots a_m$ instead of $(a_1,\ldots,a_m)$.
A path $p = (a_1,\ldots,a_m)$ of length $m \ge 1$ is a
\emph{cycle} in $Q$, or more precisely an
$m$-\emph{cycle} in $Q$, if $s(p) = t(p)$.

Let $KQ$ be the path algebra of $Q$, and let $\m$ be the ideal
generated by the arrows of $Q$.
An ideal $I$ of $KQ$ is \emph{admissible} if there exists some
$m \ge 2$ such that
$\m^m \subseteq I \subseteq \m^2$.
In this case, we call $A := KQ/I$ a \emph{basic algebra}.
Clearly, basic algebras are finite-dimensional.
By a Theorem of Gabriel, each finite-dimensional $K$-algebra is Morita equivalent to a basic algebra.

A \emph{relation} in $KQ$ is a linear combination 
$$
\sum_{i=1}^s \lambda_i p_i
$$
where the $p_i$ are pairwise different paths of length at
least $2$ in $Q$ with $s(p_i) = s(p_j)$ and $t(p_i) = t(p_j)$ 
for all $1 \le i,j \le s$ and $\lambda_i \in K^*$ for all $i$.

Each admissible ideal is generated by a finite set of relations.

Let $A = KQ/I$ be a basic algebra.
Up to isomorphism, there are $n$ simple $A$-modules
$S_1,\ldots,S_n$ corresponding to the vertices
of $Q$.
Let $P_1,\ldots,P_n$ (resp. $I_1,\ldots,I_n$) be the projective
covers (resp. injective envelopes) of the simple modules $S_1,\ldots,S_n$.

A \emph{representation} of a quiver $Q = (Q_0,Q_1,s,t)$ is a 
tuple $M = (M_i,M_a)_{i \in Q_0,a \in Q_1}$, where  
$M_i$ is
a finite-dimensional $K$-vector space for each $i \in Q_0$, and $M_a\colon M_{s(a)} \to M_{t(a)}$ is a $K$-linear map for 
each arrow $a \in Q_1$.

For a path $p = (a_1,\ldots,a_m)$ in $Q$ and a representation 
$M$ as above, let 
$$
M_p := M_{a_1} \circ \cdots \circ M_{a_m}.
$$

We call 
$$
\dimv(M) := (\dim(M_1),\ldots,\dim(M_n))
$$ 
the
\emph{dimension vector} of $M$, and let
$$
\dim(M) := \dim(M_1) + \cdots + \dim(M_n)
$$ 
be the
\emph{dimension} of $M$. 
The $i$th entry
$\dim(M_i)$ of $\dimv(M)$ 
equals the
Jordan-H\"older multiplicity $[M:S_i]$ of $S_i$ in $M$.

A \emph{representation} of a basic algebra $A = KQ/I$ is a representation $M$ of $Q$, which is annihilated by the ideal $I$, i.e.
for each relation 
$$
\sum_{j=1}^s \lambda_j p_j 
$$
in $I$ we demand that
$$
\sum_{j=1}^s \lambda_j M(p_j) = 0.
$$

In the usual way,
we identify the category $\rep(A)$ of representations of $A$ 
with the category $\md(A)$ of finite-dimensional left $A$-modules.

For $\bd = (d_1,\ldots,d_n) \in \N^n$ let
$\md(A,\bd)$ be the affine scheme of representations
of $A$
with dimension vector $\bd$.
By definition the $K$-rational points of $\md(A,\bd)$ are the representations $M = (M_i,M_a)_{i \in Q_0,a \in Q_1}$ of $A$
with $M_i = K^{d_i}$ for all $i \in Q_0$.
When there is no danger of confusion, we will just write $\md(A,\bd)$
for the set of $K$-rational points of $\md(A,\bd)$.
One can regard $\md(A,\bd)$ as
a Zariski closed subset of the affine space
$$
\md(Q,\bd) := \prod_{a \in Q_1} 
\Hom_K(K^{d_{s(a)}},K^{d_{t(a)}}).
$$
The group $\GL_\bd(K)$ acts on 
the $K$-rational points of
$\md(A,\bd)$ by conjugation.
More precisely, 
for $g = (g_1,\ldots,g_n) \in \GL_\bd(K)$ 
and $M \in \md(A,\bd)$ let
$$
g.M := (M_i,g_{t(a)}^{-1}M_ag_{s(a)})_{i \in Q_0,a \in Q_1}.
$$
For $M \in \md(A,\bd)$ let $\cO_M$ be the $\GL_\bd(K)$-orbit
of $M$. 
The $\GL_\bd(K)$-orbits are in bijection with the isomorphism
classes of representations of $A$ with dimension vector $\bd$.

For $M \in \md(A,\bd)$ we denote the tangent
space of $M$ at $\md(A,\bd)$ by
$T_M$.
Let $T_M(\cO_M)$ be the tangent space of
$M$ at $\cO_M$.
Since the $\GL_\bd(K)$-orbit $\cO_M$ is smooth, we
have 
$$
\dim T_M(\cO_M) = \dim \cO_M = \dim \GL_\bd(K) - \dim \End_A(M).
$$

The following lemma is obvious.

\begin{Lem}\label{lem:schemes1}
For $M \in \md(A,\bd)$ the following are equivalent:
\begin{itemize}

\item[(i)]
$\cO_M$ is open.

\item[(ii)]
The Zariski closure $\overline{\cO_M}$ is an irreducible component
of $\md(A,\bd)$.

\end{itemize}
\end{Lem}

For the following proposition we refer to Gabriel 
\cite[Proposition~1.1]{Ga} and Voigt \cite{V}.

\begin{Lem}[Voigt's Lemma]\label{lem:Voigt}
For $M \in \md(A,\bd)$ there is an isomorphism
$$
T_M/T_M(\cO_M) \to \Ext_A^1(M,M).
$$
of $K$-vector spaces.
\end{Lem}

\begin{Cor}\label{cor:schemes2}
Let $M \in \md(A,\bd)$ be rigid.
Then $\cO_M$ is open.
\end{Cor}

The converse of Corollary~\ref{cor:schemes2} is in general wrong.

\begin{Cor}\label{cor:rigidsmooth}
Let $M \in \md(A,\bd)$ be rigid.
Then $M$ is smooth.
\end{Cor}

\begin{Cor}\label{cor:schemes3}
For $M \in \md(A,\bd)$ the following are equivalent:
\begin{itemize}

\item[(i)]
$M$ is rigid.

\item[(ii)]
The Zariski closure of $\cO_M$ is a generically reduced component of $\md(A,\bd)$.

\end{itemize}
\end{Cor}

\begin{Lem}\label{lem:smoothreduced}
Let $M \in \md(A,\bd)$ be smooth.
Then $M$ is reduced.
\end{Lem}

\begin{proof}
This is clear, since for each irreducible component
$Z$ with $M \in Z$ we have
$$
\dim(Z) \le \dim T_M^\rd \le \dim T_M.
$$
\end{proof}

The following three results are well known and can be 
extracted e.g. from \cite{H,Sh1,Sh2}.

\begin{Prop}\label{prop:reducedtangent}
Let $Z \in \irr(A,\bd)$.
Then there is a dense open subset $U \subseteq Z$ such that
$$
\dim T_M^\rd = \dim(Z) 
$$
for all $M \in U$.
\end{Prop}

\begin{Prop}\label{prop:singularclosed}
Let $Z \in \irr(A,\bd)$.
Then the smooth points in $Z$ form a (possibly empty) 
open subset of $Z$.
\end{Prop}

\begin{Prop}\label{prop:intersecsingular}
Let $M \in \md(A,\bd)$ be contained in at least two different
irreducible components.
Then $M$ is singular.
\end{Prop}

The following statement is proved in \cite[Proposition~3.7]{G}.
It relies on results from \cite{GP}.

\begin{Prop}\label{prop:rigidsmooth2}
Let $M \in \md(A,\bd)$ with $\Ext_A^2(M,M) = 0$.
Then $M$ is smooth.
\end{Prop}

\subsection{Canonical decompositions of irreducible
components}
An irreducible component $Z \in \irr(A,\bd)$ is \emph{indecomposable} if 
$$
\ind(Z) := \{ M \in Z \mid M \text{ is indecomposable} \}
$$
is dense in $Z$.
Let $\bd$ and $\bd_1,\ldots,\bd_t$ be dimension vectors
with $\bd = \bd_1 + \cdots + \bd_t$.
For $Z_i \in \irr(A,\bd_i)$ with $1 \le i \le t$ let
$$
Z_1 \oplus \cdots \oplus Z_t
$$
be the image of the morphism
\begin{align*}
\GL_\bd(K) \times Z_1 \times \cdots \times Z_t &\to \md(A,\bd)
\\
(g,M_1,\ldots,M_t) &\mapsto
g.(M_1 \oplus \cdots \oplus M_t).
\end{align*}
For each $Z \in \irr(A)$ there are uniquely determined (up to renumbering) indecomposable irreducible components 
$Z_1,\ldots,Z_t \in \irr(A)$ such that
$$
Z = \overline{Z_1 \oplus \cdots \oplus Z_t},
$$
see \cite[Theorem~1.1]{CBS}.
This is called the \emph{canonical decomposition} of $Z$.
For $Z \in \irr(A,\bd)$ set $\dimv(Z) := \bd$.
For $Z_1,Z_2 \in \irr(A)$ let
$$
\ext_A^1(Z_1,Z_2) := \min \{ \dim \Ext_A^1(M_1,M_2) \mid M_1\in Z_1,\; M_2 \in Z_2 \}.
$$

\begin{Thm}[{\cite[Theorem~1.2]{CBS}}]\label{thm:CBS}
Let $A$ be a finite-dimensional $K$-algebra.
For $Z_1,\ldots,Z_t \in \irr(A)$ the following are equivalent:
\begin{itemize}

\item[(i)]
$\overline{Z_1 \oplus \cdots \oplus Z_t} \in \irr(A)$;

\item[(ii)]
$\ext_A^1(Z_i,Z_j) = 0$ for all $i \not= j$.

\end{itemize}
\end{Thm}

For each $Z \in \irr^\tau(A)$ there are uniquely determined (up to renumbering) indecomposable components 
$Z_1,\ldots,Z_t \in \irr^\tau(A)$ such that
$$
Z = \overline{Z_1 \oplus \cdots \oplus Z_t}.
$$

\begin{Thm}[{\cite[Theorem~5.11]{CLFS}}]\label{thm:CLFS}
Let $A$ be a finite-dimensional $K$-algebra.
For $Z_1,\ldots,Z_t \in \irr^\tau(A)$ the following are equivalent:
\begin{itemize}

\item[(i)]
$\overline{Z_1 \oplus \cdots \oplus Z_t} \in \irr^\tau(A)$;

\item[(ii)]
$h_A(Z_i,Z_j) = 0$ for all $i \not= j$.

\end{itemize}
\end{Thm}


\section{Generically $\tau$-reduced components for
tame algebras}\label{sec:tame}


In this section, we characterize the indecomposable
$\tau$-reduced components for tame algebras.
The proof consists basically of combining some known results
in a straightforward manner.

Let $A$ be a finite-dimensional $K$-algebra.
Then $A$ is a \emph{tame algebra} if for each dimension
$d$ there exists a finite number $M_1,\ldots,M_t$ of
$A$-$K[X]$-bimodules $M_i$, which are free of rank 
$d$ as $K[X]$-modules, such that
all but finitely many $d$-dimensional $A$-modules are 
isomorphic to 
$$
M_i \otimes_{K[X]} K[X]/(X-\lambda) 
$$
for some $1 \le i \le t$ and some $\lambda \in K$.

The following lemma is well known folklore.
A proof  can be found in \cite[Section~2.2]{CC}.

\begin{Lem}\label{lem:CCtame}
Let $A$ be a tame algebra, and let $Z \in \irr(A)$ be an
indecomposable irreducible component.
Then $c_A(Z) \in \{ 0,1 \}$.
Furthermore, the following hold:
\begin{itemize}

\item[(i)]
$c_A(Z) = 0$ if and only if $Z$ contains an indecomposable module
$M$ with 
$$
Z = \overline{\cO_M}.
$$

\item[(ii)]
$c_A(Z) = 1$ if and only if $Z$ contains a rational curve $C$
such that the points of $C$ are pairwise non-isomorphic indecomposable modules with
$$
Z = \overline{\bigcup_{M \in C} \cO_M}.
$$ 

\end{itemize}
\end{Lem}

\begin{Thm}\label{thm:tame}
Let $A$ be a tame algebra, and let $Z \in \irr(A)$ be an
indecomposable irreducible component.
Then the following hold:
\begin{itemize}

\item[(i)]
For $c_A(Z) = 0$ the following are equivalent:
\begin{itemize}

\item[(a)]
$Z$ is generically $\tau$-reduced.

\item[(b)]
$Z$ contains an indecomposable $\tau$-rigid
module $M$.

\end{itemize}
In this case, 
$$
Z = \overline{\cO_M}.
$$

\item[(ii)]
For $c_A(Z) = 1$ the following are equivalent:
\begin{itemize}

\item[(a)]
$Z$ is generically $\tau$-reduced.

\item[(b)]
$Z$ contains a rational curve $C$
such that the points of $C$ are pairwise non-isomorphic 
bricks.

\item[(c)]
$Z$ contains infinitely many pairwise non-isomorphic bricks.

\end{itemize}
In this case,
$$
Z = \overline{\bigcup_{M \in C} \cO_M}. 
$$ 

\end{itemize}
\end{Thm}

\begin{proof}
Part (i) follows directly from the definitions.
Thus assume $c_A(Z) = 1$.
By Lemma~\ref{lem:CCtame},
$Z$ contains a rational curve $C$
such that the points of $C$ are pairwise non-isomorphic indecomposable modules with
$$
Z = \overline{\bigcup_{M \in C} \cO_M}.
$$ 
Now a deep result by Crawley-Boevey
\cite[Theorem~D]{CB2} says that
$\tau_A(M) \cong M$ for all but finitely many $M \in C$.
Thus
$Z$ is generically $\tau$-reduced if and only if 
$h_A(Z) = 1$ if and only if $\dim \Hom_A(M,\tau_A(M)) = 
\dim \End_A(M) = 1$ for all but finitely many $M \in C$.
Thus (a) and (b) are equivalent.
In a brick component, the bricks always form a dense open subset.
Keeping in mind Lemma~\ref{lem:CCtame}, this implies the
equivalence of (b) and (c).
\end{proof}

For an arbitrary finite-dimensional $K$-algebra $A$,
each generically $\tau$-reduced component
$Z \in \irr(A)$ is a direct sum of indecomposable
generically $\tau$-reduced components.
This is explained in Section~\ref{subsec:decomp}.


\section{$\rho$-block decomposition and tangent spaces}\label{sec:blocks}


Let $A = KQ/I$, where $KQ$ is a path algebra and 
$I$ is an admissible ideal generated by a set $\rho = \{ \rho_1,\ldots,\rho_m \}$ of
relations.

For each
$$
\rho_k = \sum_{i=1}^s \lambda_i p_i
$$
with $1 \le k \le m$
let $Q(\rho_k)$ be the smallest subquiver of $Q$ containing the
paths $p_i$.
Of course, these subquivers might overlap for different relations.

For arrows $a,b \in Q_1$ write $a \sim b$ if there is some $k$ with
$a,b \in Q(\rho_k)$.
Let $\sim$ be the smallest equivalence relation on $Q_1$ respecting this rule.
In particular,
each $a \in Q_1$ which is not contained in any of the $Q(\rho_k)$
forms its own equivalence class.
 
Each equivalence class in $Q_1$ with respect to $\sim$ gives rise to a subquiver of $Q$ and also to a subalgebra of $A$.
These subalgebras are the \emph{$\rho$-blocks} of $A$.
Each vertex $i \in Q_0$, which has no arrow attached to it
yields a $1$-dimensional subalgebra (with basis $e_i$).
Such subalgebras are also called \emph{$\rho$-blocks} of $A$.

Not under this name and for a different purpose (tameness proofs), $\rho$-blocks
appear already in \cite{B}, see also \cite{AB}.
We thank Thomas Br\"ustle for pointing this out. 

Let us remark that each arrow of $Q$ belongs to exactly one $\rho$-block, whereas a standard idempotent $e_i$ can belong to several
$\rho$-blocks.
For
an arrow $a$ which does not appear in any of the relations in $\rho$,
the path algebra of the quiver
$$
\xymatrix{
s(a) \ar[r]^a & t(a)
}
$$
is a $\rho$-block.
For example, let $Q$ be the quiver
$$
\xymatrix{
1 \ar[r]^{a_1} & 2 \ar[r] & \cdots \ar[r]^{a_{n-1}} & n
}
$$
and let $A = KQ$. (In this trivial example, we have $I = 0$ and
$\rho = \varnothing$.)
For $n \ge 2$
the $\rho$-blocks of $A$ are the path algebras of the subquivers
$$
\xymatrix{
i \ar[r]^<<<<{a_i} & {i+1}
}
$$
where $1 \le i \le n-1$.
For $n=1$ there is only one $\rho$-block, namely $A$ itself.

As another example,
let $Q$ be the quiver 
$$
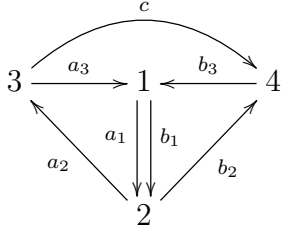
\xymatrix@+1pc{
3\ar@/^2pc/[rr]^c \ar[r]^{a_3} & 1 \ar@<-0.5ex>[d]_(.4){a_1} \ar@<0.5ex>[d]^(.4){b_1}
&4 \ar[l]_{b_3}\\
& 2 \ar[ul]^{a_2}\ar[ur]_{b_2}
}
$$
and let $I$ be the ideal in $KQ$ generated by
$\rho = \{ a_1a_3,a_2a_1,a_3a_2,b_1b_3,b_2b_1,b_3b_2 \}$. 
Then $KQ/I$ is a gentle Jacobian algebra, and
there are two $\rho$-blocks with three vertices and one 
$\rho$-block with two vertices.
(This algebra arises from a torus with one boundary component
and one marked point on the boundary.)

Our $\rho$-blocks are in general very different from the classically defined blocks of an algebra.
However, on the geometric level there is at least some resemblance.
This will be explained at the end of this subsection.

Now let $A_1,\ldots,A_t$ be the $\rho$-blocks of $A$.
For each dimension vector $\bd \in \N^n$ and
$1 \le i \le t$ let $\pi_i(\bd)$ denote the corresponding dimension vector
for $A_i$.
Each $M \in \md(A,\bd)$ induces via restriction modules
$\pi_i(M) \in \md(A_i,\pi_i(\bd))$ for $1 \le i \le t$ in the obvious way.

For each $\bd$ we obtain an isomorphism 
\begin{align*}
\md(A,\bd) &\to \md(A_1,\pi_1(\bd)) \times \cdots \times \md(A_t,\pi_t(\bd))
\\
M &\mapsto (\pi_1(M),\ldots,\pi_t(M))
\end{align*}
of affine schemes and therefore a bijection
\begin{align*}
\irr(A,\bd) &\to 
\irr(A_1,\pi_1(\bd)) \times \cdots \times \irr(A_t,\pi_t(\bd))
\\
Z &\mapsto (\pi_1(Z),\ldots,\pi_t(Z)).
\end{align*}

\begin{Prop}\label{prop:blocks1}
Let $A$ and $A_1,\ldots,A_t$ be defined as above.
For $M \in \md(A,\bd)$
the following hold:
\begin{itemize}

\item[(i)]
$T_M \cong \prod_{i=1}^n T_{\pi_i(M)}$;

\item[(ii)]
$T_M^\rd \cong \prod_{i=1}^n T_{\pi_i(M)}^\rd$.

\end{itemize}
\end{Prop}

\begin{proof}
(i): 
Obvious.

(ii): 
For a ring $R$ let $\nil(R)$ be its ideal of nilpotent elements.
For $R$ commutative and finitely generated, let $\Spec(R)$ be
as usual its prime ideal spectrum, which is an affine scheme. 

We have an isomorphism of affine schemes
$$
\md(A,\bd) \cong \prod_{i=1}^t \md(A_i,\pi_i(\bd)).
$$
Let $R_i$ be the coordinate algebra of 
$\md(A_i,\pi_i(\bd))$ for $1 \le i \le t$.
We get an isomorphism of affine schemes
$$
\md(A,\bd) \cong \Spec(R_1 \otimes \cdots \otimes R_t).
$$
Furthermore, we have
$$
\md(A,\bd)^\rd \cong 
\Spec((R_1 \otimes \cdots \otimes R_t)/\nil(R_1 \otimes \cdots \otimes R_t)).
$$

Let $B$ and $C$ be finitely generated commutative
$K$-algebras.
Then one easily shows that
$$
\nil(B \otimes C) = \nil(B) \otimes C + B \otimes \nil(C).
$$
This yields
$$
(B \otimes C)/\nil(B \otimes C) \cong 
B/\nil(B) \otimes C/\nil(C).
$$
Applying this via induction to the situation above, we get
$$
(R_1 \otimes \cdots \otimes R_t)/\nil(R_1 \otimes \cdots \otimes R_t)
\cong
R_1/\nil(A_1) \otimes \cdots \otimes R_t/\nil(A_t).
$$
We get
$$
\md(A,\bd)^\rd \cong \prod_{i=1}^t \md(A_i,\pi_i(\bd))^\rd,
$$
which implies (ii).
\end{proof}

Proposition~\ref{prop:blocks1} allows us to study the tangent spaces
of $\md(A,\bd)$ in terms of the often easier to compute tangent
spaces of $\md(A_i,\pi_i(\bd))$.

\begin{Cor}\label{cor:blocks2}
Let $M \in \md(A,\bd)$.
Then the following are equivalent:
\begin{itemize}

\item[(i)]
$M$ is smooth.

\item[(ii)]
$\pi_i(M)$ is smooth for all $1 \le i \le t$.

\end{itemize}
\end{Cor}

\begin{Cor}\label{cor:blocks3}
Let $M \in \md(A,\bd)$.
Then the following are equivalent:
\begin{itemize}

\item[(i)]
$M$ is reduced.

\item[(ii)]
$\pi_i(M)$ is reduced for all $1 \le i \le t$.

\end{itemize}
\end{Cor}

\begin{Cor}\label{cor:blocks4}
Let $Z \in \irr(A)$.
Then the following are equivalent:
\begin{itemize}

\item[(i)]
$Z$ is generically reduced.

\item[(ii)]
$\pi_i(Z)$ is generically reduced for all $1 \le i \le t$.

\end{itemize}
\end{Cor}

For the basic algebra $A = KQ/I$, let $Q(1),\ldots,Q(m)$ be the connected components of the quiver
$Q$. For $1 \le i \le m$ let $I(i) := I \cap KQ(i)$.
Then $I(i)$ is generated by a subset $\rho(i)$ of $\rho$.
With $B_i := KQ(i)/I(i)$ we get an algebra isomorphism
$$
A \cong B_1 \times \cdots \times B_m.
$$
The $B_i$ are indecomposable algebras, i.e. they are not isomorphic to the product of two algebras of smaller dimension.
In other words, the $B_i$ are the classical \emph{blocks} of $A$.
For a dimension vector ${\bf d} \in \N^n$ let $\bd(i)$ be the corresponding dimension vector for $B_i$.
We get an isomorphism
$$
\md(A,\bd) \to \md(B_1,\bd(1)) \times \cdots \times \md(B_m,\bd(m))
$$
of affine schemes.
The $\rho$-blocks of $A$ are the disjoint union of the
$\rho(i)$-blocks of the $B_i$.


\section{Modules over gentle algebras}\label{sec:modules}


Throughout this section, let $A = KQ/I$ be a gentle algebra
with $Q = (Q_0,Q_1,s,t)$.

\subsection{The maps $\sigma$ and $\varepsilon$}
We need two maps
$$
\sigma,\varepsilon\df Q_1 \to \{ \pm 1 \}
$$
satisfying the following properties:
\begin{itemize}

\item[(i)]
If $a_1,a_2 \in Q_1$ with $a_1 \not= a_2$ and $s(a_1) = s(a_2)$,
then $\sigma(a_1) = - \sigma(a_2)$.

\item[(ii)]
If $b_1,b_2 \in Q_1$ with $b_1 \not= b_2$ and $t(b_1) = t(b_2)$,
then $\varepsilon(b_1) = - \varepsilon(b_2)$.

\item[(iii)]
If $a,b \in Q_1$ with $s(a) = t(b)$ and $ab \notin I$, then
$\sigma(b) = -\varepsilon(\gamma)$.

\end{itemize}
It is straightforward to see that such maps $\sigma$ and $\varepsilon$ exist.
We fix $\sigma$ and $\varepsilon$ for the rest of this section.

Later on we will define $1$-sided and $2$-sided standard homomorphisms.
To make this anambiguous, we need the functions $\sigma$ and $\tau$.

\subsection{Strings}
For each arrow $a \in Q_1$ we introduce a formal inverse $a^-$.
We extend the maps $s,t$ by defining $s(a^-) := t(a)$ and
$t(a^-) := s(a)$.
We also set $(a^-)^- = a$.
Let $Q_1^- = \{ a^- \mid a \in Q_1 \}$ be the set of
\emph{inverse arrows}.
Now a \emph{string} $C$ of length $l(C) := m \ge 1$ is an $m$-tuple
$$
C = (c_1,\ldots,c_m)
$$
such that the following hold:
\begin{itemize}

\item
$c_i \in Q_1 \cup Q_1^-$ for all $1 \le i \le m$;

\item
$s(c_i) = t(c_{i+1})$ for all $1 \le i \le m-1$;

\item
$c_i \not= c_{i+1}^-$ for all $1 \le i \le m-1$.

\item
$\{ c_ic_{i+1},\; c_{i+1}^-c_i^- \mid 1 \le i \le m-1 \} \cap I = \varnothing$.

\end{itemize}
We often just write $C = c_1 \cdots c_m$ instead of $C = (c_1,\ldots,c_m)$.
Let $C^- := (c_m^-,\ldots,c_1^-)$ be the \emph{inverse} of $C$, which is obviously again a string.

Additionally, for each vertex $i \in Q_0$ there are two strings
$1_{i,t}$ with $t = \pm 1$ of length $l(1_{i,t}) := 0$. 
We set $s(1_{i,t}) = t(1_{i,t}) = i$ and $1_{i,t}^- = 1_{i,-t}$.

Sometimes we will just write $1_i$ instead of $1_{(i,t)}$.

We extend the maps $\sigma$ and $\varepsilon$ to strings 
as follows:
\begin{itemize}

\item[(i)]
For $a \in Q_1$ define $\sigma(a^-) := \varepsilon(a)$ and $\varepsilon(a^-) := \sigma(a)$.

\item[(ii)]
For a string $C = (c_1,\ldots,c_m)$ of length $m \ge 1$, let
$\sigma(C) := \sigma(c_m)$ and $\varepsilon(C) := \varepsilon(c_1)$.

\item[(iii)]
$\sigma(1_{i,t}) := -t$ and $\varepsilon(1_{i,t}) := t$.

\end{itemize}

For strings $C = (c_1,\ldots,c_p)$ and
$D = (d_1,\ldots,d_q)$ of length $p,q \ge 1$, 
the composition of $C$ and $D$ is defined, provided
$(c_1,\ldots,c_p,d_1,\ldots,d_q)$ is again a string.
We write then $CD = c_1 \cdots c_pd_1 \cdots d_q$.

Now let $C$ be any string.
The composition of $1_{(u,t)}$ and $C$ is defined 
if $t(C) = i$ and $\varepsilon(C) = t$.
In this case, we write $1_{(i,t)}C = C$.
The composition of $C$ and $1_{(i,t)}$ is defined if
$s(C) = i$ and $\sigma(C) = -t$.
In this case we write $C1_{(i,t)} = C$.

If $C$ and $D$ are arbitrary strings such that the composition
$CD$ is defined, then $\sigma(C) = - \varepsilon(D)$.

For a string $C$ we write $C \sim C^-$.
This defines an equivalence relation on the set of all strings.
Let $\cS$ denote a set of representatives of all equivalence classes
of strings for $A$.

A string $C$ is a \emph{direct string} if $C$ is of length
$0$ or if it does not contain any inverse arrows. 
A direct string $C$ is \emph{right-bounded} 
(resp. \emph{left-bounded})
if $Ca \in I$ (resp. $aC \in I$) for all $a \in Q_1$.

When visualizing a string we draw an arrow $a \in Q_1$
often pointing from northeast to southwest:
$$
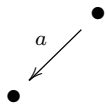
\xymatrix@-1ex{
&  \bullet \ar[dl]_a
\\
\bullet &
}
$$
Instead of the bullets one often displays the numbers $i := s(a)$ and
$j := t(a)$:
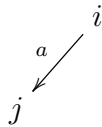
$$
\xymatrix@-1ex{
&  i \ar[dl]_a
\\
j &
}
$$
On the other hand, an inverse arrow $a^- \in Q_1^-$ is pointing from
northwest to southeast:
$$
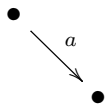
\xymatrix@-1ex{
\bullet \ar[dr]^a &
\\
& \bullet
}
$$
Note that in this picture the arrow $a^-$ carries just the label $a$.

\subsection{Example}\label{example:1}
Let again $A = KQ/I$, where $Q$ is the quiver
$$
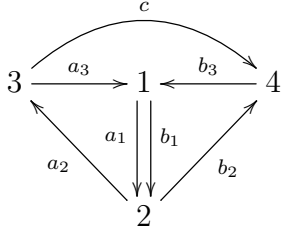
\xymatrix@+1pc{
3\ar@/^2pc/[rr]^c \ar[r]^{a_3} & 1 \ar@<-0.5ex>[d]_(.4){a_1} \ar@<0.5ex>[d]^(.4){b_1}
&4 \ar[l]_{b_3}\\
& 2 \ar[ul]^{a_2}\ar[ur]_{b_2}
}
$$
and $I$ is the ideal in $KQ$ generated by
$a_1a_3,a_2a_1,a_3a_2,b_1b_3,b_2b_1,b_3b_2$. 
Then
$$
C = a_1^-b_1a_3c^-b_2a_1b_1^- = (a_1^-,b_1,a_3,c^-,b_2,a_1,b_1^-)
$$
is a string, which looks as follows:
$$
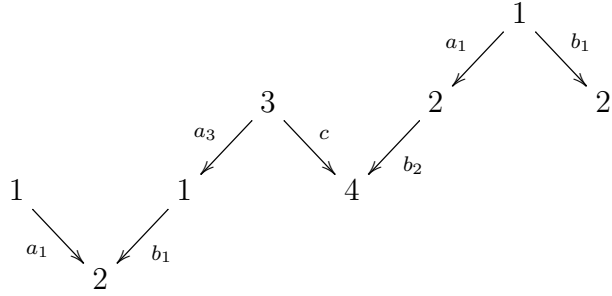
\xymatrix@-1ex{
&&&&&& 1 \ar[dl]_{a_1}\ar[dr]^{b_1}
\\
&&& 3 \ar[dl]_{a_3}\ar[dr]^c && 2 \ar[dl]^{b_2} && 2
\\
1 \ar[dr]_{a_1} && 1 \ar[dl]^{b_1} && 4
\\
& 2
}
$$

\subsection{String modules}
Let $C = (c_1,\ldots,c_m)$ be a string of length $m \ge 1$.
We define a \emph{string module} $M(C)$ as follows:
The module $M(C)$ has a \emph{standard basis}
$(b_1,\ldots,b_{m+1})$.
The generators of the algebra $A$ act on this basis as follows:
For $i \in Q_0$ and $1 \le j \le m+1$ we have
$$
e_ib_j := 
\begin{cases}
b_j & \text{if $t(c_j) = i$ and $1 \le j \le m$},
\\
b_j & \text{if $s(c_m) = i$ and $j=m+1$},
\\
0 & \text{otherwise}.
\end{cases}
$$
and for $a \in Q_1$ and $1 \le j \le m+1$ we have
$$
ab_j := 
\begin{cases}
b_{j-1} & \text{if $a = c_{j-1}$ and $2 \le j \le m+1$},
\\
b_{j+1} & \text{if $a^- = c_{j+1}$ and $1 \le j \le m$},
\\
0 & \text{otherwise}.
\end{cases}
$$

For strings $E_1$ and $E_2$ with $E_1 \sim E_2$, let
$$
\phi_{E_1,E_2}\df M(E_1) \to M(E_2)
$$ 
be the obvious canonical isomorphism.
(If $E_1 = E_2$, then $\phi_{E_1,E_2}$ is just the identity.
Let $E_1 = E_2^-$, and let $(b_1,\ldots,b_m)$ 
(resp. $(b_1',\ldots,b_m')$) be the
standard
basis of $M(E_1)$ (resp. $M(E_2)$).
Then $\phi_{E_1,E_2}(b_i) = b_{m-i+1}'$ for
$1 \le i \le m$.)

\subsection{Bands}
A \emph{band} for $A$ is a string
$B$ such that
the following hold:
\begin{itemize}

\item
$l(B) \ge 2$;

\item
$B^t$ is a string for all $t \ge 1$;

\item
$B$ is not of the form $C^s$ for some string $C$ and some 
$s \ge 2$.

\end{itemize}

Let $B$ be a band, and let $C$ and $D$ be strings such that
$B = CD$.
Then $DC$ is a \emph{rotation} of $B$.
Obviously, any rotation of $B$ is again a band.
We write 
$$
B \sim_b DC \sim_b B^-.
$$
This yields an equivalence relation on the set of all bands for $A$.
Let $\cB$ be a set of representatives of all equivalence classes of bands for $A$.

As an example, let $A = KQ/I$ as in Section~\ref{example:1}.
Then 
$$
B = c^-b_3^-a_1^-b_1a_1^-b_1a_3
$$ 
is a band.

\subsection{Band modules}
Now let $B = (c_1,\ldots,c_m)$ be a band, and let $\lambda \in K^*$.
We define a \emph{band module} $M(B,\lambda,1)$ as follows:
The module $M(B,\lambda,1)$ has a \emph{standard basis}
$(b_1,\ldots,b_{m})$.
The generators of the algebra $A$ act on this basis as follows:
For $i \in Q_0$ and $1 \le j \le m$ we have
$$
e_ib_j := 
\begin{cases}
b_j & \text{if $t(c_j) = i$ and $1 \le j \le m$},
\\
0 & \text{otherwise}.
\end{cases}
$$
and for $a \in Q_1$ and $1 \le j \le m$ we have
$$
ab_j := 
\begin{cases}
b_{j-1} & \text{if $a = c_{j-1}$ and $2 \le j \le m$},
\\
\lambda b_m & \text{if $a = c_m$ and $j = 1$},
\\
b_{j+1} & \text{if $a^- = c_{j+1}$ and $1 \le j \le m-1$},
\\
\lambda b_1 & \text{if $a^- = c_m$ and $j = m$},
\\
0 & \text{otherwise}.
\end{cases}
$$
For $q \ge 2$ and $\lambda \in K^*$ there are also band 
modules $M(B,\lambda,q)$.
They do not play a major role in this article, so we omit their
definition.
Let us just mention that they form Auslander-Reiten sequences
$$
0 \to M(B,\lambda,1) \to M(B,\lambda,2) \to M(B,\lambda,1) \to 0
$$
and
$$
0 \to M(B,\lambda,q) \to M(B,\lambda,q-1) \oplus 
M(B,\lambda,q+1) \to M(B,\lambda,q) \to 0
$$
for $q \ge 2$.
For $q \ge 1$, we say that $M(B,\lambda,q)$ has \emph{quasi-length}
$q$.

\subsection{Classification of modules}

The following classification theorem was first proved by Wald
and Waschb\"usch \cite{WW} using covering techniques. 
There is an alternative proof by Butler and Ringel \cite{BR} using
functorial filtrations. 
Both articles \cite{BR} and \cite{WW} also contain a
combinatorial description of all Auslander-Reiten sequences for string algebras.
Recall that all gentle algebras are string algebras.

\begin{Thm}
Let $A = KQ/I$ be a gentle algebra.
The modules $M(C)$ and $M(B,\lambda,q)$
with $C \in \cS$, $B \in \cB$, $\lambda \in K^*$ and
$q \ge 1$ are a complete set of pairwise non-isomorphic representatives of 
isomorphism classes of indecomposable modules in $\md(A)$.
\end{Thm}

For string modules 
we have $M(C_1) \cong M(C_2)$ if and only if $C_1 \sim C_2$, and
for band modules we have
$M(B_1,\lambda_1,q_1) \cong M(B_2,\lambda_2,q_2)$ if and only
if $B_1 \sim_b B_2$, $\lambda_1 = \lambda_2$ and $q_1 = q_2$.

\subsection{Homomorphisms}
For a string $C$ we define
$\cS(C)$ as the set of triples $(D,E,F)$ such that
the following hold:
\begin{itemize}

\item[(i)]
$C = DEF$;

\item[(ii)]
Either $l(D) = 0$, or $D = D'a^-$ for some $a \in Q_1$ and some
string $D'$;

\item[(iii)]
Either $l(F) = 0$, or $F = bF'$ for some $b \in Q_1$ and some
string $F'$.

\end{itemize}
Following our convention for displaying strings, a triple
$(D,E,F) \in \cS(C)$ with $l(D),l(F) \ge 1$ yields the following picture,
where
the left (resp. right) hand red line stands for the string $D'$
(resp. $F'$), and the blue line stands for $E$.
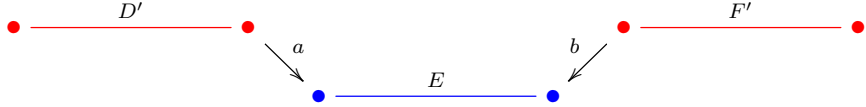
$$
\xymatrix@-2.0ex{
\red\bullet \ar@[red]@{-}[rrrr]^{D'} &&&& \red\bullet
\ar[dr]^a
&&&&&& 
\red\bullet \ar[dl]_b\ar@[red]@{-}[rrrr]^{F'} &&&& \red\bullet
\\
&&&&& \blue\bullet \ar@[blue]@{-}[rrrr]^E &&&& \blue\bullet
}
$$
We clearly see that $M(C)$ has a submodule isomorphic to
$M(E)$ and that the corresponding factor module 
is isomorphic to $M(D') \oplus M(F')$.
Let 
$$
\iota_{(D,E,F)}\df M(E) \to M(C)
$$ 
be the obvious canonical inclusion.

Dually, for a string $C$ we define
$\cF(C)$ as the set of triples $(D,E,F)$ such that
the following hold:
\begin{itemize}

\item[(i)]
$C = DEF$;

\item[(ii)]
Either $l(D) = 0$, or $D = D'a$ for some $a \in Q_1$ and some string $D'$;

\item[(iii)]
Either $l(F) = 0$, or $F = b^-F'$ for some $b \in Q_1$ and some string $F'$.

\end{itemize}
For such a
$(D,E,F) \in \cF(C)$ with $l(D),l(F) \ge 1$ we get the
following picture, where
the left (resp. right) hand green line stands for the string $D'$
(resp. $F'$), and the blue line stands for $E$.
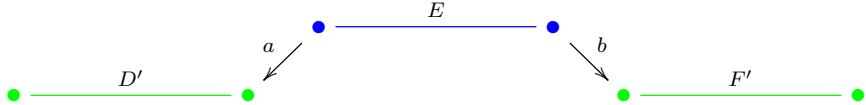
$$
\xymatrix@-2.0ex{
&&&&& \blue\bullet \ar[dl]_a\ar@[blue]@{-}[rrrr]^E &&&& 
\blue\bullet \ar[dr]^b
\\
\green\bullet \ar@[green]@{-}[rrrr]^{D'} &&&& \green\bullet
&&&&&& 
\green\bullet \ar@[green]@{-}[rrrr]^{F'} &&&& \green\bullet
}
$$
Then $M(C)$ has a submodule isomorphic to
$M(D') \oplus M(F')$ and the corresponding factor module
is isomorphic to $M(E)$.
Let 
$$
\pi_{(D,E,F)}\df M(C) \to M(E)
$$
be the obvious canonical projection.

For a pair $(C_1,C_2)$ of strings we call a pair
$$
h = ((D_1,E_1,F_1),(D_2,E_2,F_2)) \in \cF(C_1) \times \cS(C_2)
$$
\emph{admissible} if $E_1 = E_2$ or $E_1 = E_2^-$.

Suppose that $h$ is admissible.
For $E_1=E_2$,
$h$ is \emph{$2$-sided} if $l(D_i) \ge 1$ and
$l(F_j) \ge 1$ for at least one $i \in \{1,2\}$ and at least one $j \in \{1,2\}$.
For $E_1 = E_2^-$, $h$ is \emph{$2$-sided} if
$((D_1,E_1,F_1),(F_2^-,E_2^-,D_2^-))$ is $2$-sided.

Let $h$ be admissible as above, and let
$$
f_h := 
\iota_{(D_2,E_2,F_2)} \circ \phi_{E_1,E_2} \circ \pi_{(D_1,E_1,F_1)}\df
M(C_1) \to M(C_2)
$$
be the associated \emph{standard homomorphism}.
We call $f_h$ \emph{oriented} if $E_1 = E_2$.
Furthermore, $f_h$ is \emph{$2$-sided} if $h$ is $2$-sided.
Otherwise, $f_h$ is \emph{$1$-sided}.

The following picture describes $f_h$ for the case
$E_1 = E_2$ and $l(D_i),l(F_i) \ge 1$ for $i=1,2$.
$$
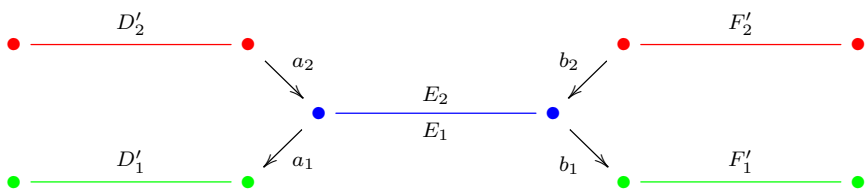
\xymatrix@-2.0ex{
\red\bullet \ar@[red]@{-}[rrrr]^{D_2'} &&&& \red\bullet
\ar[dr]^{a_2}
&&&&&& 
\red\bullet \ar[dl]_{b_2}\ar@[red]@{-}[rrrr]^{F_2'} &&&& \red\bullet
\\
&&&&& \blue\bullet \ar[dl]^{a_1}\ar@[blue]@{-}[rrrr]^{E_2}_{E_1} &&&& 
\blue\bullet \ar[dr]_{b_1}
\\
\green\bullet \ar@[green]@{-}[rrrr]^{D_1'} &&&& \green\bullet
&&&&&& 
\green\bullet \ar@[green]@{-}[rrrr]^{F_1'} &&&& \green\bullet
}
$$
Thus we have 
$C_1 = D_1'a_1E_1b_1^-F_1'$ and 
$C_2 = D_2'a_2^-E_2b_2F_2'$.
Furthermore, it follows that $a_1a_2, b_1b_2 \in I$.

Depending if some of the four strings $D_1,F_1,D_2,F_2$
are of length $0$ or not, there are $16$ different types of
oriented
standard homomorphisms.

\begin{Thm}[{\cite{CB1}}]\label{thm:CB}
For $M$ and $N$ string modules, the set of standard homomorphisms
$M \to N$ is a basis
of $\Hom_A(M,N)$.
\end{Thm}

In this article, we are mainly concerned with the question if certain
homomorphism spaces $\Hom_A(M,N)$ are zero or not.
The actual dimension of these spaces does not matter.

For a band module $M = M(B,\lambda,q)$ 
and an arbitrary indecomposable $A$-module $N$,
the conditions $\Hom_A(M,N) \not= 0$ and $\Hom_A(N,M) \not= 0$
do not depend on the quasi-length $t$.
(This follows from the description of the Auslander-Reiten sequences
involving band modules, see for example \cite{BR}.)
Therefore we can restrict our attention to band modules of
quasi-length $1$.

Krause \cite{K} extended Theorem~\ref{thm:CB} to homomorphisms
also involving band modules.
We just recall a special case here, where we only consider band
modules of quasi-length $1$.

For a band $B$ let 
\begin{align*}
\cS^\infty(B) &:= \{ (D,E,F) \in \cS(B^t) \mid 1 \le l(D),l(F) \le l(B),\;
t \ge 1 \},
\\
\cF^\infty(B) &:= \{ (D,E,F) \in \cF(B^t) \mid 1 \le l(D),l(F) \le l(B),\;
t \ge 1 \}.
\end{align*}

Let $B_1$ and $B_2$ be bands, and let $C$ be a string.
Let 
$$
h = ((D_1,E_1,F_1),(D_2,E_2,F_2)) 
$$
be an element in
$\cF^\infty(B_1) \times \cS(C)$,
$\cF(C) \times \cS^\infty(B_1)$ or
$\cF^\infty(B_1) \times \cS^\infty(B_2)$.
Then $h$ is \emph{admissible} if $E_1 = E_2$ or $E_1 = E_2^-$.
In this case, one can again define a \emph{standard homomorphism}
$f_h\df M(B_1,\lambda_1,1) \to M(C)$, $f_h\df M(C) \to M(B_1,\lambda_1,1)$ or
$f_h\df M(B_1,\lambda_1,1) \to M(B_2,\lambda_2,1)$, respectively.
All of these are \emph{$2$-sided}.
This involves of course a choice of scalars $\lambda_1$ 
and/or $\lambda_2$, in case we deal with $B_1$ and/or $B_2$.
For a band module $M(B,\lambda,1)$, the identity
is also called a \emph{standard homomorphism}.
Similarly as before, we call $f_h$ \emph{oriented} if $E_1 = E_2$.
For further details we refer to \cite{K}.

\begin{Thm}[{\cite{K}}]\label{thm:Krause}
For $M$ and $N$ string modules or band modules of quasi-length $1$,
the set of standard homomorphisms $M \to N$ 
is a basis of $\Hom_A(M,N)$.
\end{Thm}

\subsection{Auslander-Reiten translation of string modules}\label{subsec:taustring}
\begin{figure}[H]
$$
\xymatrix@!@-8.6ex{
&&&&M&&&&&&&
&&&&\tau_A(M)
\\
1.&\green\bullet\ar[dr]^{a_1}\ar@[green]@{-}_{E_1}[ddl] &&&&&&\green\bullet\ar[dl]_{a_2}
\ar@[green]@{-}^{E_2}[ddr]
\\
&&\red\bullet \ar@[red]@{-}[rrrr] &&&& \red\bullet&&& &&&&
\red\bullet \ar@[red]@{-}[rrrr] &&&& \red\bullet
\\
\green\bullet&&&&&&&& \green\bullet
\\
&&&&&&&
&&&&
\blue\bullet&&&&&&&
&\blue\bullet
\\
2.&&\red\bullet \ar@[red]@{-}[rrrr] &&&& \red\bullet&&& &&&&
\red\bullet \ar[dl]^{b_1}\ar@[red]@{-}[rrrr] &&&& \red\bullet\ar[dr]_{b_2}
\\
&&&&&&&&&&&
&\blue\bullet\ar@[blue]@{-}[luu]^{F_1} &&&&&& \blue\bullet\ar@[blue]@{-}[ruu]_{F_2}
\\
&&&&&&&\green\bullet\ar[dl]_{a_2}
\ar@[green]@{-}[ddr]^{E_2} &&&&
\blue\bullet
\\
3.&&\red\bullet \ar@[red]@{-}[rrrr] &&&& \red\bullet&&& &&&&
\red\bullet \ar[dl]^{b_1}\ar@[red]@{-}[rrrr] &&&& \red\bullet
\\
&&&&&&&&\green\bullet&&&&\blue\bullet\ar@[blue]@{-}[uul]^{F_1}
\\
4.&\green\bullet\ar[dr]^{a_1}\ar@[green]@{-}[ddl]_{E_1} &&&&&&&&&&
&&&&&&&&\blue\bullet
\\
&&\red\bullet \ar@[red]@{-}[rrrr] &&&& \red\bullet&&& &&&&
\red\bullet \ar@[red]@{-}[rrrr] &&&& \red\bullet\ar[dr]_{b_2}
\\
\green\bullet&&&&&&&&&&&
&&&&&&&\blue\bullet\ar@[blue]@{-}[ruu]_{F_2}
\\
5.&\green\bullet\ar@[green]@{-}[ldd]_{E_1}&&&&&&&&&&
&&&&&&&&\blue\bullet\ar@[blue]@{-}[ldd]^{F_2}\\
&&&&&&\\
\green\bullet&&&&&&&&&&&&&&&&&&\blue\bullet
}
$$
\caption{The Auslander-Reiten translation of string modules}\label{fig:ARtranslation}
\end{figure}
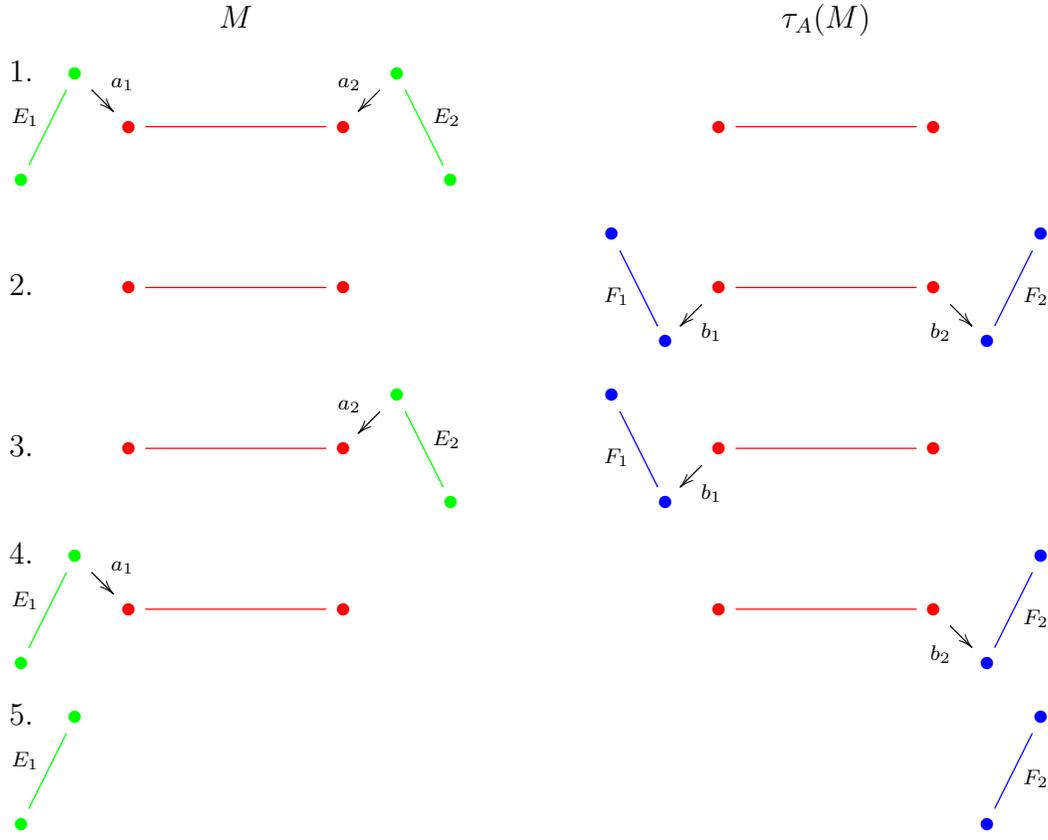

Let $A$ be a gentle algebra, and let
$M \in \md(A)$ be a non-projective string module.
It follows that $\tau_A(M)$ is also a string module,
and that we are in one
of the five situations displayed in Figure~\ref{fig:ARtranslation}, see \cite[Section~3]{BR}.
(We use here the same way of illustrating strings and 
string modules as in \cite[Section~3]{Sch}.)
The subfactor of $M$ and $\tau_A(M)$ defined by the string between the two red points is called the \emph{core} of $M$.
(In the 5th case, the core is just the $0$-module.)
The core of $M$ does not change under the Auslander-Reiten translation.

The strings $E_i$ in Figure~\ref{fig:ARtranslation} are left-bounded direct strings, 
and the strings $F_i$ are right-bounded direct strings.  
The strings $E_1a_1^-$ and $a_2E_2^-$ are 
\emph{hooks} in the sense of \cite{BR}, and
the strings $F_1^-b_1$ and $b_2^-F_2$ are \emph{cohooks}
in the sense of \cite{BR}.

For each arrow $a  = a_1 = b_2 \in Q_1$ there is exactly one Auslander-Reiten
sequence of type 5. 
In this case, there is a string 
$$
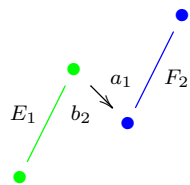
\xymatrix@!@-3.2ex{
&&& \blue\bullet \ar@[blue]@{-}[ldd]^{F_2}
\\
&\green\bullet\ar@[green]@{-}[ldd]_{E_1}\ar[dr]_{b_2}^{a_1}
\\
&& \blue\bullet
\\
\green\bullet &
}
$$
which yields the middle term of an Auslander-Reiten sequence
$$
0 \to M(F_2) \to M(E_1a^-F_2) \to M(E_1) \to 0.
$$
All other Auslander-Reiten sequences involving string modules are of types $1,\ldots,4$,
and their
middle terms are a direct sum of two indecomposable string modules.
For details we refer to \cite{BR}.

\subsection{Auslander-Reiten formulas}
The following is a well known statement from Auslander-Reiten
theory, see for example \cite{ARS,ASS,R}.

\begin{Thm}[Auslander,Reiten]\label{ARformula1}
Let $A$ be a finite-dimensional basic algebra.
For $M,N \in \md(A)$ the following hold:
\begin{itemize}

\item[(i)]
$\Ext_A^1(M,N) \cong D\overline{\Hom}_A(N,\tau_A(M)) 
\cong D\sHom_A(\tau_A^{-1}(N),M)$.

\item[(ii)]
If $\pdim(M) \le 1$, then 
$\Ext_A^1(M,N) \cong D\Hom_A(N,\tau_A(M))$.

\item[(iii)]
If $\idim(N) \le 1$, then
$\Ext_A^1(M,N) \cong D\Hom_A(\tau^{-1}(N),M)$.
\end{itemize}
\end{Thm}

\begin{Lem}\label{lem:ARband}
Let $A$ be a gentle algebra.
For any band module $M \in \md(A)$ the following hold:
\begin{itemize}

\item[(i)]
$\pdim(M) \le 1$ and $\idim(M) \le 1$;

\item[(ii)]
$\tau_A(M) \cong M$.

\end{itemize}
\end{Lem}

\begin{proof}
(i): This is well known, see for example \cite[Corollary~3.6]{BS}.

(ii): This is proved for example in \cite[Section~3]{BR}.
\end{proof}

Note that
part (ii) of the above lemma holds also for all
string algebras $A$.

\begin{Cor}\label{cor:ARformula2}
Let $A$ be a gentle algebra, and let
$M,N \in \md(A)$.
If $M$ is a band module, then
$$
\Ext_A^1(N,M) \cong 
D\sHom_A(\tau^{-1}(M),N) \cong 
D\Hom_A(\tau_A^{-1}(M),N) \cong
D\Hom_A(M,N)
$$
and 
$$
\Ext_A^1(M,N) 
\cong D\overline{\Hom}_A(N,\tau_A(M)) \cong D\Hom_A(N,\tau_A(M))
\cong D\Hom_A(N,M).
$$
\end{Cor}

\subsection{Rank functions for gentle algebras}\label{subsec:rank}

Let $A = KQ/I$ be a gentle algebra, and
let $\bd \in \N^n$ be a dimension vector.
A map $r\df Q_1 \to \N$ is a \emph{rank function} for $(A,\bd)$
if the following hold:
\begin{itemize}

\item[(i)]
$r(a) \le \min\{ d_{s(a)},d_{t(a)} \}$ for all $a \in Q_1$;

\item[(ii)]
Let $a,b \in Q_1$ with $s(a) = t(b)$ and $ab \in I$.
Then $r(a) + r(b) \le d_{s(a)}$.

\end{itemize}
(Using a slightly different wording, this definition appears in
\cite[Section~5]{CC}.)

For $M \in \md(A)$ the \emph{rank function} 
of $M$ is defined by
\begin{align*}
r_M\df Q_1 &\to \N
\\
a &\mapsto \rk(M_a).
\end{align*}
One easily checks that $r_M$ is a rank function for $(A,\bd)$ where
$\bd = \dimv(M)$.
Furthermore, each rank function for $(A,\bd)$ is obtained in this
way.

The following lemma is well known and follows directly from the
definitions of string and band modules.

\begin{Lem}\label{lem:stringbandrank}
Let $A$ be a gentle algebra.
The number of string modules
in a direct sum decomposition of $M \in \md(A)$ into indecomposable modules is
$$
\dim(M) - \sum_{a \in Q_1} r_M(a).
$$
\end{Lem}

\begin{proof}
It follows directly from the definition of a string module
$M$ that 
$$
\dim(M) - \sum_{a \in Q_1} r_M(a) = 1.
$$
For a band module $M$ we have 
$$
\dim(M) - \sum_{a \in Q_1} r_M(a) = 0.
$$
Since each $A$-module is isomorphic to a direct sum of string modules and
band modules, the claim follows.
\end{proof}

Let $r$ and $r'$ be rank functions for $(A,\bd)$.
We write $r \le r'$ if $r(a) \le r'(a)$ for all $a \in Q_1$.
This defines a partial order on the set of rank functions for $(A,\bd)$.

For a rank function $r$ for $(A,\bd)$ let
$$
\md(A,\bd,r) := \{ M \in \md(A,\bd) \mid r_M \le r \}.
$$
This is a non-empty closed subset of $\md(A,\bd)$.


\section{Schemes of complexes}\label{sec:complex}


As already mentioned in the introduction, the 
study of schemes of modules over gentle algebras can (to a large extent) be reduced to schemes of complexes.
This section deals with all necessary results on schemes of complexes.

\subsection{Definition of schemes of complexes}\label{subsec:complex1}
For $n \ge 1$ let 
$$
C_n := KQ/I,
$$ 
where $Q$ is the quiver 
$$
\xymatrix{
1 \ar[r]^{a_1} & 2 \ar[r]^{a_2} & \cdots \ar[r]^<<<<<{a_{n-2}} & n-1
\ar[r]^<<<<<{a_{n-1}} & n
}
$$
and $I$ is the ideal generated by all paths of length $2$.
(For $n=1$, $Q$ has just one vertex and no arrows. 
For $n=1,2$, we set $I = 0$.)

For $n \ge 1$ let 
$$
\widetilde{C}_n := KQ/I,
$$ 
where $Q$ is the quiver 
$$
\xymatrix{
1 \ar[r]^{a_1} & 2 \ar[r]^{a_2} & \cdots \ar[r]^<<<<<{a_{n-2}} & n-1
\ar[r]^<<<<<{a_{n-1}} & n \ar@/^2pc/[llll]^{a_n}
}
$$
and $I$ is the ideal generated by all paths of length $2$.
For $\widetilde{C}_n$ we adopt the convention that all indices
are meant modulo $n$.

Let $A$ be one of the algebras $C_n$ or $\widetilde{C}_n$.
By \emph{scheme of complexes} we mean the affine schemes 
$\md(A,\bd)$ with $\bd \in \N^n$.
This definition is a bit more general than the one used by
De Concini and Strickland \cite{DS}, who consider only the case
$C_n$.

The representation theory of $A$ is extremely well understood.
Obviously, $A$ is a representation-finite gentle algebra.
So all its indecomposable modules are string modules.
For each vertex $i \in Q_0$ there is a simple module $S_i$ and
an indecomposable projective modules $P_i$, and these are all
indecomposable $A$-modules up to isomorphism.
The modules $S_1,\ldots,S_n,P_1,\ldots,P_n$  
are pairwise non-isomorphic, with the exception of $P_n$ being
equal to $S_n$ in case $A = C_n$.
Using the usual notation for string modules, for each $i \in Q_0$ 
we have
$S_i = M(e_i)$ and 
$$
P_i = 
\begin{cases}
M(e_i) & \text{if $A = C_n$ and $i=n$},
\\
M(a_i) & \text{otherwise}.
\end{cases}
$$
It is straightforward to compute homomorphism spaces and extension groups between $A$-modules.
All this can be proved in an elementary fashion using mainly Linear Algebra.
The next two lemmas contain all the homological data we need.

\begin{Lem}\label{lem:complexhom}
Let $A$ be one of the algebras $C_n$ or $\widetilde{C}_n$.
The only pairs $(X,Y)$ of indecomposable $A$-modules with $\Hom_A(X,Y) \not= 0$ are
$$
(S_i,S_i),\quad (P_i,P_i),\quad (P_i,S_i), \quad 
(S_{t(a)},P_{s(a)}), \quad (P_{t(a)},P_{s(a)}),  
$$
where $i \in Q_0$ and $a \in Q_1$.
In these cases, we have $\dim \Hom_A(X,Y) = 1$ with only one exception for $A = \widetilde{C}_1$, where we have 
$\dim \Hom_A(P_1,P_1) = 2$.
\end{Lem}

\begin{Lem}\label{lem:complexextgroups}
Let $A$ be one of the algebras $C_n$ or $\widetilde{C}_n$.
The only pairs $(X,Y)$ of indecomposable $A$-modules with
$\Ext_A^1(X,Y) \not= 0$ are 
$$
(S_{s(a)},S_{t(a)}),
$$
where $a \in Q_1$.
In these cases, we have $\dim \Ext_A^1(X,Y) = 1$.
\end{Lem}

Let $A$ be one of the algebras $C_n$ or $\widetilde{C}_n$.
Let $\bd = (d_1,\ldots,d_n) \in \N^n$ be a dimension vector, and let $r$ be a
rank function for $(A,\bd)$.
Then there exists a unique (up to isomorphism) $A$-module $M = M_{\bd,r}$ with $\dimv(M) = \bd$ and $r_M = r$.
More precisely, we have
$$
M_{\bd,r} = \bigoplus_{a \in Q_1} P_{s(a)}^{r(a)} \oplus 
\bigoplus_{i \in Q_0} S_i^{d_i-r_i}
$$
where
$$
r_i :=
\begin{cases}
r(a_i) + r(a_{i-1}) & \text{if $2 \le i \le n-1$},
\\
r(a_1) & \text{if $i=1$ and $A = C_n$},
\\
r(a_1) + r(a_n) & \text{if $i=1$ and $A = \widetilde{C}_n$},
\\
r(a_{n-1}) & \text{if $i=n$ and $A = C_n$},
\\
r(a_n) + r(a_{n-1}) & \text{if $i=n$ and $A = \widetilde{C}_n$}.
\end{cases}
$$

\begin{Prop}\label{prop:orbitclosure}
Let $A$ be  one of the algebras $C_n$ or $\widetilde{C}_n$, and
let $\bd \in \N^n$. 
For each rank function $r$ for $(A,\bd)$ we have
$$
\md(A,\bd,r) = \overline{\cO_{M_{\bd,r}}}.
$$
\end{Prop}

\begin{proof}
For $M \in \md(A,\bd)$ and $a \in Q_1$ we have
$$
r_M(a) = d_{s(a)} - \dim \Hom_A(S_{s(a)},M).
$$
Now the claim follows from
\cite[Theorem~1 and its Corollary]{Z}.
\end{proof}

\begin{Cor}\label{cor:rankmaximal}
Let $A$ be  one of the algebras $C_n$ or $\widetilde{C}_n$, and
let $\bd \in \N^n$. 
For each $M = M_{\bd,r}$ the following are equivalent:
\begin{itemize}

\item[(i)]
$\cO_M$ is open;

\item[(ii)]
The rank function $r$ is maximal.

\end{itemize}
\end{Cor}

\begin{Cor}\label{cor:complexcomp}
Let $A$ be  one of the algebras $C_n$ or $\widetilde{C}_n$, and
let $\bd \in \N^n$. 
Then
$$
\irr(A,\bd) = \{ \md(A,\bd,r) \mid \text{$r$ is a maximal rank function
for $(A,\bd)$} \}.
$$
\end{Cor}

\begin{Lem}\label{lem:simplesummands}
Let $A = KQ/I$ be  one of the algebras $C_n$ or $\widetilde{C}_n$. 
For $M \in \md(A)$ we have 
$$
\sum_{a \in Q_1} r_M(a) \le \frac{1}{2} \dim(M).
$$
Furthermore, this becomes an equality if and only if 
$M$ does not have a simple direct summand.
\end{Lem}

\begin{proof}
We have
$$
M \cong \bigoplus_{a \in Q_1} P_{s(a)}^{m_a} \oplus
\bigoplus_{i \in Q_0} S_i^{m_i}
$$
for some $m_a,m_i \ge 0$.
Thus for $a \in Q_1$ we have $r_M(a) = m_a$.
This implies
$$
\dim(M) = \sum_{a \in Q_1} 2r_M(a) + \sum_{i \in Q_0} m_i
\geq \sum_{a \in Q_1} 2r_M(a).
$$
The claim follows.
\end{proof}

\subsection{Rigid and $\tau$-rigid modules}\label{subsec:Rigid}

\begin{Prop}[Rigid modules]\label{prop:complexrigid}
Let $A$ be one of the algebras $C_n$ or $\widetilde{C}_n$, and
let $\bd \in \N^n$.
For $M \in \md(A,\bd)$ the following are equivalent:
\begin{itemize}

\item[(i)]
$M$ is rigid;

\item[(ii)]
$M$ does not have a direct summand isomorphic to
$$
S_a := \bigoplus_{i \in \{ s(a),t(a) \}} S_i 
$$
for some $a \in Q_1$.

\end{itemize}
For $A = \widetilde{C}_1$ we assume now additionally that
$\bd = (d_1)$ with $d_1$ even.
Then the two conditions above are equivalent to the following:
\begin{itemize}

\item[(iii)]
$\cO_M$ is open.

\end{itemize}
\end{Prop}

\begin{proof}
The equivalence (i) $\iff$ (ii) follows from
Lemma~\ref{lem:complexextgroups}.
The implication (i) $\implies$ (iii) is true in general and
follows from Voigt's Lemma~\ref{lem:Voigt}.

(iii) $\implies$ (ii):
Assume that (ii) does not hold.
Thus there is an arrow $a$ such that $S_a$ is isomorphic
to a direct summand of $M$. 
For $A \not= \widetilde{C}_1$ there is a non-split short 
exact sequence
$$
0 \to S_{t(a)} \to P_{s(a)} \to S_{s(a)} \to 0.
$$
Thus $M$ is properly contained in the orbit closure of
$$
N := P_{s(a)} \oplus M/S_a. 
$$
For $A = \widetilde{C}_1$ and $\bd = (d_1)$ with $d_1$ even,
we get that $M$ has a direct summand isomorphic to 
$S_{s(a)} \oplus S_{s(a)}$. 
(Here we used that $d_1$ is even.)
We get a 
non-split short 
exact sequence
$$
0 \to S_{s(a)} \to P_{s(a)} \to S_{s(a)} \to 0.
$$
Thus $M$ is properly contained in the orbit closure of
$$
N := P_{s(a)} \oplus M/(S_{s(a)} \oplus S_{s(a)}). 
$$
In both case, this shows that $\cO_M$ is not open.
\end{proof}

The module $S_a$ in Proposition~\ref{prop:complexrigid}(ii) is
a \emph{critical summand} of type $I$ of $M$.
In Proposition~\ref{prop:complexrigid}(ii)  we have
$$
|\{ s(a),t(a) \}| 
=
\begin{cases}
1 & \text{if $A = \widetilde{C}_1$},
\\
2 & \text{otherwise}.
\end{cases}
$$
Consequently, we have
$$
S_a = 
\begin{cases}
S_1 & \text{if $A = \widetilde{C}_1$},
\\
S_{s(a)} \oplus S_{t(a)} & \text{otherwise}.
\end{cases}
$$

Recall that a $\tau$-rigid module is automatically rigid.
Thus to get a decription of all $\tau$-rigid modules, it suffices
to look at rigid modules.

\begin{Prop}[$\tau$-rigid modules]\label{prop:complextaurigid}
Let $A$ be one of the algebras $C_n$ or $\widetilde{C}_n$, and
let $\bd \in \N^n$.
For a rigid $M \in \md(A,\bd)$ the following are equivalent:
\begin{itemize}

\item[(i)]
$M$ is $\tau$-rigid;

\item[(ii)]
$M$ has no direct summand isomorphic to
$$
P_a := P_{t(a)} \oplus S_{s(a)}
$$
for some $a \in Q_1$.

\end{itemize}
\end{Prop}

\begin{proof}
We have
$\tau_A(P_i) = 0$ for $i \in Q_0$ and $\tau_A(S_{s(a)}) = S_{t(a)}$
for $a \in Q_1$.

(i) $\implies$ (ii):
Assume that $M$ has a direct summand isomorphic to
$P_a$.
Then 
$$
\dim \Hom_A(M,\tau_A(M)) \ge 
\dim \Hom_A(P_{t(a)},\tau_A(S_{s(a)})) =
\dim \Hom_A(P_{t(a)},S_{t(a)}) = 1.
$$

(ii) $\implies$ (i):
Assume that $\Hom_A(M,\tau_A(M)) \not= 0$.
Thus there are indecomposable direct summands $X$ and $Y$ of
$M$ with
$\Hom_A(X,\tau_A(Y)) \not= 0$.
We get $Y \cong S_{s(a)}$ and $\tau_A(Y) \cong S_{t(a)}$ for some
$a \in Q_1$.
This implies $X \cong S_{t(a)}$ or $X \cong P_{t(a)}$.
If $X \cong S_{t(a)}$, then the rigid module $M$ has a direct summand isomorphic to $S_a$, a contradiction to Proposition~\ref{prop:complexrigid}.
If $X \cong P_{t(a)}$, then $X \oplus Y \cong P_a$.
This proves the claim.
\end{proof}

The module $P_a$ in Proposition~\ref{prop:complextaurigid}(ii) is
a \emph{critical summand} of type $II$ of $M$.

\subsection{Generic reducedness and singular locus}

\begin{Prop}\label{prop:complexreduced}
Let $A = KQ/I$ be one of the algebras 
$C_n$ or $\widetilde{C}_n$, and let
$\bd = (d_1,\ldots,d_n) \in \N^n$.
For $Z \in \irr(A,\bd)$ the following are equivalent:
\begin{itemize}

\item[(i)]
$Z$ is not generically reduced.

\item[(ii)]
$A = \widetilde{C}_1$ and $d_1$ is odd.

\end{itemize}
\end{Prop}

\begin{proof}
(i) $\implies$ (ii):
Suppose that (ii) does not hold.
Then it follows from Proposition~\ref{prop:complexrigid} that $Z$ contains a rigid module $M$.
Then $Z = \overline{\cO_M}$ and $Z$ is
generically reduced by Corollary~\ref{cor:schemes3}.

(ii) $\implies$ (i):
Assume that (ii) holds.
Then 
$$
Z = \overline{\cO_M} = \md(A,\bd) 
$$
with $M = S_1 \oplus P_1^{(d_1-1)/2}$.
In particular, $M$ is not rigid and therefore $Z$ is not generically
reduced, again by Corollary~\ref{cor:schemes3}.
\end{proof}

Proposition~\ref{prop:complexreduced} is not really original.
Using very different methods,
it is shown in
\cite[Theorem~1.7]{DS} that $\md(C_n,\bd)$ is reduced for all $\bd$.
Reducedness is in general a much stronger and harder to prove property than being generically reduced.
Also the schemes $\md(\widetilde{C}_n,\bd)$ should be
reduced provided $n \ge 2$.
A proof for $n=2$ is in \cite[Proposition~1.3]{St}.

\begin{Prop}\label{prop:complexsingular}
Let $A = KQ/I$ be one of the algebras $C_n$ or $\widetilde{C}_n$, and
let $\bd \in \N^n$.
For $M \in \md(A,\bd)$ the following are equivalent:
\begin{itemize}

\item[(i)]
$M$ is singular;

\item[(ii)]
There exist arrows $a,b \in Q_1$ with $s(a) = t(b)$ such that the module
$$
S_{ab}:= \bigoplus_{k \in \{ s(a),t(a),s(b) \}} S_k 
$$
is isomorphic to a direct summand of $M$.

\end{itemize}
\end{Prop}

\begin{proof}
Let $r$ be a rank function for $(A,\bd)$ and let
$$
M = M_{\bd,r} = \bigoplus_{a \in Q_1} P_{s(a)}^{r(a)} 
\oplus \bigoplus_{i \in Q_0} S_i^{q_i}.
$$ 
For
$A = C_n$ we adopt the convention that $P_j = S_j = 0$ and 
$r_j = q_j = 0$ for all $j \notin Q_0$, and for
$A = \widetilde{C}_n$ we use all indices modulo $n$.

Case 1: $A = C_1$ or $A = C_2$.
In this case, $\md(A,\bd)$ is always an affine space.
Therefore all modules $M$ are smooth.
On the other hand condition (ii) is never satisfied. 
This proves (i) $\iff$ (ii).

Case 2: $A = \widetilde{C}_1$.
In this case, $M$ is of the form
$$
M = P_1^{r_1} \oplus S_1^{q_1}.
$$
We have $\dim \Ext_A^1(M,M) = q_1^2$ and
$\dim \cO_M = 2r_1^2 + 2r_1q_1$.
Thus
$$
\dim T_M = \dim \cO_M + \dim \Ext_A^1(M,M) =
2r_1^2 + 2r_1q_1 + q_1^2.
$$
Now $Z = \overline{\cO_N} = \md(A,\bd)$ is irreducible, 
where 
$$
N = 
\begin{cases}
P_1^{r_1+q_1/2} & \text{if $q_1$ is even},
\\
P_1^{r_1+(q_1-1)/2} \oplus S_1& \text{if $q_1$ is odd}.
\end{cases}
$$
We get
$$
\dim(Z) =
\begin{cases}
2r_1^2 + 2r_1q_1 + 1/2 q_1^2& \text{if $q_1$ is even},
\\
2r_1^2 + 2r_1q_1 + 1/2 q_1^2-1/2 & \text{if $q_1$ is odd}.
\end{cases}
$$

This shows that $M$ is smooth if and only if $q_1 = 0$.
(Thus if $\bd = (d_1)$ is odd, then $\md(A,\bd)$ does not contain any smooth module,
and if $\bd$ is even, then there is only one smooth module up to isomorphism,
namely $M = P_1^{d_1/2}$.)
Note that (ii) holds if and only if $q_1 \ge 1$.
This proves (i) $\iff$ (ii).

Case 3: $A = \widetilde{C}_2$.
In this case, $M$ is of the form
$$
M = P_1^{r_1} \oplus P_2^{r_2} \oplus S_1^{q_1} \oplus S_2^{q_2}.
$$
Assume that $q_1 = 0$ or $q_2 = 0$.
Then $M$ is rigid and therefore smooth.
Next, assume that $q_1,q_2 \ge 1$. 
Then $M$ is contained in the intersection of at least two different irreducible components $Z_1$ and $Z_2$, with
maximal rank functions $r_1$ and $r_2$, respectively, which are defined by
\begin{align*}
r_1(a_1) &:= r_1+\min\{q_1,q_2\}, & r_1(a_2) &:= r_2,
\\
r_2(a_1) &:= r_1, & r_2(a_2) &:= r_2+\min\{q_1,q_2\}.
\end{align*}
Thus $M$ is singular.

This shows that $M$ is singular if and only if $q_1,q_2 \ge 1$.
But this condition is equivalent to (ii).

Case 4: $n \ge 3$.
Let 
\begin{align*}
H_2 &:= \{ 1 \le i \le n \mid q_i,q_{i+1} \ge 1 \text{ and }
q_{i-1} = q_{i+2} = 0 \},
\\
H_3 &:= \{ 1 \le i \le n \mid q_i,q_{i+1},q_{i+2} \ge 1 \text{ and }
q_{i-1} = q_{i+3} = 0 \}.
\end{align*}

Case 4(a):
Assume that $q_i,q_{i+1},q_{i+2} \ge 1$ and
$q_i+q_{i+2} > q_{i+1}$ for some $i$.
Similarly as in Case 3 one shows that $M$ is contained in at least two
different irreducible components of $\md(A,\bd)$.
Thus $M$ is singular.

Case 4(b): 
Assume that for all $i$ with $q_i,q_{i+1},q_{i+2} \ge 1$
we have $q_i+q_{i+2} \le q_{i+1}$.
It follows immediately that $q_{i-1} = q_{i+3} = 0$ for all such $i$.
In other words, we have $i \in H_3$.

We get that $M$ is contained in exactly one irreducible component
$Z = \overline{\cO_N}$, where $N$ is obtained from $M$ as follows:
For each $i \in H_2$ replace 
$$
S_i^{q_i} \oplus S_{i+1}^{q_{i+1}}
\text{\qquad by \qquad} 
P_i^{\min\{q_i,q_{i+1}\}} \oplus S_{i+1}^{|q_i-q_{i+1}|}.
$$
Furthermore, for each $i \in H_3$ replace
$$
S_i^{q_i} \oplus S_{i+1}^{q_{i+1}} \oplus S_{i+2}^{q_{i+2}}
\text{\qquad by \qquad}
P_i^{q_i} \oplus P_{i+1}^{q_{i+2}} \oplus S_{i+1}^{q_{i+1}-q_i-q_{i+2}}.
$$
The module $N$ is rigid and therefore smooth.

Now $M$ is smooth if and only if
$$
\dim \cO_N = \dim(Z) = \dim(T_M) = \dim \cO_M + \dim\Ext_A^1(M,M).
$$
(Note that the first and third equality always hold.)
Thus $M$ is smooth if and only if
\begin{equation}\label{eq1}
\dim \End_A(M) - \dim \End_A(N) = \dim \Ext_A^1(M,M).
\end{equation}
We have
$$
\dim \Ext_A^1(M,M) = \sum_{a \in Q_1} q_{s(a)}q_{t(a)}.
$$
Now a straightforward but lengthy calculation shows that Equation~(\ref{eq1}) holds if and only if $H_3 = \varnothing$.
More precisely, one gets that
$$
\dim \End_A(M) - \dim \End_A(N) = \dim \Ext_A^1(M,M) + 
\sum_{i \in H_3} q_iq_{i+2}.
$$
Thus $M$ is smooth if and only if $H_3 = \varnothing$.
This finishes the proof.
\end{proof}

In Proposition~\ref{prop:complexsingular}(ii)  we have
$$
|\{ s(a),t(a),s(b) \}| 
=
\begin{cases}
1 & \text{if $A = \widetilde{C}_1$},
\\
2 & \text{if $A = \widetilde{C}_2$},
\\
3 & \text{otherwise}.
\end{cases}
$$
Consequently, we have
$$
S_{ab} = 
\begin{cases}
S_1 & \text{if $A = \widetilde{C}_1$},
\\
S_1 \oplus S_2 & \text{if $A = \widetilde{C}_2$},
\\
S_{s(a)} \oplus S_{t(a)} \oplus S_{s(b)} & \text{otherwise}.
\end{cases}
$$

The singularities of the closures of the $\GL_\bd(K)$-orbits 
of the schemes
$\md(C_n,\bd)$ have been described by Lakshmibai \cite{L}
for $n=3$ and by Gonciulea \cite{Go} for arbitrary $n$.
Note the difference to Proposition~\ref{prop:complexsingular},
where we look at the singularities of the whole scheme.

\subsection{$\rho$-blocks of gentle Jacobian algebras}\label{subsec:blockssurface}
Let $A = KQ/I$ be a gentle Jacobian algebra.
It follows from the definitions that the $\rho$-blocks of $A$ 
are isomorphic to $C_1$, $C_2$ or $\widetilde{C}_3$.
We call them \emph{$1$-blocks}, $2$-\emph{blocks} or
$3$-\emph{blocks}, respectively.

A $1$-block can only occur if $A = C_1$.
Here we used that gentle Jacobian algebras are by definition connected.

Now let $A_s$ be a $1$-block or $2$-block.
Then the schemes
$\md(A_s,\bd)$ are obviously just affine spaces.
In particular, they are irreducible, and all modules $M \in \md(A_s,\bd)$ are smooth and reduced.
Furthermore, $\md(A_s,\bd)$ contains a unique $\tau$-rigid module.
In particular, $\md(A_s,\bd)$ is generically $\tau$-reduced.

Next, let $A_s$ be a $3$-block of $A$.
For convenience, we assume that $A = \widetilde{C}_3 = KQ/I$,
where $Q$ is the quiver 
$$
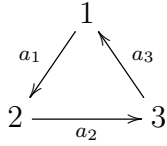
\xymatrix@-1.7pc{
&&1 \ar[dddll]_{a_1}\\
&&&&\\
&&&&\\
2 \ar[rrrr]_{a_2} &&&& 3 \ar[lluuu]_{a_3}
}
$$
and 
$I$ is generated by
the paths $a_2a_1$, $a_3a_2$ and $a_1a_3$.

For later use, we define
$$
I_3 := \{ (1,3,2), (2,1,3), (3,2,1) \}.
$$

\begin{Lem}\label{lem:3blockunique}
Let $A$ be a $1$-block, $2$-block or $3$-block as above.
For $\tau$-rigid $A$-modules $M$ and $N$ the following are
equivalent:
\begin{itemize}

\item[(i)]
$M \cong N$;

\item[(ii)]
$\dimv(M) = \dimv(N)$.

\end{itemize}
\end{Lem}

\begin{proof}
By the discussion above, the statement is clear for $1$-blocks and
$2$-block.
Thus assume $A$ is a $3$-block as above.

(i) $\implies$ (ii):
This is trivial.

(ii) $\implies$ (i):
By Proposition~\ref{prop:complextaurigid}
there are four types of $\tau$-rigid $A$-modules:
\begin{align*}
P_1^{r_1} \oplus P_2^{r_2} \oplus P_3^{r_3} && \text{type 0} 
\\
P_i^{r_i} \oplus P_j^{r_j} \oplus S_i^{s_i} && \text{type $i$} && \text{$(i,j,k) \in I_3$}
\end{align*}
where $r_i \ge 0$ and $s_i \ge 1$ for all $i$.

First, let $M$ be of type $0$ with $\dimv(M) = \bd = (d_1,d_2,d_3)$.
It follows that
\begin{align*}
r_1 + r_3 &= d_1,
\\
r_2 + r_1 &= d_2,
\\
r_3 + r_2 &= d_3.
\end{align*}
For a fixed $\bd$, this system of linear equations has exactly one
solution.
This proves (ii) $\implies$ (i) for modules of type $0$.

Next, let $M$ be of type $i$ for some $1 \le i \le 3$ 
with $\dimv(M) = \bd = (d_1,d_2,d_3)$.
It follows that
\begin{align*}
r_i + r_j + s_i &= d_i,
\\
r_j &= d_j,
\\
r_i &= d_k.
\end{align*}
For a fixed $\bd$, this system of linear equations has exactly one
solution.
This proves (ii) $\implies$ (i) for modules of type $i$.

Finally, we observe that modules of different types have always
different dimension vectors.
This finishes the proof.
\end{proof}

\begin{Lem}\label{lem:3blocksingular}
Let $A$ be a $1$-block, $2$-block or $3$-block as above.
For $M \in \md(A,\bd)$ the following are equivalent:
\begin{itemize}

\item[(i)]
$M$ is singular;

\item[(ii)]
$M$ is contained in at least two different irreducible components
of $\md(A,\bd)$.

\end{itemize}
\end{Lem}

\begin{proof}
By the discussion above, the statement is clear for $1$-blocks and
$2$-block. Thus assume $A$ is a $3$-block as above.

(i) $\implies$ (ii):
Assume $M$ is singular.
Now Proposition~\ref{prop:complexsingular} implies that
$$
M = \bigoplus_{i=1}^3 P_i^{r_i} \oplus \bigoplus_{i=1}^3 S_i^{q_i}
$$
with $q_1,q_2,q_3 \ge 1$.
Without loss of generality assume that 
$$
q_1 = \min\{ q_1,q_2,q_3 \}.
$$
It follows that $q_2 + q_3 > q_1$.
Now one proceeds as in the proof of
Proposition~\ref{prop:complexsingular} to show that
$M$ is contained in at least two different
irreducible components.

(ii) $\implies$ (i): This holds for arbitrary finite-dimensional $K$-algebras, see Proposition~\ref{prop:intersecsingular}.
\end{proof}


\section{Irreducible components for gentle algebras}
\label{sec:schemesgentle}


\subsection{Irreducible components}
Finding the irreducible components of schemes of modules over
gentle algebras is rather easy, since each of these schemes is
isomorphic to a product of schemes of complexes.

Let $A = KQ/I$ be a gentle algebra, and let $A_1,\ldots,A_t$ be its $\rho$-blocks.
For each $\rho$-block $A_s$ there is a unique 
$$
A_s' \in \{ C_n,\widetilde{C}_n \mid n \ge 1 \}
$$ 
such that there exists 
an algebra homomorphism
$$
f_s\df A_s' \to A_s
$$
with the following properties:
\begin{itemize}

\item[(i)]
$f_s$ sends vertices to vertices and arrows to arrows.

\item[(ii)]
$f_s$ is bijective on the sets of arrows.

\end{itemize}
(In (i) we think of the vertices as standard idempotents.)
This follows directly from the definition of a gentle algebra
and from the definition of a $\rho$-block.
We say that $A_s$ is of \emph{type} $A_s'$.
Let $n_s$ (resp. $n_s'$) be the number of vertices of $A_s$
(resp. $A_s'$).
For each dimension vector $\bd = (d_1,\ldots,d_{n_s})$, the homomorphism $f_s$
induces an isomorphism
$$
f_{s,\bd}\df \md(A_s,\bd) \to \md(A_s',\bd')
$$
of affine schemes,
where 
$$
\bd' = (d_{f_s(1)},\ldots,d_{f_s(n_s')}).
$$

For example, let $A = KQ$, where $Q$ is the quiver
$$
\xymatrix{
1 \ar[r]^a & 2 \ar[r]^b & 3
}
$$
So here we have $I = 0$ and $\rho = \varnothing$.
There are two $\rho$-blocks $A_1$ and $A_2$ of type $C_2$,
i.e. $A_1' = A_2' = C_2$.
Define
$f_1\df A_1' \to A_1$ by 
$1 \mapsto 1$, $2 \mapsto 2$, $a_1 \mapsto a$, and
define $f_2\df A_2' \to A_2$ by $1 \mapsto 2$, $2 \mapsto 3$ and
$a_1 \mapsto b$.
For $s=1,2$ and
a dimension vector $\bd$ for $A_s$ we have $\bd' = \bd$.

As a less trivial example, let $A = KQ/I$, where $Q$ is the quiver
$$
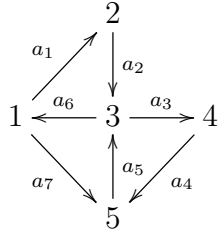
\xymatrix{
& 2 \ar[d]^{a_2}
\\
1 \ar[ur]^{a_1}\ar[dr]_{a_7} & 3 \ar[r]^{a_3}\ar[l]_{a_6} & 4 \ar[dl]^{a_4}
\\
&5 \ar[u]_{a_5}&
}
$$
and $I$ is generated by the paths $\{ a_{i+1}a_i \mid 1 \le i \le 6 \}$.
Then $A$ has only one $\rho$-block, namely $A_1 = A$, which is of type
$C_8$.
Define $f_s\df A_1' \to A_1$ by
$$
f_s(i) :=
\begin{cases}
i & \text{if $1 \le i \le 5$},
\\
3 & \text{if $i=6$},
\\
1 & \text{if $i=7$},
\\
5 & \text{if $i=8$},
\end{cases}
$$
and $f_s(a_i) := a_i$ for $1 \le i \le 7$.

For $\bd = (d_1,d_2,d_3,d_4,d_5) \in \N^5$ we get an isomorphism
$$
\md(A_s,\bd) \to \md(A_1',\bd')
$$
of affine schemes, where 
$\bd' = (d_1,d_2,d_3,d_4,d_5,d_3,d_1,d_5)$.

The following result follows almost immediately from \cite{DS}, 
see also \cite[Propositions~3.4 and 5.2]{CW}. 
Note that Carroll and Weyman \cite{CW} only consider the class
of gentle algebras admitting a colouring.
However, the result holds in general.

\begin{Prop}[{\cite{CW,DS}}]\label{prop:gentlecomponents}
Let $A$ be a gentle algebra, and let $\bd \in \N^n$.
Then we have
$$
\irr(A,\bd) = \{ \md(A,\bd,r) \mid \text{$r$ is a maximal
rank function for $(A,\bd)$} \}.
$$
\end{Prop}

\begin{proof}
Let $A_1,\ldots,A_t$ be the $\rho$-blocks of $A$.
Recall that for each $\bd$ we have an isomorphism
$$
\md(A,\bd) \to \md(A_1,\pi_1(\bd)) \times \cdots \times \md(A_t,\pi_t(\bd))
$$
which yields a bijection
$$
\irr(A,\bd) \to \irr(A_1,\pi_1(\bd)) \times \cdots \times 
\irr(A_t,\pi_t(\bd)).
$$
Now the isomorphisms 
$$
f_{s,\pi_s(\bd)}\df \md(A_s,\pi_s(\bd)) \to
\md(A_s',\pi_s(\bd)')
$$
and the description of irreducible components of varieties
of complexes (see Corollary~\ref{cor:complexcomp})
yield the result.
\end{proof}

\subsection{String and band components and generic decompositions}
Let $A = KQ/I$ be a gentle algebra.
An indecomposable irreducible component $Z$ of $\md(A,\bd)$ is
a \emph{string component} provided there is a string $C$ 
such that the orbit $\cO_{M(C)}$ is dense
in $Z$. 
In this case, $C$ is (up to equivalence of strings) uniquely determined by $Z$,
and we write $Z = Z(C)$.

An indecomposable component $Z \in \irr(A,\bd)$ is
a \emph{band component} provided there is a band $B$ 
such that the union 
$$
\bigcup_{\lambda \in K^*} \cO_{M(B,\lambda,1)}
$$
is dense in $Z$. 
In this case, $B$ is (up to equivalence of bands) 
uniquely determined by $Z$, and we write $Z = Z(B)$.
(The band modules $M(B,\lambda,q)$ are contained in the 
closure of the union
$$
\bigcup_{\lambda_1,\ldots,\lambda_q \in K^*} 
\cO_{M(B,\lambda_1,1) \oplus \cdots \oplus M(B,\lambda_t,1)},
$$
so they do no play a role here.)

Any indecomposable component $Z \in \irr(A)$ is either a string or
a band component.

For $Z \in \irr(A,\bd)$ let
$$
Z = \overline{Z(C_1) \oplus \cdots \oplus Z(C_p) \oplus Z(B_1) \oplus \cdots \oplus
Z(B_q)}
$$
be the canonical decomposition of $Z$.
Then $M$ is \emph{generic} in $Z$, if
$$
M \cong M(C_1) \oplus \cdots \oplus M(C_p) \oplus
M(B_1,\lambda_1,1) \oplus \cdots \oplus M(B_q,\lambda_q,1)
$$
with pairwise different $\lambda_1,\ldots,\lambda_q \in K^*$.

\begin{Lem}\label{lem:genericdecomp}
Let $A$ be a gentle algebra.
For $Z \in \irr(A,\bd)$ let
$$
Z = \overline{Z(C_1) \oplus \cdots \oplus Z(C_p) \oplus Z(B_1) \oplus \cdots \oplus
Z(B_q)}
$$
be the canonical decomposition of $Z$.
Then 
$c_A(Z) = q$.
\end{Lem}

\begin{proof}
Let 
$f\df\GL_\bd(K) \times (K^*)^q \to \md(A,\bd)$
be defined by
$$
(g,(\lambda_1,\ldots,\lambda_q)) \mapsto
g.(M(C_1) \oplus \cdots \oplus M(C_p) \oplus
M(B_1,\lambda_1,1) \oplus \cdots \oplus M(B_q,\lambda_q,1)).
$$
For $M \in \Ima(f)$ 
the fibre 
$f^{-1}(M)$ is obviously isomorphic to the automorphism group
${\rm Aut}_A(M)$ of $M$.
This implies
$$
\dim f^{-1}(M) = \dim \End_A(M).
$$
Thus we have 
$$
\dim \cO_M + \dim f^{-1}(M) = \dim(G).
$$
By definition 
$$
c_A(Z) = \dim(Z) - \dim \cO_M
$$
where $M$ is generic in $Z$.
By Chevelley's Theorem we have
$$
\dim(Z) + \dim f^{-1}(M) = \dim(G \times (K^*)^q) = \dim(G) + q
$$
where $M$ is again generic in $Z$.
Combining these equations yields $c_A(Z) = q$.
\end{proof}

\begin{Cor}\label{cor:comp1}
Let $A$ be a gentle algebra.
For $Z \in \irr(A,\bd)$ the following hold:
\begin{itemize}

\item[(i)]
If $Z$ is a string component, then $c_A(Z) = 0$.

\item[(ii)]
If $Z$ is a band component, then $c_A(Z) = 1$.

\end{itemize}
\end{Cor}

Note that Corollary~\ref{cor:comp1} is just a special case of
Lemma~\ref{lem:CCtame}.

\subsection{Generically reduced components}
\label{subsec:genreduced}

\begin{Thm}\label{thm:main2b}
Let $A$ be a gentle algebra, and let $A_1,\ldots,A_t$ be its $\rho$-blocks.
For $\bd = (d_1,\ldots,d_n) \in \N^n$ and $Z \in \irr(A,\bd)$ the following are
equivalent:
\begin{itemize}

\item[(i)]
$Z$ is generically reduced;

\item[(ii)]
For each loop $a \in Q_1$, the number $d_{s(a)}$ is even.

\end{itemize}
\end{Thm}

\begin{proof}
We know from Corollary~\ref{cor:blocks4} that
$Z$ is generically reduced if and only if $\pi_i(Z)$ is
generically reduced for all $1 \le i \le t$.
Now the result follows from
Proposition~\ref{prop:complexreduced}.
\end{proof}

\begin{Cor}\label{cor:genreduced2}
Let $A$ be a gentle algebra without loops.
Then each $Z \in \irr(A)$ is generically reduced.
\end{Cor}

Note that Corollary~\ref{cor:genreduced2} is exactly the statement
of Theorem~\ref{thm:main2}.

\subsection{Singular locus}\label{subsec:singular}
The following theorem describes the singular locus of schemes of modules over gentle algebras.
It turns out that the rank function of a module determines
completely if this module is singular or not.

\begin{Thm}\label{thm:main1b}
Let $A = KQ/I$ be a gentle algebra.
Let $M \in \md(A,\bd)$, and let $r = r_M\df Q_1 \to Q_0$ 
be the rank function of $M$.
The following are equivalent:
\begin{itemize}

\item[(i)]
$M$ is singular;

\item[(ii)]
There exist $a,b \in Q_1$ with $s(a) = t(b)$ and $ab \in I$
such that the following hold:
\begin{itemize}
\item[(1)]
$r(a) < d_{t(a)}$, $r(b) < d_{s(b)}$ and
$r(a) + r(b) < d_{s(a)}$.
\item[(2)]
If $a' \in Q_1$ with $s(a') = t(a)$ and $a'a \in I$, then
$r(a') + r(a) < d_{t(a)}$.
\item[(3)]
If $b' \in Q_1$ with $t(b') = s(b)$ and $bb' \in I$, then
$r(b) + r(b') < d_{s(b)}$.
\end{itemize}

\end{itemize}
\end{Thm}

$$
\xymatrix{
\ar@{-->}[r]^{b'} & s(b) \ar[r]^b & s(a) \ar[r]^a & 
t(a) \ar@{-->}[r]^{a'} & 
}
$$

\begin{proof}
Let $A_1,\ldots,A_t$ be the $\rho$-blocks of $A$.
For $M \in \md(A,\bd)$ we know from Corollary~\ref{cor:blocks2} that
$M$ is smooth if and only if $\pi_i(M)$ is smooth for all $1 \le i \le t$.
Now for each $\rho$-block $A_i$ and each dimension vector $\bd$
there is an algebra $A_i' = C_{n_i'}$ or $A_i' = \widetilde{C}_{n_i'}$ 
and an isomorphism
$$
f_{s,\pi_i(\bd)}\df \md(A_i,\pi_i(\bd)) \to \md(A_i',\pi_i(\bd)').
$$
of affine schemes.
In particular, $\pi_i(M)$ is singular if and only if
$f_{s,\pi_i(\bd)}(\pi_i(M))$ is singular.

By Proposition~\ref{prop:complexsingular} we know all
singular points of $\md(A_i',\pi_i(\bd)')$.
The conditions Theorem~\ref{thm:main1b}(ii) and Proposition~\ref{prop:complexsingular}(ii) are equivalent.
More precisely, let $A_i$ be the $\rho$-block containing the arrows $a$ and $b$.
Then $f_{i,\pi_i(\bd)}(\pi_i(M))$ has a direct summand isomorphic
to $S_{ab}$ if and only if condition Theorem~\ref{thm:main1b}(ii)
holds.
This finishes the proof.
\end{proof}

\begin{Thm}\label{thm:singularlocus}
Let $A$ be a gentle Jacobian algebra.
For $M \in \md(A,\bd)$ the following are equivalent:
\begin{itemize}

\item[(i)]
$M$ is singular;

\item[(ii)]
$M$ is contained in at least two different irreducible components
of $\md(A,\bd)$.

\end{itemize}
\end{Thm}

\begin{proof}
Let $A_1,\ldots,A_t$ be the $\rho$-blocks of $A$.
We know that $M$ is singular if and only if $\pi_i(M)$ is
singular for some $1 \le i \le t$.

We also know that $M$ is contained in two different components
if and only if $\pi_i(M)$ is contained in two different components.

Now the claim follows from Lemma~\ref{lem:3blocksingular}.
\end{proof}

\begin{Cor}\label{cor:smoothgentlesurface}
Let $A$ be a gentle Jacobian algebra.
For each $\bd$ we have
$$
\smooth(A,\bd) = \bigcup_{Z \in \irr(A,\bd)} Z^\circ.
$$
\end{Cor}

Note that Corollary~\ref{cor:smoothgentlesurface} is 
exactly the statement of Theorem~\ref{thm:main1}.

\subsection{Band components}\label{subsec:bandcomp}

\begin{Prop}\label{prop:bandcomp1b}
Let $A$ be a gentle algebra, and let $M \in \md(A,\bd)$
be a direct sum of band modules.
Then $M$ is smooth.
\end{Prop}

\begin{proof}
By Lemma~\ref{lem:ARband}(i) we have 
$\pdim(M) \le 1$. This implies $\Ext_A^2(M,M) = 0$.
Now Proposition~\ref{prop:rigidsmooth2} yields that $M$ is 
smooth.
\end{proof}

\begin{Cor}\label{cor:bandcomp3}
Let $A$ be a gentle algebra, and let $Z \in \irr(A)$ be a direct
sum of band components.
Then $Z$ is generically reduced.
\end{Cor}

\begin{proof}
In a direct sum of band components, the direct sums of band modules
form a dense open subset.
Now the statement follows from Proposition~\ref{prop:bandcomp1b}
combined with Lemma~\ref{lem:smoothreduced}.
\end{proof}

\begin{Prop}\label{prop:bandcomp4}
Let $A$ be a gentle algebra.
For any band component $Z \in \irr(A)$ we have
$$
c_A(Z) = e_A(Z) = h_A(Z) = 1.
$$
In particular, $Z$ is a brick component.
\end{Prop}

\begin{proof}
Let $Z$ be a band component.
Thus there is a band $B$ such that the union 
$$
\bigcup_{\lambda \in K^*} \cO_{M(B,\lambda,1)}
$$
forms a dense subset of $Z$.
Let $M = M(B,\lambda,1)$ for some $\lambda \in K^*$.

By Corollary~\ref{cor:comp1} we have $c_A(Z) = 1$.
Now Corollary~\ref{cor:bandcomp3} implies $e_A(Z) = 1$.
In other words, we have
$$
c_A(Z) = e_A(Z) = \dim \Ext_A^1(M,M) = 1.
$$
Now Lemma~\ref{lem:ARband}(ii) together with Corollary~\ref{cor:ARformula2} imply that
$$
\dim \Hom_A(M,\tau_A(M)) = \dim \End_A(M) = 1.
$$
In other words,
$h_A(Z) = 1$ and $M$ is a brick.
It follows that $Z$ is a brick component.
\end{proof}

Note that Proposition~\ref{prop:bandcomp4} yields Theorem~\ref{thm:main6}.

\begin{Cor}\label{cor:bandcomp5}
Let $A$ be a gentle algebra, and let
$Z \in \irr(A)$ be a direct sum of band components.
Then $Z$ is generically $\tau$-reduced.
\end{Cor}

\begin{proof}
We have 
$$
Z = \overline{Z_1 \oplus \cdots \oplus Z_m}
$$ 
for some band components $Z_i = Z(B_i)$, $1 \le i \le m$.
By Lemma~\ref{lem:genericdecomp} we have $c_A(Z) = m$.

By Theorem~\ref{thm:CBS} we get
${\rm ext}_A^1(Z_i,Z_j) = 0$ for all $i \not= j$.
Let 
$$
M = M(B_1,\lambda_1,1) \oplus \cdots \oplus 
M(B_t,\lambda_m,1)
$$ 
with pairwise different $\lambda_1,\ldots,\lambda_m$.
In other words, $M$ is generic in $Z$.
For brevity we set $M_i := M(B_i,\lambda_i,1)$.
It follows that
$$
0 = \dim \Ext_A^1(M_i,M_j) =
\dim \overline{\Hom}_A(M_j,\tau_A(M_i)) 
= \dim \Hom_A(M_j,\tau_A(M_i)) 
$$
for all $i \not= j$.
For the last equality we used again 
Corollary~\ref{cor:ARformula2}.
By Proposition~\ref{prop:bandcomp4} we have 
$$
h_A(Z_i) = 
\dim \Hom_A(M_i,\tau_A(M_i)) = 1
$$
for all $i$.
Combining this, we get 
$$
h_A(Z) = \dim \Hom_A(M,\tau_A(M)) = m.
$$
In other words, $c_A(Z) = h_A(Z)$, thus $Z$ is generically $\tau$-reduced.
\end{proof}

\begin{Thm}\label{thm:bandcomp6}
Let $A$ be a gentle algebra.
For $Z \in \irr(A,\bd)$ the following are equivalent:
\begin{itemize}

\item[(i)]
$Z$ is a direct sum of band components. 

\item[(ii)]
$\dim(Z) = \dim(\GL_\bd(K))$.

\end{itemize}
\end{Thm}

\begin{proof}
(i) $\implies$ (ii):
Let 
$$
Z = \overline{Z(B_1) \oplus \cdots \oplus Z(B_m)}
$$
be a direct sum of band components, and let
$$
M = M(B_1,\lambda_1,1) \oplus \cdots \oplus M(B_m,\lambda_m,1)
$$
be generic in $Z$.
It follows from Proposition~\ref{prop:bandcomp4} and the proof of Corollary~\ref{cor:bandcomp5} that
$$
\dim \End_A(M) = \dim \Ext_A^1(M,M) = m.
$$
By Proposition~\ref{prop:bandcomp1b}, $M$ is smooth.
Thus 
\begin{align*}
\dim(Z) &= \dim \cO_M + \dim \Ext_A^1(M,M) 
\\
&= 
\dim(\GL_\bd(K))-  \dim \End_A(M) + \dim \Ext_A^1(M,M)
\\
&= \dim(\GL_\bd(K)).
\end{align*}

(ii) $\implies$ (i):
Let 
$$
Z = \overline{Z(C_1) \oplus \cdots \oplus Z(C_p) \oplus Z(B_1) \oplus \cdots
\oplus Z(B_q)}
$$
be a direct sum of string and band components.
For a generic $M \in Z$ we get
$c_A(Z) = q$, see Lemma~\ref{lem:genericdecomp}.
In other words,
$$
\dim(Z) = q + \dim \cO_M = \dim(\GL_\bd(K)) - \dim \End_A(M) + q.
$$
Clearly $\dim \End_A(M) \ge p+q$.
So $\dim(Z) = \dim(\GL_\bd(K))$ implies $p=0$.
In other words, $Z$ is a direct sum of band components.
This finishes the proof.
\end{proof}

Combining Corollary~\ref{cor:bandcomp5} and Theorem~\ref{thm:bandcomp6} proves Theorem~\ref{thm:main5}.

The following theorem is a combination of 
\cite[Corollary~10]{CC} and \cite[Proposition~11]{CC},
see also \cite[Theorem~2]{C}.
Proposition~\ref{prop:bandcomp4} generalizes 
Theorem~\ref{thm:CC1}(ii) to arbitrary gentle algebras, whereas
Theorem~\ref{thm:CC1}(i) fails in general.

\begin{Thm}[{\cite{CC}}]\label{thm:CC1}
Let $A$ be an acyclic gentle algebra.
Then the following hold:
\begin{itemize}

\item[(i)]
For each dimension vector $\bd$ there exists at most one
band component $Z$ in $\irr(A,\bd)$.

\item[(ii)]
Each band component $Z \in \irr(A,\bd)$ is a brick component.

\end{itemize}
\end{Thm}

For acyclic gentle algebras $A$,
a combinatorial construction of generic modules for each irreducible
component of $\md(A,\bd)$ is described in \cite{C}.

\subsection{Examples}

\subsubsection{}
Let $A = KQ/I$, where $Q$ is the quiver
$$
\xymatrix{
1 \ar@(ul,ur)^a
}
$$
and $I$ is generated by $\{ a^2 \}$.
Obviously, $A$ is gentle.
Let $\bd = (1)$.
Then $\md(A,\bd)$ has just one $K$-rational point, corresponding
to the simple $A$-module $M = S_1$.
Clearly, $M$ is not smooth and not reduced.

\subsubsection{}
Let $A = KQ/I$, where $Q$ is the quiver
$$
\xymatrix{
1 & 2 \ar[l]_a & 3 \ar[l]_b
}
$$
and $I$ is generated by $\{ ab \}$.
Clearly, $A$ is a gentle algebra.
Let $\bd = (1,1,1)$.
Then $\md(A,\bd)$ has $2$ irreducible components.
The module $M = S_1 \oplus S_2 \oplus S_3 \in \md(A,\bd)$ is reduced, but not smooth.
For $\bd = (1,2,1)$, the affine scheme $\md(A,\bd)$ is
irreducible, reduced, but not smooth.

\subsubsection{}
Let $A = KQ/I$, where $Q$ is the quiver
$$
\xymatrix{
1 \ar@/_/[r]_b & 2 \ar@/_/[l]_a
}
$$
and $I$ is generated by $\{ ab \}$.
Then $A$ is a gentle algebra, which does not admit a colouring in
the sense of \cite{CC}.

\subsubsection{}
Let $A = KQ/I$, where $Q$ is the quiver
$$
\xymatrix{
1 \ar@(ul,dl)_a & 2 \ar[l]_b \ar@/_/[r]_{c_2} & 
3 \ar@/_/[l]_{c_1} \ar[r]^d & 4
\ar@(ur,dr)^e
}
$$
and $I$ is generated by $\{ a^2,\; e^2,\; c_1c_2,\; c_2c_1 \}$.
This is a gentle algebra admitting a colouring.
For $\bd = (2,2,2,2)$, the affine scheme $\md(A,\bd)$
has $3$ irreducible components, and all of these are band
components.


\section{Generically $\tau$-reduced components for gentle Jacobian 
algebras}\label{sec:taureduced}


In this section, we concentrate on the description of 
generically $\tau$-reduced components for gentle Jacobian algebras.
Some of this can be generalized to arbitrary gentle algebras.
We leave this endeavor to the reader.

\subsection{Simple summands of restrictions}
Let $A = KQ/I$ be a gentle Jacobian algebra and let $A_1,\ldots,A_t$ be its $\rho$-blocks.
For $a \in Q_0 \cup Q_1$ and $1 \le s \le t$ let
$$
\delta_{a,A_s} :=
\begin{cases}
1 & \text{if $a$ belongs to $A_s$},
\\
0 & \text{otherwise}.
\end{cases}
$$

\begin{Lem}\label{lem:stringprojection}
Let $A = KQ/I$ be a gentle Jacobian algebra and let $A_1,\ldots,A_t$ be its $\rho$-blocks.
For a string module $M = M(C) \in \md(A)$ and any $\rho$-block $A_s$,
the $A_s$-module $\pi_s(M)$ has a simple direct summand if and only if
one of the following hold:
\begin{itemize}

\item
$C = 1_i$ and $i \in A_s$;

\item 
$C = (c_1,\ldots,c_r)$ with $s(C) \in A_s$ and
$c_r \notin A_s$;

\item 
$C = (c_1,\ldots,c_r)$ with $t(C) \in A_s$ and
$c_1 \notin A_s$.
 
\end{itemize}
\end{Lem}

\begin{proof}
For $C = 1_i$ the claim is clear.
Thus let $C = (c_1,\ldots,c_r)$.
For each $1 \le i \le r$ we have
$c_i = a_i^\pm$ for some $a_i \in Q_1$.
For $M = M(C)$ we get
$$
\pi_s(M) \cong
S_{s(c_r)}^{\delta_{s(c_r),A_s}(1-\delta_{a_r,A_s})} 
\oplus
S_{t(c_1)}^{\delta_{t(c_1),A_s}(1-\delta_{a_1,A_s})}
\oplus
\bigoplus_{i=1}^r P_{s(a_i)}^{\delta_{a_i,A_s}}.
$$
This follows directly from the definition of a string module.
The claim follows.
\end{proof}

\begin{Lem}\label{lem:bandprojection}
Let $A = KQ/I$ be a gentle Jacobian algebra, and let $A_1,\ldots,A_t$ 
be its $\rho$-blocks.
For a band module $M \in \md(A)$ and any $\rho$-block $A_s$,
the module $\pi_s(M)$ has no simple direct summand. 
In particular, $\pi_s(M)$ is a projective $A_s$-module.
\end{Lem}

\begin{proof}
Let $M = M(B,\lambda,q)$ be a band module where
$B = (c_1,\ldots,c_r)$.
For each $1 \le i \le r$ we have
$c_i = a_i^\pm$ for some $a_i \in Q_1$.
We get 
$$
\pi_s(M) \cong \bigoplus_{i=1}^r P_{s(a_i)}^{q\delta_{a_i,A_s}}.
$$
This follows directly from the definition of a band module.
The claim follows.
\end{proof}

\subsection{Non-vanishing of $\Hom_A(M,\tau_A(M))$}
Let $A$ be a gentle Jacobian algebra, and let $A_1,\ldots,A_t$
be its $\rho$-blocks.
Recall from Section~\ref{subsec:Rigid} the definition of critical summands of type $I$ or $II$
for modules over $C_n$ or $\widetilde{C}_n$.
We say that $M \in \md(A,\bd)$ has a \emph{critical summand}
of type $I$ (resp. type $II$) if there exists some $1 \le i \le t$ such
that $\pi_i(M)$ has a critical summand of type $I$ (resp. of type $II$).

\begin{Lem}\label{lem:nonvanish1}
Let $A$ be a gentle Jacobian algebra, and let $A_1,\ldots,A_t$ be 
its $\rho$-blocks.
For $M \in \md(A)$ the following are equivalent:
\begin{itemize}

\item[(i)]
$M$ does not have a critical summand of type $I$ or $II$.

\item[(ii)]
$\pi_i(M)$ is $\tau$-rigid for all $1 \le i \le t$.

\end{itemize}
\end{Lem}

\begin{proof}
This follows from Propositions~\ref{prop:complexrigid} 
and \ref{prop:complextaurigid}
\end{proof}

\begin{Lem}\label{lem:nonvanish2a}
Let $A$ be a gentle Jacobian algebra, and let $A_1,\ldots,A_t$
be its $\rho$-blocks.
Let $M_1,M_2 \in \md(A)$ such that the following hold: 
There exists a $\rho$-block $A_i$ containing an arrow  
$a \in Q_1$ such that
$S_{s(a)}$ is (up to isomorphism) a direct summand of $\pi_i(M_1)$,
and $S_{t(a)}$ is (up to isomorphism) a direct summand of $\pi_i(M_2)$.
Then $\Ext_A^1(M_1,M_2) \not= 0$.
\end{Lem}

\begin{proof}
We can assume that $M_1$ and $M_2$ are both indecomposable.
By Lemma~\ref{lem:bandprojection} we know that
$M_1 = M(C_1)$ and $M_2 = M(C_2)$ are both string modules.
By Lemma~\ref{lem:stringprojection} we can assume without
loss of generality that $s(C_1) = s(a)$ and $t(C_2) = t(a)$
and that $C_1a^{-1}C_2$ is a string.
We obtain a non-split short exact sequence
$$
0 \to M(C_2) \to M(C_1a^{-1}C_2) \to M(C_1) \to 0.
$$
Thus $\Ext_A^1(M_1,M_2) \not= 0$.
\end{proof}

\begin{Cor}\label{cor:nonvanish2a}
Let $A$ be a gentle Jacobian algebra, and let $A_1,\ldots,A_t$
be its $\rho$-blocks.
Let $M_1,M_2 \in \md(A)$ such that $\Ext_A^1(M_1,M_2) = 0$.
Then
$$
\Ext_A^1(\pi_i(M_1),\pi_i(M_2)) = 0
$$ 
for all $1 \le i \le t$.
\end{Cor}

\begin{proof}
Combine Lemma~\ref{lem:nonvanish2a} with
Proposition~\ref{prop:complexrigid}.
\end{proof}

\begin{Lem}\label{lem:nonvanish2b}
Let $A$ be a gentle Jacobian algebra, and let $A_1,\ldots,A_t$
be its $\rho$-blocks.
Let $M_1,M_2 \in \md(A)$ such that
the following hold: 
There exists a $3$-block $A_s$ containing an arrow $a \in Q_1$
such that
$S_{s(a)}$ is (up to isomorphism) a direct summand of $\pi_s(M_1)$,
and $P_{t(a)}$ is (up to isomorphism) a direct summand of $\pi_s(M_2)$.
Then
$$
\dim \Hom_A(M_2,\tau_A(M_1)) \not= 0.
$$
\end{Lem}

\begin{proof}
We can assume that $M_1$ and $M_2$ are both indecomposable.
We know that $M_1 = M(C_1)$ for some string $C_1$ 
(see Lemma~\ref{lem:bandprojection}) and
$M_2 = M(C_2)$ or $M_2 = M(C_2,\lambda,q)$ for some
string or band $C_2$, respectively.

If $M_2$ is a band module, then there is a surjective
homomorphism
$M(C_1,\lambda,q) \to M(C_1,\lambda,1)$.
Thus in this case we can assume without loss of generality that 
$q=1$.

We can assume that the
$3$-block $A_s$ is of the form
$$
\xymatrix@-1.7pc{
&&1 \ar[dddll]_{a_1}\\
&&\\
&&\\
2 \ar[rrrr]_{a_2} &&&& 3 \ar[lluuu]_{a_3}
}
$$
with $a_2a_1, a_3a_2, a_1a_3 \in I$ and $a = a_1$.

Without loss of generality we can assume that $s(C_1) = 1$ and that
either $l(C_1) = 0$ or
$C_1 = (c_1,\ldots,c_m)$ such that $c_m \notin A_2$.
We can also assume that
$C_2 = C'a_2C''$ for some strings $C'$ and $C''$ and we can
assume that $C' = (c_1',\ldots,c_r')$ with $c_1' \in Q_1^{-1}$.

We want to construct a non-zero homomorphism 
$$
M_2 \to \tau_A(M_1).
$$

Let $E$ be a path of maximal length such that $a_1^{-1}E$ is a string.
It follows that $\tau_A(M(C_1)) = M(E'E)$ for some string $E'$, where
$E'$ is either of length $0$ or of the form $E' = E''a_1^{-1}$ for
some string $E''$,
compare Section~\ref{subsec:taustring}.

Let $F$ be a path of maximal length such that $FF' = C''$.
Thus $C_2 = C'a_2FF'$.
It follows that $F'$ is of length $0$ or of the form $F' = b^{-1}F''$ for some $b \in Q_1$ and some string $F''$.
This yields a surjective homomorphism
$$
f_1\df M_2 \to M(F).
$$
Furthermore, we have $E = FG'$ for some direct string $G'$.
We get a standard homomorphism
$$
f_2 = f_{(1_{t(F)},F,1_{s(F)}),(E',F,G')}\df M(F) \to \tau_A(M(C_1)).
$$
Thus
$$
f_2 \circ f_1\df M_2 \to \tau_A(M_1)
$$
is the desired non-zero homomorphism. 
This finishes the proof.
\end{proof}

\begin{Cor}\label{cor:nonvanish2c}
Let $A$ be a gentle Jacobian algebra, and 
let $A_1,\ldots,A_t$ be its $\rho$-blocks.
Let $M_1,M_2 \in \md(A)$ such that 
$\Hom_A(M_2,\tau_A(M_1)) = 0$.
Then
$$
\Hom_{A_i}(\pi_i(M_2),\tau_{A_i}(\pi_i(M_1))) = 0
$$
for all $1 \le i \le t$.
\end{Cor}

\begin{proof}
Combine 
Lemma~\ref{lem:nonvanish2b}, Corollary~\ref{cor:nonvanish2a}
and Proposition~\ref{prop:complextaurigid}.
\end{proof}

\begin{Cor}\label{cor:nonvanish2d}
Let $A$ be a gentle Jacobian algebra, and assume that
$M \in \md(A)$ has a critical summand of type $I$ or $II$.
Then
$$
\dim \Hom_A(M,\tau_A(M)) \not= 0.
$$
\end{Cor}

\subsection{Proof of Theorem~\ref{thm:main4}}
Let $A = KQ/I$ be a gentle Jacobian algebra, and
let $A_1,\ldots,A_t$ be its $\rho$-blocks.
Let $Z \in \irr(A)$.
We want to show that the following
are equivalent:
\begin{itemize}

\item[(i)]
$Z \in \irr^\tau(A)$;

\item[(ii)]
$\pi_s(Z) \in \irr^\tau(A_s)$ for all $1 \le s \le t$.

\end{itemize}

Throughout, let
$$
M = M(C_1) \oplus \cdots \oplus M(C_p) \oplus
M(B_1,\lambda_1,1) \oplus \cdots \oplus M(B_q,\lambda_q,1)
$$
be generic in $Z$.

For $1 \le i \le p$ let $N_i := M(C_i)$, and for
$1 \le j \le q$ let $N_{p+j} := M(B_j,\lambda_j,1)$.  

(i) $\implies$ (ii):
Assume that $Z \in \irr^\sr(A)$.
Then 
Theorem~\ref{thm:decomp} yields that
$\Hom_A(N_i,\tau_A(N_j)) = 0$ for all $i \not= j$.
Furthermore we have
$$
\dim \Hom_A(N_i,\tau_A(N_i)) =
\begin{cases}
0 & \text{if $1 \le i \le p$},
\\
1 & \text{if $p+1 \le i \le p+q$}.
\end{cases}
$$
Now it follows from Corollary~\ref{cor:nonvanish2c} that
$$
\Hom_{A_s}(\pi_s(N_i),\tau_{A_s}(\pi_s(N_j))) = 0
$$
for all $i \not = j$, and also for all $i=j$ with $1 \le i \le p$.
Since $N_{p+1},\ldots,N_{p+q}$ are band modules,
we get from Lemma~\ref{lem:bandprojection} that also in this case
$$
\Hom_{A_s}(\pi_s(N_i),\tau_{A_s}(\pi_s(N_i))) = 0.
$$
This proves that $\pi_s(M)$ is a $\tau$-rigid $A_s$-module
for all $s$.
Thus $\pi_s(Z) \in \irr^\tau(A_s)$.

(ii) $\implies$ (i):
Assume that $\pi_s(Z) \in \irr^\sr(A_s)$ for
all $1 \le s \le t$.

We have
$$
0 = \Ext_A^1(N_{p+j},N_k) \cong 
\overline{\Hom}_A(N_k,\tau_A(N_{p+j})) 
= \Hom_A(N_k,\tau_A(N_{p+j})) 
$$
for all $1 \le j \le q$ and $1 \le k \le p+q$.
For the third equality we used Corollary~\ref{cor:ARformula2}.

Thus $Z$ is generically $\tau$-reduced if and only if
$\Hom_A(N_k,\tau_A(N_i)) = 0$ for all $1 \le i \le p$ and
$1 \le k \le p+q$.
To get a contradiction, assume that 
$$
\Hom_A(N_k,\tau_A(N_i)) \not= 0
$$
for some $1 \le i \le p$ and some $1 \le k \le p+q$.
On the other hand, we know that
$$
0 = \Ext_A^1(N_i,N_k) \cong \overline{\Hom}_A(N_k,\tau_A(N_i)).
$$

Let 
$$
f\colon N_k \to \tau_A(N_i)
$$ 
be a non-zero homomorphism. 
We know that $f$ factors through some injective
$A$-module.
Without loss of generality, we can assume that
this injective module equals $I_r$ for some $r \in Q_0$.
Thus we have $f = f_1 \circ f_2$ with
$f_1 \in \Hom_A(I_r,\tau_A(N_i))$ and $f_2 \in \Hom_A(N_k,I_r)$.
Again without loss of generality we can assume that
$$
f_1 = f_{(E,F,G),(E',F,G')}\df I_r \to \tau_A(N_i)
$$ 
is a standard homomorphism.
(Here we use the same notation and terminology as in \cite{Sch}.)

The module $I_r$ is of the form
$$
I_r = M(D^{-1}C),
$$
where $C$ and $D$ are direct strings in $Q$ 
such that $C\gamma,D\gamma \in I$ for all $\gamma \in Q_1$.
Since $\tau_A(M_i)$ is not injective, we know that $f_1$ cannot be a
monomorphism.
Thus $I_r$ is not simple and we can assume without loss of
generality that
$C = \alpha_1 \cdots \alpha_v$
and that 
$$
\Ker(f_1) = M(D^{-1}\alpha_1 \ldots \alpha_k)
$$
for some $1 \le k \le v$.
Thus we have
\begin{align*}
E &= D^{-1}\alpha_1 \cdots \alpha_k,
\\
F &= \alpha_{k+1} \cdots \alpha_v & \text{if $k < v$}, 
\\ 
F &= 1_{s(C)} & \text{if $k = v$},
\\
G &= 1_{s(C)}.
\end{align*}
Since $A$ is a gentle algebra, we also know that 
$\alpha_v\gamma \in I$ for all $\gamma \in Q_1$.
This implies $E' = 1_{s(\alpha)}$.
Set $a := \alpha_k$.

The following picture shows $I_r$, where $I_r/\Ker(f_1)$ is given
by the string $F$ between the blue vertices.
$$
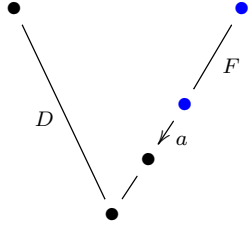
\xymatrix@-4.5ex{
\bullet \ar@{-}[rrrrdddddddd]_D
&&&&&&&& \blue\bullet \ar@{-}[lldddd]^F
\\
&&\\
&&\\
&&\\
&&&&&&\blue\bullet \ar[ldd]^a
\\
&&&&\\
&&&&&\bullet\ar@{-}[ldd]
\\
&\\
&&&& \bullet
}
$$

By the properties of $f_1$ discussed above, we see that we must
be in the 2nd, 4th or $5$th case and that $F$ coincides with
the subfactor of $\tau_A(M)$ marked by the two rightmost 
blue points. 
Here we refer to Section~\ref{subsec:taustring} for the description
of $\tau_A(M)$.

We get $ab \in I$.
Thus there exists a third arrow $c \in Q_1$ with $s(c) = t(a)$
and $t(c) = s(b)$.
So the arrows $a,b,c$ form a $3$-block, say $A_s$, of $A$.
So we are in the following situation:
$$
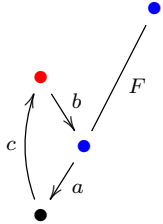
\xymatrix@-4.0ex{
&&& \blue\bullet \ar@{-}[lldddd]^F
\\
&&\\
\red\bullet\ar[rdd]^b&&\\
&&\\
&\blue\bullet \ar[ldd]^a
\\
&&&&\\
\bullet\ar@/^/[uuuu]^c
}
$$
(In the 5th case, the red bullet in this picture should be green.)

Clearly $\pi_s(N_k)$ contains $M(a)$ as a direct summand,
and $\pi_s(N_i)$ contains $S_{s(b)}$ as a direct summand.

It follows that $\pi_s(N_k \oplus N_i)$ has a direct summand
isomorphic to $S_{s(b)} \oplus P_{s(a)}$.
Now Proposition~\ref{prop:complextaurigid} implies that
$\pi_s(N_k \oplus N_i)$ and therefore also $\pi_s(M)$ 
is not $\tau$-rigid in $\md(A_s)$.
This finishes the proof.

\subsection{Proof of Theorem~\ref{thm:main3}}
Let $A$ be a gentle Jacobian algebra, and let $A_1,\ldots,A_t$
be its $\rho$-blocks.
Let $Z_1,Z_2 \in \irr^\tau(A)$.
We want to show that the following are
equivalent:
\begin{itemize}

\item[(i)]
$\dimv(Z_1) = \dimv(Z_2)$;

\item[(ii)]
$Z_1 = Z_2$.

\end{itemize}

(ii) $\implies$ (i):
This direction is trivial.

(i) $\implies$ (ii):
Assume that $\dimv(Z_1) = \dimv(Z_2)$.
We know from Theorem~\ref{thm:main4} that
$\pi_i(Z_1)$ and $\pi_i(Z_2)$ are generically $\tau$-reduced for all $1 \le i \le t$.
In particular, $\pi_i(Z_1)$ and $\pi_i(Z_2)$ both contain a $\tau$-rigid
$A_i$-module.
We clearly have $\dimv(\pi_i(Z_1)) = \dimv(\pi_i(Z_2))$ for 
all $i$.
Note that for gentle Jacobian algebras, we have $A_i = A_i'$ for all
$i$.
Now the statement follows from Lemma~\ref{lem:3blockunique}.


\section{Schemes of decorated modules}\label{sec:decrep}


\subsection{Decorated modules}
Let $A = KQ/I$ be a basic algebra.
A \emph{decorated $A$-module} 
is a pair $\cM = (M,V)$, where $M \in \md(A)$ and 
$V = (V_1,\ldots,V_n)$ is
a tuple of finite-dimensional $K$-vector spaces.

One defines morphisms and direct sums of decorated modules in
the obvious way.
Let $\decrep(A)$ be the abelian category of decorated $A$-modules.

For $1 \le i \le n$ set $\cS_i := (S_i,0)$, and let
$\cS_i^- := (0,V)$, where
$V_i = K$ and $V_j = 0$ for all $j \not= i$.
The decorated modules $\cS_i$ and $\cS_i^-$ are the \emph{simple}
and \emph{negative simple} 
decorated $A$-modules, respectively.

\subsection{Schemes of decorated modules}
For $(\bd,\bv) \in \N^n \times \N^n$ let
$\decrep(A,(\bd,\bv))$ be the affine scheme
of decorated $A$-modules $\cM = (M,V)$
with $M \in \md(A,\bd)$ and
$V = K^\bv := (K^{v_1},\ldots,K^{v_n})$,
where $\bv = (v_1,\ldots,v_n)$.
Note that $\md(A,\bd) \cong \decrep(A,(\bd,\bv))$ for all
$(\bd,\bv)$.

For $\cM = (M,V) \in \decrep(A,(\bd,\bv))$ let
$g.\cM := (g.M,V)$.
This defines a $\GL_\bd(K)$-action on $\decrep(A,(\bd,\bv))$.
The $\GL_\bd(K)$-orbit of $\cM$ is denoted by $\cO_\cM$.

\subsection{$E$-invariants and $g$-vectors of decorated modules}
Let $\cM = (M,V)$ be a decorated $A$-module, and let
$$
\bigoplus_{i=1}^n P_i^{m_i} \to \bigoplus_{i=1}^n P_i^{n_i}
\to M \to 0
$$
be a minimal projective presentation of $M$.
The $g$-\emph{vector} of $\cM$ is defined as
$$
\bg_A(\cM) := (g_1,\ldots,g_n)
$$
with 
$$
g_i := g_i(\cM) := m_i-n_i + \dim(V_i)
$$
for $1 \le i \le n$.

For decorated $A$-modules $\cM = (M,V)$ and
$\cN = (N,W)$ let
$$
E_A(\cM,\cN) := \dim \Hom_A(N,\tau_A(M)) + \sum_{i=1}^n
\dim(V_i)\dim(N_i).
$$
For finite-dimensional Jacobian algebras $A$ arising from quivers with potentials,
$E_A(\cM,\cN)$ coincides with $E^\proj(\cM,\cN)$ 
as defined in \cite[Section~10]{DWZ2}.

The $E$-\emph{invariant} of $\cM$ is defined as
$E_A(\cM) := E_A(\cM,\cM)$.
The decorated module $\cM$ is called $E$-\emph{rigid} if $E_A(\cM) = 0$.

For $\cM = (M,0)$ we also write $\bg_A(M)$ and $E_A(M)$ instead of 
$\bg_A(\cM)$ and $E_A(\cM)$, respectively.

Dualizing the results from \cite[Section~10]{DWZ2} (for Jacobian algebras $A$) and \cite[Section~3]{CLFS} (for arbitrary $A$), for 
decorated $A$-modules $\cM = (M,V)$ and $\cN = (N,W)$ 
we have
$$
E_A(\cM,\cN) = \dim \Hom_{A}(M,N) 
+  \sum_{i=1}^n g_i(\cM)\dim(N_i).
$$
Note that in \cite{DWZ2} and
\cite{CLFS}, this equation is used as a definition.

\subsection{Generically $\tau$-reduced decorated components}
Let $A = KQ/I$ be a basic algebra, and let
$(\bd,\bv) \in \N^n \times \N^n$.
By $\decirr(A,(\bd,\bv))$
we denote the set of irreducible components
of $\decrep(A,(\bd,\bv))$.
For $Z \in \decirr(A,(\bd,\bv))$ we write 
$\dimv(Z) := (\bd,\bv)$.
Let 
$$
\decirr(A) := \bigcup_{(\bd,\bv) \in \N^n \times \N^n} 
\decirr(A,(\bd,\bv)).
$$

For $Z \in \decirr(A,(\bd,\bv))$ 
set $Z' := \{ M \in \md(A,\bd) \mid (M,K^\bv) \in Z \}$.
We clearly have $Z' \in \irr(A,\bd)$, and write
$Z = (Z',K^\bv)$. 
Define $c_A(Z) := c_A(Z')$ and $e_A(Z) := e_A(Z')$.

For $Z,Z_1,Z_2 \in \decirr(A)$ there are dense open subsets
$U \subseteq Z$ and $U' \subseteq Z_1 \times Z_2$ such that
the maps $\bg_A(-)$, $E_A(-)$ and $E_A(-,-)$ are constant on $U$ and
$U'$, respectively.
These generic values are denoted by $\bg_A(Z)$, $E_A(Z)$ and $E_A(Z_1,Z_2)$,
respectively.

For $Z \in \decirr(A)$ we have
$$
c_A(Z) \le e_A(Z) \le E_A(Z).
$$
An irreducible component $Z \in \decirr(A)$ is
\emph{generically reduced} if $c_A(Z) = e_A(Z)$ and
\emph{generically $\tau$-reduced} provided
$$
c_A(Z) = e_A(Z) = E_A(Z).
$$
Let $\decirr^\sr(A,(\bd,\bv))$ be the set of all generically $\tau$-reduced components of $\decrep(A,(\bd,\bv))$, and let
$$
\decirr^\sr(A) := \bigcup_{(\bd,\bv) \in \N^n \times \N^n}
\decirr^\sr(A,(\bd,\bv)).
$$
It follows from the definitions that
$$
\decirr^\sr(A,(\bd,\bv)) = \{ (Z,K^\bv) \mid 
Z \in \irr^\sr(A,\bd),\, d_1v_1 + \cdots + d_nv_n = 0 \},
$$
where $\bd = (d_1,\ldots,d_n)$ and $\bv = (v_1,\ldots,v_n)$.

The following beautiful result due to Plamondon shows that
the generic $g$-vectors parametrize the 
generically $\tau$-reduced decorated components.

\begin{Thm}[{\cite[Theorem~1.2]{P1}}]\label{parametrization}
Let $A$ be a basic algebra.
Then the map
\begin{align*}
\bg_A\df \decirr^\sr(A) &\to \Z^n
\\
Z &\mapsto \bg_A(Z)
\end{align*}
is bijective.
\end{Thm}

\subsection{Decomposition of generically $\tau$-reduced components}\label{subsec:decomp}
An irreducible component $Z \in \decirr(A,(\bd,\bv))$ is called
\emph{indecomposable} if 
there exists a dense open
subset $U \subseteq Z$, which contains only indecomposable decorated modules.
This is the case if and only if
$Z = (Z',0)$ with
$Z' \in \irr(A,\bd)$ indecomposable or $Z = \{ \cS_i^- \}$ for some
$i$.
In particular, if $Z \in \decirr(A,(\bd,\bv))$ is indecomposable,
then either $\bd = 0$ or $\bv = 0$.

Given irreducible components $Z_i$ of 
$\decrep(A,(\bd_i,\bv_i))$
for $1 \le i \le t$, let $(\bd,\bv) := (\bd_1,\bv_1) + \cdots + (\bd_t,\bv_t)$, and let 
$$
Z_1 \oplus \cdots \oplus Z_t
$$
be the image of the morphism
\begin{align*}
\GL_\bd(K) \times Z_1 \times \cdots \times Z_t &\to
\decrep(A,(\bd,\bv))
\\
(g,(\cM_1,\ldots,\cM_t)) &\mapsto g.(\cM_1 \oplus \cdots \oplus \cM_t).
\end{align*}

The Zariski closure $\overline{Z_1 \oplus \cdots \oplus Z_t}$ of
$Z_1 \oplus \cdots \oplus Z_t$ is
an irreducible closed subset of $\decrep(A,(\bd,\bv))$ and is
called the \emph{direct sum} of $Z_1,\ldots,Z_t$.
Note that $\overline{Z_1 \oplus \cdots \oplus Z_t}$ is
in general not an irreducible component.

\begin{Thm}[{\cite[Theorem~1.3]{CLFS}}]\label{thm:decomp}
For $Z_1,\ldots,Z_t \in \decirr(A)$ the following are equivalent:
\begin{itemize}

\item[(i)]
$\overline{Z_1 \oplus \cdots \oplus Z_t}$ is a generically $\tau$-reduced component.

\item[(ii)]
Each $Z_i$ is generically $\tau$-reduced and
$E_A(Z_i,Z_j) = 0$ for all $i \not= j$.

\end{itemize}
\end{Thm}

Each  $Z \in \decirr^\sr(A)$
is a direct sum of indecomposable generically $\tau$-reduced components,
which are uniquely determined up to reordering.


\section{Laminations and generically $\tau$-reduced decorated components}\label{sec:laminations}


\subsection{Marked surfaces}
By an \emph{unpuntured marked surface} $\SM$
we mean a connected, compact, oriented surface $\bbS$ with non-empty boundary $\partial \bbS$ together with a finite set $\bbM$ of \emph{marked points} 
on the boundary.
We set $\bbS^\circ := \bbS \setminus \partial\bbS$.
We assume that there is at least one marked point on each boundary component.
We also require that $\SM$ is not a monogon, digon or triangle.
(This ensures the existence of non-trivial triangulations.)

\subsection{Curves and loops}
A \emph{curve in} $\SM$ is a map 
$$
\gamma\df [0,1] \to \bbS
$$ 
of differentiability class $C^1$, with derivative vanishing in
at most finitely many points of $[0,1]$,
such that the following hold:
\begin{itemize}

\item[(A1)]
$\gamma(0)$ and $\gamma(1)$ are in $\bbM$;

\item[(A2)]
$\Ima(\gamma) \setminus \{ \gamma(0),\gamma(1) \}$ is disjoint from $\partial \bbS$;

\item[(A3)]
$\Ima(\gamma)$ does not cut out a monogon or digon.

\end{itemize}

A curve $\gamma$ in $\SM$ is \emph{simple} if 
additionally the following holds:
\begin{itemize}

\item[(A4)]
$\gamma$ is injective on the open interval $(0,1)$, i.e.
$\gamma$ does not intersect itself, except that $\gamma(0)$ and
$\gamma(1)$ may coincide.

\end{itemize}
Simple curves in $\SM$ are also called \emph{arcs}.

Let $S^1$ be the unit circle in $\C$. 
A \emph{loop in} $\SM$ is a map
$$
\gamma\df S^1 \to \bbS
$$ 
of differentiability class $C^1$, with derivative vanishing in
at most finitely many points of $S^1$,
such that
the following hold:
\begin{itemize}

\item[(L1)]
$\Ima(\gamma)$ is disjoint from $\partial \bbS$;

\item[(L2)]
$\Ima(\gamma)$ is non-contractible.

\end{itemize}
A loop $\gamma$ in $\SM$ is \emph{simple} if additionally
the following holds:
\begin{itemize}

\item[(L3)]
$\gamma$ is injective, i.e.
$\gamma$ does not intersect itself.

\end{itemize}

Let $A\SM$ be the set of curves in $\SM$
up to homotopy (relative to $\gamma(0)$ and $\gamma(1)$) 
and up to the equivalence 
$\gamma \sim \gamma^{-1}$.
For a curve $\gamma$ let $[\gamma]$ be its class in
$A\SM$.

Let $L\SM$ be the set of loops in $\SM$ up to homotopy and up
to the equivalence 
$\gamma \sim \gamma^{-1}$.
For a loop $\gamma$ let $[\gamma]$ be its class
in $L\SM$.

For a curve or loop $\gamma$ in $\SM$ 
we just write $\gamma$ for the image $\Ima(\gamma)$.

For
$$
[\gamma],[\delta] \in A\SM \cup L\SM
$$
let 
$$
\Int([\gamma],[\delta]) := \min\{ 
|\gamma' \cap \delta' \cap \bbS^\circ| 
\mid \gamma' \in [\gamma],\; \delta' \in [\delta] \}. 
$$ 
Note that for a simple curve or loop $\gamma$ we get $\Int([\gamma],[\gamma]) = 0$.

From now on we will not distinguish between a curve
or loop $\gamma$ and its class $[\gamma]$.

A loop $\gamma$ is \emph{primitive} of it is not of the form 
$\gamma = \theta^m$ for some loop $\theta$ and some $m \ge 2$.
Here $\theta^m(z) := \theta(z^m)$ for all $z \in S^1$.
Let
$$
L\SM^{\rm prim} \subset L\SM
$$
be the subset of primitive loops.

Let
$\pi\df \widetilde{S}^1 \to S^1$ be the universal cover of $S^1$.
For a loop $\gamma\df S^1 \to \bbS$ in $\SM$ let 
$$
\widetilde{\gamma} := \gamma \circ \pi \df \widetilde{S}^1 \to \bbS.
$$
We call this the \emph{periodic curve} associated with $\gamma$.

For later convenience, for $\gamma \in A\SM$ we set
$\widetilde{\gamma}:= \gamma$.

For $\gamma_1,\gamma_2 \in A\SM \cup L\SM$, define
$\Int(\widetilde{\gamma}_1,\gamma_2)$ and 
$\Int(\widetilde{\gamma}_1,\widetilde{\gamma}_2)$
similarly as above.
Note that the value $\infty$ might occur in this situation.

The following lemma is straightforward

\begin{Lem}\label{lem:opencurves1}
For $\gamma_1,\gamma_2 \in A\SM \cup L\SM^{\rm prim}$ the following are
equivalent:
\begin{itemize}

\item
$\Int(\gamma_1,\gamma_2) = 0$;

\item
$\Int(\widetilde{\gamma}_1,\gamma_2) = 0$;

\item
$\Int(\widetilde{\gamma}_1,\widetilde{\gamma}_2) = 0$.

\end{itemize}
\end{Lem}

\subsection{Laminations and triangulations} \label{ssec:LamTri}
By a \emph{lamination} of $\SM$ we mean a
pair $L = (\gamma,m)$, where
$\gamma$
is a (finite) subset of $A\SM \cup L\SM$ 
such that $\Int(\gamma_i,\gamma_j) = 0$ for all
$\gamma_i,\gamma_j \in \gamma$, and $m\df \gamma \to \Z_{>0}$
is a map.
Instead of $L = (\gamma,m)$ we also write
$L = \{ (\gamma_1,m_1),\ldots,(\gamma_t,m_t) \}$, where
$\gamma = \{ \gamma_1,\ldots,\gamma_t \}$ and $m_i = m(\gamma_i)$ for $1 \le i \le t$.
By abuse of terminology, we also say that $\gamma$ is a
lamination.
Note that each element in $\gamma$ is a
simple curve or a  
simple loop.
We think of $m_i$ as the multiplicity of $\gamma_i$ in the lamination
$L$.
Let $\Lam\SM$ be the set of laminations of $\SM$.
Note that in \cite[Definition~3.17]{MSW2}, the set $\Lam\SM$ is denoted by $\cC^\circ\SM$.

Each boundary component of $\SM$ with $m$ marked points has
$m$ \emph{boundary segments}, each connecting two consecutive marked points.

\begin{figure}[H]
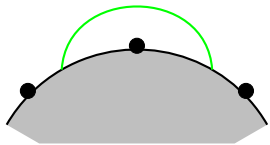
 
  \tikz{
    \fill[fill=lightgray] (30:1.5) -- (30:2)
        arc[start angle=30, end angle=150, radius=2] -- (150:1.5) -- cycle;
    \filldraw (45:2.05) circle[radius=.1];
    \filldraw (90:2.05) circle[radius=.1];
    \filldraw (135:2.05) circle[radius=.1];
    \draw[thick] (30:2) arc[start angle=30, end angle=150, radius=2];
    \draw[thick,green] (60:2) .. controls (72:3) and (108:3) .. (120:2);
  }
\caption{The green curve cuts out a monogon, so it is not a
\emph{curve in} $\bbS \setminus \bbM$.}\label{fig:forbidden}
\end{figure}

Next, a \emph{curve in} $\bbS \setminus \bbM$ is
a map 
$$
\gamma\df [0,1] \to \bbS \setminus \bbM
$$ 
of differentiability class $C^1$, with derivative vanishing in
at most finitely many points of $[0,1]$,
such that the following hold:
\begin{itemize}

\item[(A1)]
$\gamma(0)$ and $\gamma(1)$ are in 
$\partial \bbS \setminus \bbM$;

\item[(A2)]
$\Ima(\gamma) \setminus \{ \gamma(0),\gamma(1) \}$ is disjoint from $\partial \bbS$;

\item[(A3)]
$\Ima(\gamma)$ is non-contractible (with respect to the relative homotopy described below) and does not cut out a monogon, see
Figure~\ref{fig:forbidden}.

\end{itemize}
A curve $\gamma$ in $\bbS \setminus \bbM$ is \emph{simple} 
if additionally the following holds:
\begin{itemize}

\item[(A4)]
$\gamma$ is injective, i.e.
$\gamma$ does not intersect itself.

\end{itemize}

Let $A(\bbS \setminus \bbM)$ be the set of curves in $\bbS \setminus \bbM$
up to homotopy, such that $\gamma(0)$ and $\gamma(1)$ never
leave their respective boundary segment, and up to the equivalence 
$\gamma \sim \gamma^{-1}$.
More precisely, we consider here homotopies
$$
H\df [0,1] \times [0,1] \to \bbS \setminus \bbM
$$
such that for each $t \in [0,1]$ the map $H_t := H(t,-)$ is a
curve in $\bbS \setminus \bbM$ with $H_t(0)$ (resp. $H_t(1)$)
belonging to the same boundary segment as $\gamma(0)$ (resp. $\gamma(1)$) and such that $H_0 = \gamma$.

As before, we just write $\gamma$ for the class of $\gamma$ in
$A(\bbS \setminus \bbM)$.

By a \emph{classical lamination} of $\SM$ we mean a
pair $L = (\gamma,m)$, where
$\gamma$
is a (finite) subset of $A(\bbS \setminus \bbM) \cup L\SM$ 
such that $\Int(\gamma_i,\gamma_j) = 0$ for all
$\gamma_i,\gamma_j \in \gamma$, and $m\df \gamma \to \Z_{>0}$
is a map.
Here $\Int(\gamma_i,\gamma_j)$ is defined in the obvious way.
Again by abuse of terminology, we also say that $\gamma$ is a
classical lamination.
Let $\Lam(\bbS \setminus \bbM)$ be the set of classical laminations of $\SM$. 

Given a curve $\gamma \in A\SM$, let 
$\tau^{1/2}(\gamma) \in A(\bbS \setminus \bbM)$ be the curve
obained from $\gamma$ by rotating its endpoints in clockwise direction to the adjacent boundary segment.
This yields a bijection
$$
\tau^{1/2}\df \Lam\SM \to \Lam(\bbS \setminus \bbM).
$$

A \emph{triangulation} $T$ of $\SM$ consists of all boundary segments
together with a maximal collection $T^\circ$
of curves in $\SM$ such that  
$\Int(\gamma_i,\gamma_j) = 0$ for all $\gamma_i,\gamma_j \in T^\circ$.
In this case, we have
$$
|T^\circ| = 6g+3b+|\bbM|-6,
$$
where $g$ is the genus of $\bbS$ and $b$ is the number of
boundary components of $\bbS$, see for example
\cite[Proposition~2.10]{FST}.

Note that the classical laminations defined above correspond to the
\emph{$\cX$-laminations} in the sense of Fock and Goncharov
\cite{FG}.
Let $T$ be a triangulation of $\SM$, and let $A_T$ be the associated
gentle Jacobian algebra.
We refer to Section~\ref{subsec:9.5} for a precise definition of $A_T$.
To a lamination $L$ of $\SM$ we will associate a certain generic 
decorated $A_T$-module, which is a direct sum
of indecomposable $\tau$-rigid modules, of certain
band modules of quasi-length $1$, and of negative simples.
In Section~\ref{sec:bangle}, we will look at the Caldero-Chapoton
functions of these modules, which can be thought of 
as generating functions of Euler characteristics of
quiver Grassmannians.
In contrast, 
Allegretti \cite{A} works with certain $\cA$-laminations (see
\cite{A,FG} for a definition), and he associates $A_T$-modules,
which are direct sums of indecomposable $\tau$-rigid modules,
of band modules with arbitrary quasi-length, and of negative simples.
He then looks at certain generating functions of Euler characteristics
of transversal quiver Grassmannians.

\subsection{Curves and loops as combinatorial objects}
\label{ssec:curv-comb}
A triangulation cuts the surface into \emph{triangles}.
Each triangle has exactly three sides.
(Recall that we work here with unpunctured marked surfaces, i.e.
we do not have any marked points in the interior of $\bbS$.)

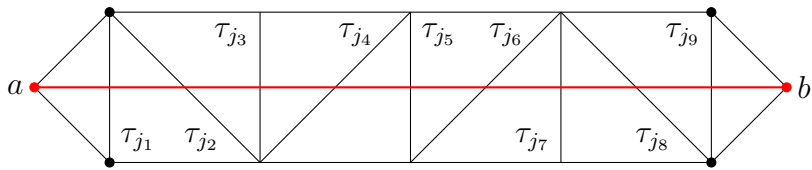
\begin{figure}[H]
\begin{tikzpicture}
\draw[-](1,0) -- (1,2);
\draw[-](1,2) -- (0,1);
\draw[-](0,1) -- (1,0);
\draw [fill,red] (0,1) circle [radius=0.06];
\draw [fill] (1,0) circle [radius=0.06];
\draw [fill] (1,2) circle [radius=0.06];

\draw[-](9,0) -- (9,2);
\draw[-](9,2) -- (10,1);
\draw[-](10,1) -- (9,0);
\draw [fill] (9,0) circle [radius=0.06];
\draw [fill] (9,2) circle [radius=0.06];
\draw [fill,red] (10,1) circle [radius=0.06];

\draw[-](1,2) -- (9,2);
\draw[-](1,0) -- (9,0);

\draw[-](1,2) -- (3,0); \node[right] at (1,0.3) {$\tau_{j_1}$};
\draw[-](3,0) -- (3,2); \node[left] at (2.6,0.3) {$\tau_{j_2}$};
\draw[-](3,0) -- (5,2); \node[left] at (3,1.7) {$\tau_{j_3}$};
\draw[-](5,2) -- (5,0); \node[right] at (3.9,1.7) {$\tau_{j_4}$};
\draw[-](5,0) -- (7,2); \node[right] at (5.9,1.7) {$\tau_{j_6}$};
\draw[-](7,2) -- (7,0); \node[right] at (5.0,1.7) {$\tau_{j_5}$};
\draw[-](7,2) -- (9,0); \node[left] at (7,0.3) {$\tau_{j_7}$};
\node[left] at (8.6,0.3) {$\tau_{j_8}$};
\node[left] at (9,1.7) {$\tau_{j_9}$};

\draw[thick,red,-](0,1) -- (10,1);
\node[left] at (0,1) {$a$};
\node[right] at (10,1) {$b$};

\end{tikzpicture}
\caption{A curve $(a,\tau_{j_1},\ldots,\tau_{j_9},b)$}\label{fig:curve}
\end{figure}

Let $T$ be a fixed triangulation of $\SM$
with $T^\circ = \{ \tau_1,\ldots,\tau_t \}$. 
Let $\gamma\df [0,1] \to \bbS$
be a curve in $\SM$, and let 
$$
m := \Int(\gamma,T) := 
\sum_{\tau \in T^\circ} \Int(\gamma,\tau)
$$
We assume that $\gamma$ is minimal in the sense that
$$
m = \sum_{\tau \in T^\circ} |\gamma \cap \tau \cap \bbS^\circ|.
$$
To $\gamma$
we associate a sequence
$$
(a,\tau_{j_1},\ldots,\tau_{j_m},b),
$$
where $a = \gamma(0)$ and 
$b = \gamma(1)$, and there exist $0 < t_1 < \cdots < t_m < 1$
such that $\gamma(t_i) \in \tau_{j_i}$.
We illustrate this in Figure~\ref{fig:curve}.
Note that the curves $\tau_{i_1},\ldots,\tau_{i_m}$ do not have to be
pairwise different.
We do have, however, $\tau_{i_j} \not= \tau_{i_{j+1}}$ for all
$1 \le i \le m-1$.
The curve $\gamma^{-1}$ yields
$(b,\tau_{j_m},\ldots,\tau_{j_1},a)$.

\begin{figure}[H]
\begin{tikzpicture}
\draw[-](1,0) -- (1,2);
\draw [fill,red] (1,1) circle [radius=0.06];
\draw [fill] (1,0) circle [radius=0.06];
\draw [fill] (1,2) circle [radius=0.06];

\draw[-](9,0) -- (9,2);
\draw [fill] (9,0) circle [radius=0.06];
\draw [fill] (9,2) circle [radius=0.06];
\draw [fill,red] (9,1) circle [radius=0.06];

\draw[-](1,2) -- (9,2);
\draw[-](1,0) -- (9,0);

\draw[-](1,2) -- (3,0); \node[right] at (1,0.3) {$\tau_{j_1}$};
\draw[-](3,0) -- (3,2); \node[left] at (2.6,0.3) {$\tau_{j_2}$};
\draw[-](3,0) -- (5,2); \node[left] at (3,1.7) {$\tau_{j_3}$};
\draw[-](5,2) -- (5,0); \node[right] at (3.9,1.7) {$\tau_{j_4}$};
\draw[-](5,0) -- (7,2); \node[right] at (5.9,1.7) {$\tau_{j_6}$};
\draw[-](7,2) -- (7,0); \node[right] at (5.0,1.7) {$\tau_{j_5}$};
\draw[-](7,2) -- (9,0); \node[left] at (7,0.3) {$\tau_{j_7}$};
\node[left] at (8.6,0.3) {$\tau_{j_8}$};
\node[left] at (9,1.7) {$\tau_{j_1}$};

\draw[thick,red,-](1,1) -- (9,1);
\node[left] at (1,1) {$a$};
\node[right] at (9,1) {$a$};

\end{tikzpicture}
\caption{A loop $(a,\tau_{j_1},\ldots,\tau_{j_8},\tau_{j_1},a)$}\label{fig:loop}
\end{figure}
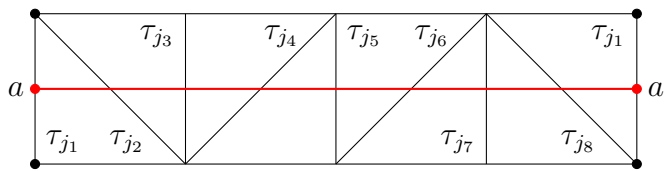

Analogously, with a loop $\gamma\df S^1\to \bbS$ in $\SM$
we associate a sequence 
$$
(a,\tau_{j_1},\tau_{j_2},\ldots,\tau_{j_m},\tau_{j_1},a),
$$ 
where $a = \gamma(1)$.
Starting in $1 \in S^1$ in clockwise orientation, we assume that 
$\gamma$ first passes through $\tau_{j_1}$, then through $\tau_{j_2}$ etc. 
We can assume here that
$a \in \tau_{j_1}$.
This is illustrated in Figure~\ref{fig:loop}.

The periodic curve $\widetilde{\gamma}\df \mathbb R \to \bbS$ associated with $\gamma$ is represented
by the infinite sequence
$$
(\cdots,\tau_{j_1},\tau_{j_2},\ldots,\tau_{j_m},
\tau_{j_1},\tau_{j_2},\ldots,\tau_{j_m},
\tau_{j_1},\tau_{j_2},\ldots,\tau_{j_m},
\cdots).
$$

Arguing similarly as in \cite[Section~4.2]{ABCP}, we can identify each 
$\gamma \in A\SM \cup L\SM$ with its associated sequence
$(a,\tau_{j_1},\ldots,\tau_{j_m},b)$
modulo the equivalence
$$
(a,\tau_{j_1},\ldots,\tau_{j_m},b) \sim (b,\tau_{j_m},\ldots,\tau_{j_1},a).
$$

\subsection{From triangulations to gentle Jacobian algebras}
\label{subsec:9.5}
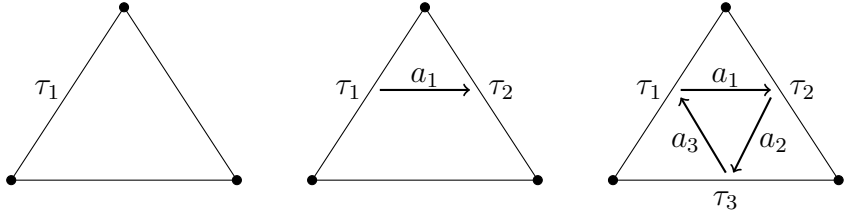
\begin{figure}[H]
\begin{tikzpicture}

\draw [fill] (0,0) circle [radius=0.06];
\draw [fill] (3,0) circle [radius=0.06];
\draw [fill] (1.5,2.3) circle [radius=0.06];
\draw [-] (0,0)--(3,0);
\draw [-] (3,0)--(1.5,2.3);
\draw [-] (1.5,2.3)--(0,0);\node[left] at (0.8,1.2) {$\tau_1$};


\draw [fill] (4,0) circle [radius=0.06];
\draw [fill] (7,0) circle [radius=0.06];
\draw [fill] (5.5,2.3) circle [radius=0.06];
\draw [-] (4,0)--(7,0);
\draw [-] (7,0)--(5.5,2.3);\node[right] at (6.2,1.2) {$\tau_2$};
\draw [-] (5.5,2.3)--(4,0);\node[left] at (4.8,1.2) {$\tau_1$};

\draw [thick,->] (4.9,1.2)--(6.1,1.2);
\node[above] at (5.5,1.1) {$a_1$};


\draw [fill] (8,0) circle [radius=0.06];
\draw [fill] (11,0) circle [radius=0.06];
\draw [fill] (9.5,2.3) circle [radius=0.06];
\draw [-] (8,0)--(11,0);\node[below] at (9.5,0) {$\tau_3$};
\draw [-] (11,0)--(9.5,2.3);\node[right] at (10.2,1.2) {$\tau_2$};
\draw [-] (9.5,2.3)--(8,0);\node[left] at (8.8,1.2) {$\tau_1$};

\draw [thick,->] (8.9,1.2)--(10.1,1.2);
\node[above] at (9.5,1.1) {$a_1$};
\draw [thick,->] (10.1,1.1)--(9.6,0.1);
\node[right] at (9.8,0.5) {$a_2$};
\draw [thick,->] (9.5,0.1)--(8.9,1.1);
\node[left] at (9.3,0.5) {$a_3$};

\end{tikzpicture}
\caption{How triangles yield arrows}\label{fig:triangles}
\end{figure}

Let $T$ be a triangulation of an unpunctured marked surface
$\SM$.
Assume that $T^\circ$ consists of $n$ curves $\tau_1,\ldots,\tau_n$.
Then $Q = Q_T$ is by definition the quiver with vertices
$1,\ldots,n$.
The arrows of $Q$ are defined as follows:
As displayed in Figure~\ref{fig:triangles},
there are three types of triangles defined by $T$,
and two of these yield arrows in $Q$, as indicated in
the picture.
Note that the non-labelled sides of the triangles are meant to be boundary segments of $\SM$, and note that our arrows point in
clockwise direction.
Other authors might choose the opposite convention.
The algebra associated to $T$ is then $A_T := KQ/I$,
where $I$ is generated by the paths
$a_2a_1,\; a_3a_2,\; a_1a_3$ arising from triangles with all three
sides in $T^\circ$.

The algebra $A_T$ was first studied by \cite{ABCP} and \cite{LF1},
where it was defined as the Jacobian algebra $\cP_\C(Q_T,W_T)$ of a quiver with potential.

\begin{Thm}[{\cite[Section~2]{ABCP}}]\label{thm1:ABCP}
The Jacobian algebras $A_T$ arising from triangulations of unpunctured marked surfaces 
are exactly the gentle Jacobian algebras.
\end{Thm}

\subsection{From curves and loops to string and band 
modules}
Let $\SM$ be an unpunctured marked surface,
and let $T$ be a fixed triangulation of 
$(\bbS,\bbM)$.

\begin{Thm}[{\cite[Propositions~4.2 and 4.3]{ABCP}}]\label{thm2:ABCP}
For $(\bbS,\bbM)$ and $T$ as above the following hold:
\begin{itemize}

\item[(i)]
There is a bijection $\gam \mapsto M_\gam$ between $A\SM \setminus T^\circ$
and the set of isoclasses of string modules in $\md(A_T)$. 

\item[(ii)]
There is a bijection $(\gam,\la) \mapsto M_{\gam,\la}$ between 
$L\SM \times K^*$
and the set of isoclasses of band modules in 
$\md(A_T)$.

\end{itemize}
\end{Thm}

The arcs in $T^\circ$ correspond bijectively to the negative simple
decorated $A_T$-modules.
Thus the isoclasses of indecomposable decorated $A_T$-modules are in bijection with 
$A\SM \cup L\SM \times K^*$.

For an indecomposable decorated $A_T$-module $\cM$
let $\gamma_{\cM}$ be the corresponding curve or
loop in $A\SM \cup L\SM$.
For $\cM = (M,0)$ we set $\gamma_M := \gamma_{\cM}$.

The string associated with the curve in Figure~\ref{fig:curve}
looks as follows:
$$
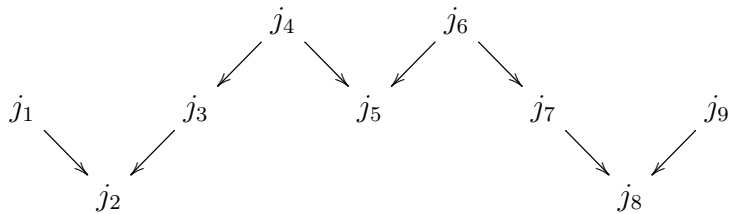
\xymatrix@-1.5ex{
&&& j_4 \ar[dl]\ar[dr] && j_6 \ar[dl]\ar[dr]
\\
j_1 \ar[dr] && j_3 \ar[dl] && j_5 && j_7 \ar[dr] && j_9 \ar[dl]
\\
& j_2 &&&&&& j_8
}
$$
The band
associated with the curve in Figure~\ref{fig:loop}
looks as in the following picture, where the two blue vertices have to be identified:
$$
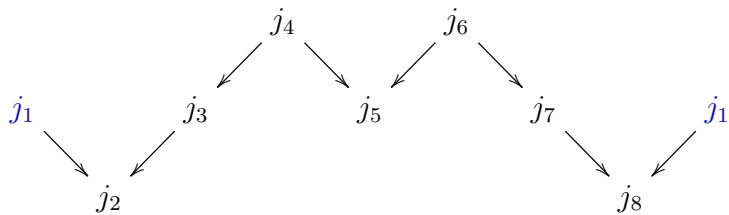
\xymatrix@-1.5ex{
&&& j_4 \ar[dl]\ar[dr] && j_6 \ar[dl]\ar[dr]
\\
{\blue j_1} \ar[dr] && j_3 \ar[dl] && j_5 && j_7 \ar[dr]
&& {\blue j_1} \ar[dl]
\\
& j_2 &&&&&& j_8
}
$$

Note that 
for an arbitrary gentle algebra $A$ there is also a geometric
model for the derived category $D^b(\md(A))$ 
(see \cite{HKK,LP,OPS}), which differs
substantially from the one for $\md(A)$ used in this article.

\begin{Thm}[{\cite[Corollary~5.4]{BZ}}]\label{thm:BZ}
Let $A = A_T$ as above, and
let
$M,N \in \md(A)$ be string modules.
Then the following are equivalent:
\begin{itemize}

\item[(i)]
$\Int(\gamma_M,\gamma_N) = 0$;

\item[(ii)]
$\Hom_A(N,\tau_A(M)) = 0 = \Hom_A(M,\tau_A(N)) = 0$.

\end{itemize}
\end{Thm}

Note that the results in \cite{BZ} are formulated in terms of the cluster
category associated with $\SM$. 
Theorem~\ref{thm:BZ} is a straightforward reformulation in terms of decorated $A_T$-modules.

In Section~\ref{subsec:proofmain7} we reprove and generalize Theorem~\ref{thm:BZ} by also including band modules.

\subsection{Rotation of curves and the Auslander-Reiten translation}
\label{subsec:rotation}
Let $T$ be a triangulation of $\SM$, and let $A = A_T$.
Let $M \in \md(A)$ be a string module, and let
$\gamma_M = (a,\tau_{j_1},\ldots,\tau_{j_t},b)$ be
the associated curve in $A\SM \setminus T^\circ$.

For the following two statements we refer to \cite[Theorem~3.6]{BZ}.
(Note that the orientation of our $Q_T$ is opposite to 
the one used in \cite{BZ}.)

We orient each boundary component of $\bbS$ by requiring that 
when following the orientation, the surface lies to the left.
We call this the \emph{induced orientation} of the boundary component.

If $M$ is non-projective, then 
$
\gamma_{\tau_A(M)} = \tau(\gamma),
$
where $\tau(\gamma)$ is obtained from $\gamma$ by rotating the
points
$a = \gamma(0)$ and $b = \gamma(1)$ of $\gamma$ 
to
the next marked point on their respective boundary component, following the induced orientation.

Dually, if $M$ is non-injective,  then 
$
\gamma_{\tau_A^{-1}(M)} = \tau^{-1}(\gamma),
$
where $\tau^{-1}(\gamma)$ is obtained from $\gamma$ by rotating 
$a$ and $b$ to
the next marked points on their respective boundary component,
following the opposite induced orientation.

The proof of these statements uses the combinatorial descriptions of $\tau_A(M)$ and $\tau_A^{-1}(M)$
given in \cite{BR} and \cite{WW}.

For more details we refer to \cite[Section~3]{BZ}.

\subsection{Three types of intersections}
Let $\SM$ be an unpunctured marked surface.
We fix a triangulation $T$ of $\SM$.
Now let $\gamma_1, \gamma_2 \in A\SM \cup L\SM$. 
Then the intersections of $\widetilde{\gamma}_1$ and 
$\widetilde{\gamma}_2$ can be
divided into three different types:
Type $I$ (resp. $II$) are displayed on the left (resp. right) in 
Figure~\ref{fig:intersec1}.
Up to symmetry there are 6 different kinds of
Type III intersections, which are pictured in the left hand column of Figure~\ref{fig:maintable1}.
(Note that the definition of intersection types depend here on our
fixed triangulation $T$.)

\begin{figure}[H]
\begin{tikzpicture}
\draw[-] (0,0)--(0,4);\node[left] at (0,2.5) {$\tau_1$};
\draw[-] (0,4)--(2,2);
\draw[-] (2,2)--(0,0);\node[below] at (1.5,1.4) {$\tau_2$};
\draw [fill] (0,0) circle [radius=0.06];
\draw [fill,red] (0,4) circle [radius=0.06];
\draw [fill,green] (2,2) circle [radius=0.06];

\draw[thick,green,-](-1,2) -- (2,2);
\node[left] at (-1,2) {$\widetilde{\gamma}_M$};
\draw[thick,red,-](0,4) -- (1.3,0.3);
\node[right] at (1.1,0) {$\widetilde{\gamma}_N$};


\draw[-] (5,0)--(5,4);\node[left] at (5,2.5) {$\tau_1$};
\draw[-] (5,4)--(7,2);\node[above] at (6.5,2.6) {$\tau_2$};
\draw[-] (7,2)--(5,0);\node[below] at (6.5,1.4) {$\tau_3$};
\draw [fill] (5,0) circle [radius=0.06];
\draw [fill] (5,4) circle [radius=0.06];
\draw [fill,red] (7,2) circle [radius=0.06];

\draw[thick,red,-](4,2) -- (7,2);
\node[left] at (4,2) {$\widetilde{\gamma}_N$};
\draw[thick,green,-](6,3.7) -- (6,0.3);
\node[above] at (6,3.7) {$\widetilde{\gamma}_M$};

\end{tikzpicture}
\caption{Intersections of types $I$ and $II$}\label{fig:intersec1}
\end{figure}
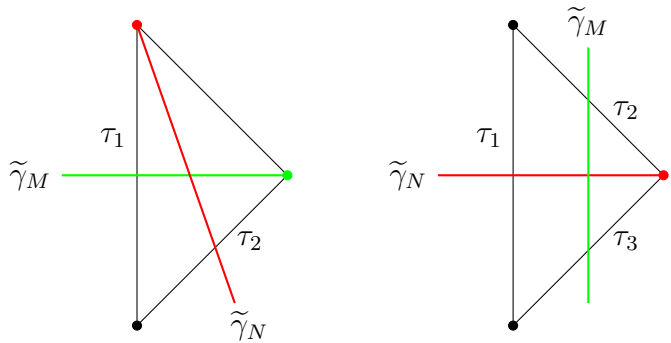
\begin{figure}[H]
\begin{tikzpicture}


\draw [fill] (0,16) circle [radius=0.06];
\draw [fill] (1,15) circle [radius=0.06];
\draw [fill] (1,17) circle [radius=0.06];
\draw[-](1,15) -- (1,17);
\draw[-](1,17) -- (0,16);
\draw[-](0,16) -- (1,15);
\draw [fill] (5,15) circle [radius=0.06];
\draw [fill] (5,17) circle [radius=0.06];
\draw [fill] (6,16) circle [radius=0.06];
\draw[-](5,15) -- (5,17);
\draw[-](5,17) -- (6,16);
\draw[-](6,16) -- (5,15);
\draw[-](1,17) -- (5,17);
\draw[-](1,15) -- (5,15);

\draw [thick,red,->] (0.6,16.4) -- (0.9,16);
\draw [thick,green,->] (0.9,15.9) -- (0.6,15.5);
\draw [thick,->] (0.5,15.6) -- (0.5,16.4);
\draw [thick,red,->] (5.4,15.6) -- (5.1,16);
\draw [thick,green,->] (5.1,16.1) -- (5.4,16.5);
\draw [thick,->] (5.5,16.4) -- (5.5,15.6);

\draw[thick,red,-](0,16.7) -- (6,15.3);
\draw[thick,green,-](0,15.3) -- (6,16.7);


\draw [fill,red] (8,16.6) circle [radius=0.06];
\draw [fill,red] (9.5,16.6) circle [radius=0.06];
\draw[thick,red,-](8,16.6) -- (9.5,16.6);
\draw [thick,red,->] (9.6,16.5) -- (9.9,16.1);
\draw [fill,red] (12.5,16.6) circle [radius=0.06];
\draw [fill,red] (14,16.6) circle [radius=0.06];
\draw[thick,red,-](12.5,16.6) -- (14,16.6);
\draw [thick,red,->] (12.4,16.5) -- (12.1,16.1);

\draw [fill,green] (8,15.4) circle [radius=0.06];
\draw [fill,green] (9.5,15.4) circle [radius=0.06];
\draw[thick,green,-](8,15.4) -- (9.5,15.4);
\draw [fill,blue] (10,16) circle [radius=0.06];
\draw [fill,blue] (12,16) circle [radius=0.06];
\draw[thick,blue,-](10,16) -- (12,16);
\draw [thick,green,->] (9.9,15.9) -- (9.6,15.5);
\draw [thick,green,->] (12.1,15.9) -- (12.4,15.5);
\draw [fill,green] (12.5,15.4) circle [radius=0.06];
\draw [fill,green] (14,15.4) circle [radius=0.06];
\draw[thick,green,-](12.5,15.4) -- (14,15.4);


\draw [fill] (0,13) circle [radius=0.06];
\draw [fill] (1,12) circle [radius=0.06];
\draw [fill] (1,14) circle [radius=0.06];
\draw[-](1,12) -- (1,14);
\draw[-](1,14) -- (0,13);
\draw[-](0,13) -- (1,12);
\draw [fill] (5,12) circle [radius=0.06];
\draw [fill] (5,14) circle [radius=0.06];
\draw [fill] (6,13) circle [radius=0.06];
\draw[-](5,12) -- (5,14);
\draw[-](5,14) -- (6,13);
\draw[-](6,13) -- (5,12);
\draw[-](1,14) -- (5,14);
\draw[-](1,12) -- (5,12);

\draw [thick,red,->] (0.6,13.4) -- (0.9,13);
\draw [thick,green,->] (0.9,12.9) -- (0.6,12.5);
\draw [thick,->] (0.5,12.6) -- (0.5,13.4);
\draw [dashed,->] (5.4,12.6) -- (5.1,13);
\draw [thick,green,->] (5.1,13.1) -- (5.4,13.5);
\draw [dashed,->] (5.5,13.4) -- (5.5,12.6);

\draw[thick,red,-](0,13.7) -- (6,13);
\draw[thick,green,-](0,12.3) -- (6,13.7);


\draw [fill,red] (8,13.6) circle [radius=0.06];
\draw [fill,red] (9.5,13.6) circle [radius=0.06];
\draw[thick,red,-](8,13.6) -- (9.5,13.6);
\draw [thick,red,->] (9.6,13.5) -- (9.9,13.1);

\draw [fill,green] (8,12.4) circle [radius=0.06];
\draw [fill,green] (9.5,12.4) circle [radius=0.06];
\draw[thick,green,-](8,12.4) -- (9.5,12.4);
\draw [fill,blue] (10,13) circle [radius=0.06];
\draw [fill,blue] (12,13) circle [radius=0.06];
\draw[thick,blue,-](10,13) -- (12,13);
\draw [thick,green,->] (9.9,12.9) -- (9.6,12.5);
\draw [thick,green,->] (12.1,12.9) -- (12.4,12.5);
\draw [fill,green] (12.5,12.4) circle [radius=0.06];
\draw [fill,green] (14,12.4) circle [radius=0.06];
\draw[thick,green,-](12.5,12.4) -- (14,12.4);


\draw [fill] (0,10) circle [radius=0.06];
\draw [fill] (1,9) circle [radius=0.06];
\draw [fill] (1,11) circle [radius=0.06];
\draw[-](1,9) -- (1,11);
\draw[-](1,11) -- (0,10);
\draw[-](0,10) -- (1,9);
\draw [fill] (5,9) circle [radius=0.06];
\draw [fill] (5,11) circle [radius=0.06];
\draw [fill] (6,10) circle [radius=0.06];
\draw[-](5,9) -- (5,11);
\draw[-](5,11) -- (6,10);
\draw[-](6,10) -- (5,9);
\draw[-](1,11) -- (5,11);
\draw[-](1,9) -- (5,9);

\draw [thick,red,->] (0.6,10.4) -- (0.9,10);
\draw [thick,green,->] (0.9,9.9) -- (0.6,9.5);
\draw [thick,->] (0.5,9.6) -- (0.5,10.4);
\draw [thick,red,->] (5.4,9.6) -- (5.1,10);
\draw [dashed,->] (5.1,10.1) -- (5.4,10.5);
\draw [dashed,->] (5.5,10.4) -- (5.5,9.6);

\draw[thick,red,-](0,10.7) -- (6,9.3);
\draw[thick,green,-](0,9.3) -- (6,10);


\draw [fill,red] (8,10.6) circle [radius=0.06];
\draw [fill,red] (9.5,10.6) circle [radius=0.06];
\draw[thick,red,-](8,10.6) -- (9.5,10.6);
\draw [thick,red,->] (9.6,10.5) -- (9.9,10.1);
\draw [fill,red] (12.5,10.6) circle [radius=0.06];
\draw [fill,red] (14,10.6) circle [radius=0.06];
\draw[thick,red,-](12.5,10.6) -- (14,10.6);
\draw [thick,red,->] (12.4,10.5) -- (12.1,10.1);

\draw [fill,green] (8,9.4) circle [radius=0.06];
\draw [fill,green] (9.5,9.4) circle [radius=0.06];
\draw[thick,green,-](8,9.4) -- (9.5,9.4);
\draw [fill,blue] (10,10) circle [radius=0.06];
\draw [fill,blue] (12,10) circle [radius=0.06];
\draw[thick,blue,-](10,10) -- (12,10);
\draw [thick,green,->] (9.9,9.9) -- (9.6,9.5);


\draw [red,fill] (0,7) circle [radius=0.06];
\draw [fill] (1,6) circle [radius=0.06];
\draw [fill] (1,8) circle [radius=0.06];
\draw[-](1,6) -- (1,8);
\draw[-](1,8) -- (0,7);
\draw[-](0,7) -- (1,6);
\draw [fill] (5,6) circle [radius=0.06];
\draw [fill] (5,8) circle [radius=0.06];
\draw [fill] (6,7) circle [radius=0.06];
\draw[-](5,6) -- (5,8);
\draw[-](5,8) -- (6,7);
\draw[-](6,7) -- (5,6);
\draw[-](1,8) -- (5,8);
\draw[-](1,6) -- (5,6);

\draw [dashed,->] (0.6,7.4) -- (0.9,7);
\draw [thick,green,->] (0.9,6.9) -- (0.6,6.5);
\draw [dashed,->] (0.5,6.6) -- (0.5,7.4);
\draw [dashed,->] (5.4,6.6) -- (5.1,7);
\draw [thick,green,->] (5.1,7.1) -- (5.4,7.5);
\draw [dashed,->] (5.5,7.4) -- (5.5,6.6);

\draw[thick,red,-](0,7) -- (6,7);
\draw[thick,green,-](0,6.3) -- (6,7.7);


\draw [fill,green] (8,6.4) circle [radius=0.06];
\draw [fill,green] (9.5,6.4) circle [radius=0.06];
\draw[thick,green,-](8,6.4) -- (9.5,6.4);
\draw [fill,blue] (10,7) circle [radius=0.06];
\draw [fill,blue] (12,7) circle [radius=0.06];
\draw[thick,blue,-](10,7) -- (12,7);
\draw [thick,green,->] (9.9,6.9) -- (9.6,6.5);
\draw [thick,green,->] (12.1,6.9) -- (12.4,6.5);
\draw [fill,green] (12.5,6.4) circle [radius=0.06];
\draw [fill,green] (14,6.4) circle [radius=0.06];
\draw[thick,green,-](12.5,6.4) -- (14,6.4);


\draw [fill] (0,4) circle [radius=0.06];
\draw [fill] (1,3) circle [radius=0.06];
\draw [fill] (1,5) circle [radius=0.06];
\draw[-](1,3) -- (1,5);
\draw[-](1,5) -- (0,4);
\draw[-](0,4) -- (1,3);
\draw [fill] (5,3) circle [radius=0.06];
\draw [fill] (5,5) circle [radius=0.06];
\draw [fill] (6,4) circle [radius=0.06];
\draw[-](5,3) -- (5,5);
\draw[-](5,5) -- (6,4);
\draw[-](6,4) -- (5,3);
\draw[-](1,5) -- (5,5);
\draw[-](1,3) -- (5,3);

\draw [thick,red,->] (0.6,4.4) -- (0.9,4);
\draw [dashed,->] (0.9,3.9) -- (0.6,3.5);
\draw [dashed,->] (0.5,3.6) -- (0.5,4.4);
\draw [thick,red,->] (5.4,3.6) -- (5.1,4);
\draw [dashed,->] (5.1,4.1) -- (5.4,4.5);
\draw [dashed,->] (5.5,4.4) -- (5.5,3.6);

\draw[thick,red,-](0,4.7) -- (6,3.3);
\draw[thick,green,-](0,4) -- (6,4);


\draw [fill,red] (8,4.6) circle [radius=0.06];
\draw [fill,red] (9.5,4.6) circle [radius=0.06];
\draw[thick,red,-](8,4.6) -- (9.5,4.6);
\draw [thick,red,->] (9.6,4.5) -- (9.9,4.1);
\draw [fill,red] (12.5,4.6) circle [radius=0.06];
\draw [fill,red] (14,4.6) circle [radius=0.06];
\draw[thick,red,-](12.5,4.6) -- (14,4.6);
\draw [thick,red,->] (12.4,4.5) -- (12.1,4.1);

\draw [fill,blue] (10,4) circle [radius=0.06];
\draw [fill,blue] (12,4) circle [radius=0.06];
\draw[thick,blue,-](10,4) -- (12,4);

\draw [fill] (0,1) circle [radius=0.06];
\draw [fill] (1,0) circle [radius=0.06];
\draw [fill] (1,2) circle [radius=0.06];
\draw[-](1,0) -- (1,2);
\draw[-](1,2) -- (0,1);
\draw[-](0,1) -- (1,0);
\draw [fill] (5,0) circle [radius=0.06];
\draw [fill] (5,2) circle [radius=0.06];
\draw [fill] (6,1) circle [radius=0.06];
\draw[-](5,0) -- (5,2);
\draw[-](5,2) -- (6,1);
\draw[-](6,1) -- (5,0);
\draw[-](1,2) -- (5,2);
\draw[-](1,0) -- (5,0);

\draw [dashed,->] (0.6,1.4) -- (0.9,1);
\draw [thick,green,->] (0.9,0.9) -- (0.6,0.5);
\draw [dashed,->] (0.5,0.6) -- (0.5,1.4);
\draw [thick,red,->] (5.4,0.6) -- (5.1,1);
\draw [dashed,->] (5.1,1.1) -- (5.4,1.5);
\draw [dashed,->] (5.5,1.4) -- (5.5,0.6);

\draw[thick,red,-](0,1) -- (6,0.3);
\draw[thick,green,-](0,0.3) -- (6,1);


\draw [fill,red] (12.5,1.6) circle [radius=0.06];
\draw [fill,red] (14,1.6) circle [radius=0.06];
\draw[thick,red,-](12.5,1.6) -- (14,1.6);
\draw [thick,red,->] (12.4,1.5) -- (12.1,1.1);

\draw [fill,green] (8,0.4) circle [radius=0.06];
\draw [fill,green] (9.5,0.4) circle [radius=0.06];
\draw[thick,green,-](8,0.4) -- (9.5,0.4);
\draw [fill,blue] (10,1) circle [radius=0.06];
\draw [fill,blue] (12,1) circle [radius=0.06];
\draw[thick,blue,-](10,1) -- (12,1);
\draw [thick,green,->] (9.9,0.9) -- (9.6,0.5);

\end{tikzpicture}
\caption{Type III intersections and 2-sided standard homomorphisms}\label{fig:maintable1}
\end{figure}
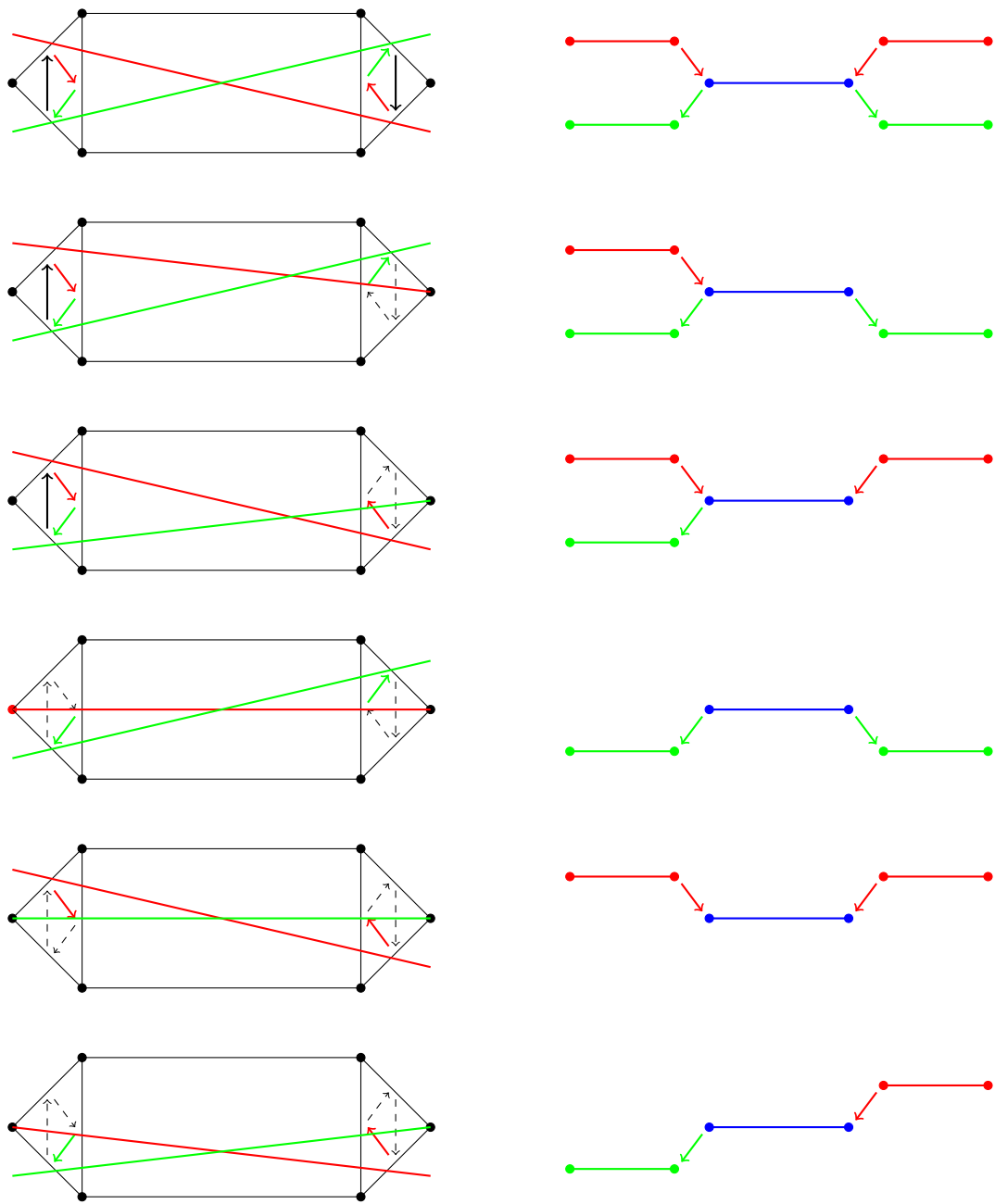
\begin{figure}[H]
\begin{tikzpicture}


\draw [fill] (0,10) circle [radius=0.06];
\draw [fill] (1,9) circle [radius=0.06];
\draw [fill] (1,11) circle [radius=0.06];
\draw[-](1,9) -- (1,11);
\draw[-](1,11) -- (0,10);
\draw[-](0,10) -- (1,9);
\draw [fill] (5,9) circle [radius=0.06];
\draw [fill] (5,11) circle [radius=0.06];
\draw [fill] (6,10) circle [radius=0.06];
\draw[-](5,9) -- (5,11);
\draw[-](5,11) -- (6,10);
\draw[-](6,10) -- (5,9);
\draw[-](1,11) -- (5,11);
\draw[-](1,9) -- (5,9);

\draw [thick,red,->] (0.6,10.4) -- (0.9,10);
\draw [thick,green,->] (0.9,9.9) -- (0.6,9.5);
\draw [thick,->] (0.5,9.6) -- (0.5,10.4);
\draw [dashed,->] (5.4,9.6) -- (5.1,10);
\draw [dashed,->] (5.1,10.1) -- (5.4,10.5);
\draw [dashed,->] (5.5,10.4) -- (5.5,9.6);

\draw[thick,red,-](0,10.7) -- (6,10);
\draw[thick,green,-](0,9.3) -- (6,10);


\draw [fill,red] (8,10.6) circle [radius=0.06];
\draw [fill,red] (9.5,10.6) circle [radius=0.06];
\draw[thick,red,-](8,10.6) -- (9.5,10.6);
\draw [thick,red,->] (9.6,10.5) -- (9.9,10.1);

\draw [fill,green] (8,9.4) circle [radius=0.06];
\draw [fill,green] (9.5,9.4) circle [radius=0.06];
\draw[thick,green,-](8,9.4) -- (9.5,9.4);
\draw [fill,blue] (10,10) circle [radius=0.06];
\draw [fill,blue] (12,10) circle [radius=0.06];
\draw[thick,blue,-](10,10) -- (12,10);
\draw [thick,green,->] (9.9,9.9) -- (9.6,9.5);


\draw [fill] (0,7) circle [radius=0.06];
\draw [fill] (1,6) circle [radius=0.06];
\draw [fill] (1,8) circle [radius=0.06];
\draw[-](1,6) -- (1,8);
\draw[-](1,8) -- (0,7);
\draw[-](0,7) -- (1,6);
\draw [fill] (5,6) circle [radius=0.06];
\draw [fill] (5,8) circle [radius=0.06];
\draw [fill] (6,7) circle [radius=0.06];
\draw[-](5,6) -- (5,8);
\draw[-](5,8) -- (6,7);
\draw[-](6,7) -- (5,6);
\draw[-](1,8) -- (5,8);
\draw[-](1,6) -- (5,6);

\draw [dashed,->] (0.6,7.4) -- (0.9,7);
\draw [thick,green,->] (0.9,6.9) -- (0.6,6.5);
\draw [dashed,->] (0.5,6.6) -- (0.5,7.4);
\draw [dashed,->] (5.4,6.6) -- (5.1,7);
\draw [dashed,->] (5.1,7.1) -- (5.4,7.5);
\draw [dashed,->] (5.5,7.4) -- (5.5,6.6);

\draw[thick,red,-](0,7) -- (6,7);
\draw[thick,green,-](0,6.3) -- (6,7);


\draw [fill,green] (8,6.4) circle [radius=0.06];
\draw [fill,green] (9.5,6.4) circle [radius=0.06];
\draw[thick,green,-](8,6.4) -- (9.5,6.4);
\draw [fill,blue] (10,7) circle [radius=0.06];
\draw [fill,blue] (12,7) circle [radius=0.06];
\draw[thick,blue,-](10,7) -- (12,7);
\draw [thick,green,->] (9.9,6.9) -- (9.6,6.5);


\draw [fill] (0,4) circle [radius=0.06];
\draw [fill] (1,3) circle [radius=0.06];
\draw [fill] (1,5) circle [radius=0.06];
\draw[-](1,3) -- (1,5);
\draw[-](1,5) -- (0,4);
\draw[-](0,4) -- (1,3);
\draw [fill] (5,3) circle [radius=0.06];
\draw [fill] (5,5) circle [radius=0.06];
\draw [fill] (6,4) circle [radius=0.06];
\draw[-](5,3) -- (5,5);
\draw[-](5,5) -- (6,4);
\draw[-](6,4) -- (5,3);
\draw[-](1,5) -- (5,5);
\draw[-](1,3) -- (5,3);

\draw [thick,red,->] (0.6,4.4) -- (0.9,4);
\draw [dashed,->] (0.9,3.9) -- (0.6,3.5);
\draw [dashed,->] (0.5,3.6) -- (0.5,4.4);
\draw [dashed,->] (5.4,3.6) -- (5.1,4);
\draw [dashed,->] (5.1,4.1) -- (5.4,4.5);
\draw [dashed,->] (5.5,4.4) -- (5.5,3.6);

\draw[thick,red,-](0,4.7) -- (6,4);
\draw[thick,green,-](0,4) -- (6,4);


\draw [fill,red] (8,4.6) circle [radius=0.06];
\draw [fill,red] (9.5,4.6) circle [radius=0.06];
\draw[thick,red,-](8,4.6) -- (9.5,4.6);
\draw [thick,red,->] (9.6,4.5) -- (9.9,4.1);

\draw [fill,blue] (10,4) circle [radius=0.06];
\draw [fill,blue] (12,4) circle [radius=0.06];
\draw[thick,blue,-](10,4) -- (12,4);


\draw [fill] (0,1) circle [radius=0.06];
\draw [fill] (1,0) circle [radius=0.06];
\draw [fill] (1,2) circle [radius=0.06];
\draw[-](1,0) -- (1,2);
\draw[-](1,2) -- (0,1);
\draw[-](0,1) -- (1,0);
\draw [fill] (5,0) circle [radius=0.06];
\draw [fill] (5,2) circle [radius=0.06];
\draw [fill] (6,1) circle [radius=0.06];
\draw[-](5,0) -- (5,2);
\draw[-](5,2) -- (6,1);
\draw[-](6,1) -- (5,0);
\draw[-](1,2) -- (5,2);
\draw[-](1,0) -- (5,0);

\draw [dashed,->] (0.6,1.4) -- (0.9,1);
\draw [dashed,->] (0.9,0.9) -- (0.6,0.5);
\draw [dashed,->] (0.5,0.6) -- (0.5,1.4);
\draw [dashed,->] (5.4,0.6) -- (5.1,1);
\draw [dashed,->] (5.1,1.1) -- (5.4,1.5);
\draw [dashed,->] (5.5,1.4) -- (5.5,0.6);

\draw[thick,red,-](0,1.05) -- (6,1.05);
\draw[thick,green,-](0,0.95) -- (6,0.95);


\draw [fill,blue] (10,1) circle [radius=0.06];
\draw [fill,blue] (12,1) circle [radius=0.06];
\draw[thick,blue,-](10,1) -- (12,1);

\end{tikzpicture}
\caption{Non-intersections and 1-sided standard homomorphisms}\label{fig:maintable2}
\end{figure}
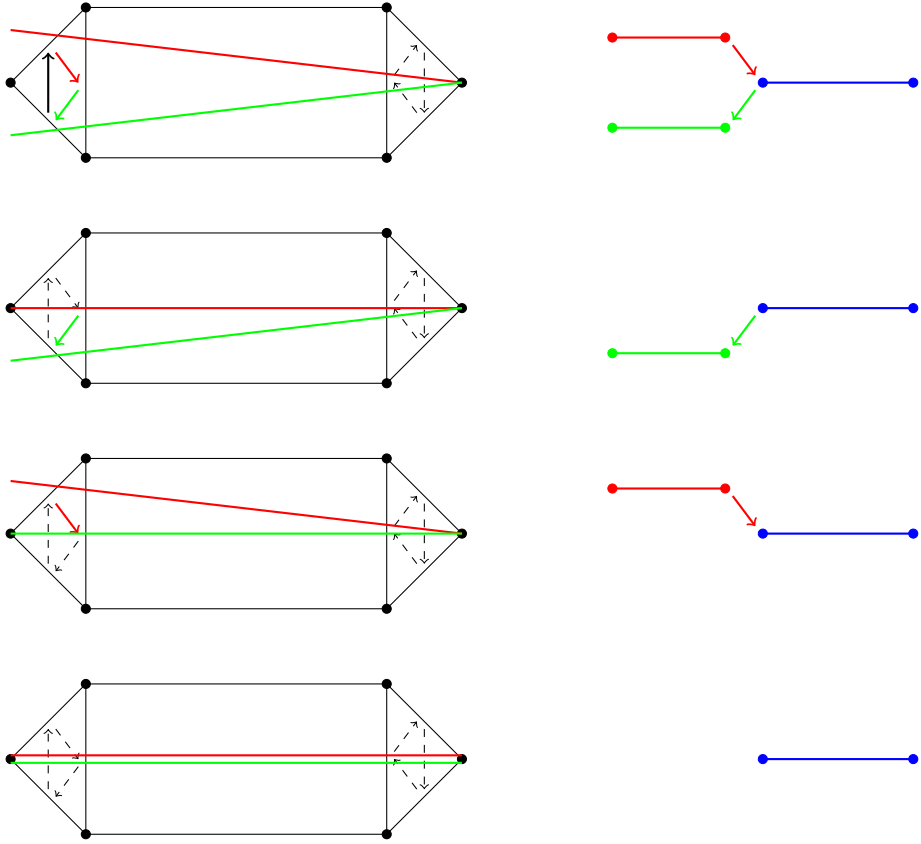

\subsection{Proof of Theorem~\ref{thm:main7}}
\label{subsec:proofmain7}
Throughout, we fix a triangulation $T$ of $\SM$.
Let $A = A_T$.

\begin{Lem}\label{lem:BZbands}
Let $M = M(B,\lambda,q) \in \md(A)$ be a band module
with $q \ge 2$.
Then $\Int(\gamma_M,\gamma_M) \not= 0$.
\end{Lem} 

\begin{proof}
Since $q \ge 2$, we have $\gamma_M = \gamma^q$ for some primitive loop $\gamma$.
It follows that $\Int(\gamma_M,\gamma_M) \not= 0$.
\end{proof}

\begin{Lem}\label{lem:negsimples}
Let $\cS_i^-$ be a negative simple decorated $A$-module, and let
$\cM$ be an indecomposable decorated $A$-module.
Then the following are equivalent:
\begin{itemize}

\item[(i)]
$\Int(\gamma_{\cS_i^-},\gamma_{\cM}) = 0$;

\item[(ii)]
$E_A(\cS_i^-,\cM) = 
E_A(\cM,\cS_i^-) = 0$.

\end{itemize}
\end{Lem}

\begin{proof}
Suppose that $\cM = \cS_j^-$ is also negative simple.
Then the equivalence of (i) and (ii) follows directly from the definitions.
Next, assume that $\cM = (M,0)$.
Let $\bd = (d_1,\ldots,d_n) = \dimv(M)$.
Then $d_i = 0$ if and only if
the simple $A$-module $S_i$ is not a composition factor of
$M$ if and only if (i) holds.
By the definition of $E_A(-,-)$, condition (ii) holds if and only if
$d_i = 0$.
This finishes the proof.
\end{proof}

In view of Lemma~\ref{lem:negsimples}, we can 
now restrict to indecomposable $A$-modules and curves
in $A\SM \setminus T^\circ$.

Recall that the notions of $\rho$-blocks and of the associated 
restriction maps $\pi_i$ were defined in Section~\ref{sec:blocks}.

\begin{Lem}\label{lem:intersec1}
Let $M$ and $N$ be indecomposable $A$-modules such that
$\widetilde{\gamma}_M$ and $\widetilde{\gamma}_N$ have a Type I or Type II intersection as shown
in Figure~\ref{fig:intersec1}.
Then 
$$
\Hom_A(M,\tau_A(N)) \not= 0.
$$
\end{Lem}

\begin{proof}
Assume we are in Type I:
Let $A_i$ be the $\rho$-block of $A$ containing  the arrow
$2 \to 1$.
Then $\pi_i(M) = S_1$ and 
$\pi_i(N) = S_2$.
By Lemma~\ref{lem:nonvanish2a} we get 
$\Ext_A^1(N,M) \not= 0$, which implies
$$
\Hom_A(M,\tau_A(N)) \not= 0.
$$

Next, assume we are in Type II:
Let $A_i$ be the $\rho$-block of $A$ containing the arrow
$1 \to 2$.
Thus $A_i$ also contains the arrows $2 \to 3$ and $3 \to 1$.
We get $\pi_i(M) = P_2$ and $\pi_i(N) = S_1$.
By Lemma~\ref{lem:nonvanish2b} this implies
$$
\Hom_A(M,\tau_A(N)) \not= 0.
$$
\end{proof}

\begin{Lem}\label{lem:intersec2}
Let $A_1,\ldots,A_t$ be the $\rho$-blocks of $A$.
Let $M$ and $N$ be indecomposable $A$-modules.
Then the following are equivalent:
\begin{itemize}

\item[(i)]
$\widetilde{\gamma}_M$ and $\widetilde{\gamma}_N$ have an intersection of type $I$ or $II$.

\item[(ii)]
For some $1 \le i \le t$,
$\pi_i(M \oplus N)$ is not $\tau$-rigid.
\end{itemize}
\end{Lem}

\begin{proof}
This is a direct consequence of Propositions~\ref{prop:complexrigid}
and \ref{prop:complextaurigid}.
\end{proof}

\begin{Lem}\label{lem:intersec3}
Let $M$ and $N$ be indecomposable $A$-modules.
Then the following are equivalent:
\begin{itemize}

\item[(i)]
$\widetilde{\gamma}_M$ and $\widetilde{\gamma}_N$ have a Type $III$ intersection,
as shown in the left column of Figure~\ref{fig:maintable1}
with $\widetilde{\gamma}_M$ green and $\widetilde{\gamma}_N$ red.

\item[(ii)]
There exists a $2$-sided standard homomorphism in
$\Hom_A(M,N)$.

\end{itemize}
\end{Lem}

\begin{proof}
This follows by a case by case inspection, which is carried out in
Figures~\ref{fig:maintable1} and \ref{fig:maintable2}.
\end{proof}

\begin{Lem}\label{lem:intersec4}
Let $M$ and $N$ be indecomposable $A$-modules.
If there exists a $2$-sided standard homomorphism 
in $\Hom_A(M,N)$, then $\Hom_A(M,\tau_A(N)) \not= 0$.
\end{Lem}

\begin{proof}
Assume that $M$ and $N$ are string modules.
It follows from \cite[Proposition~4.9]{Sch} that the existence of
a $2$-sided standard homomorphism in $\Hom_A(M,N)$ implies
$\Ext_A^1(N,M) \not= 0$.
By Theorem~\ref{ARformula1}(i), this yields 
$\Hom_A(M,\tau_A(N)) \not= 0$.

If $N$ is a band module, 
then $\tau_A(N) \cong N$, which implies the claim.

Finally, let $M$ be a band module and assume that 
$\Hom_A(M,N) \not= 0$.
Since $\tau_A(M) \cong M$ and $\idim(M) \le 1$ (see Lemma~\ref{lem:ARband}), we get from
Theorem~\ref{ARformula1}(iii) that
$$
0 \not= \Hom_A(M,N) \cong \Hom_A(\tau_A^{-1}(M),N)
\cong \Ext_A^1(N,M). 
$$
By Theorem~\ref{ARformula1}(i), this implies
$\Hom_A(M,\tau_A(N)) \not= 0$.
\end{proof}

For indecomposable $A$-modules $M$ and $N$, let
$\rad_A(M,N)$ be the non-invertible homomorphisms 
in $\Hom_A(M,N)$.
These form a subspace of $\Hom_A(M,N)$.

\begin{Lem}\label{lem:intersec5}
Let $M$ and $N$ be indecomposable $A$-modules.
Then the following hold:
\begin{itemize}

\item[(i)]
Let $N$ be a string module.
If $\Hom_A(M,\tau_A(N)) \not= 0$, then
$\Int(\widetilde{\gamma}_M,\widetilde{\gamma}_N) \not= 0$.

\item[(ii)]
Let $N$ be a band module of quasi-length $1$.
If $\rad_A(M,\tau_A(N)) \not= 0$, then
$\Int(\widetilde{\gamma}_M,\widetilde{\gamma}_N) \not= 0$.

\end{itemize}
\end{Lem}

\begin{proof}
(i)
Let $N = M(C)$ be a string module, and let
$f\df M \to \tau_A(N)$ be a standard homomorphism.
Thus, up to symmetry, $f$
is given by one of the ten pictures in Figures~\ref{fig:maintable1}
and \ref{fig:maintable2}.
The green curves in these pictures stand now for $\widetilde{\gamma}_M$ and the red curves for
$\widetilde{\gamma}_{\tau_A(N)}$.
Now $\tau^{-1}(\widetilde{\gamma}_{\tau_A(N)}) = \widetilde{\gamma}_N$ is obtained 
by a rotation in the direction opposite to the induced orientation.
By a straightforward case by case analysis we obtain $\Int(\widetilde{\gamma}_M,\widetilde{\gamma}_N) \not= 0$ in all ten cases.

(ii)
Let $N$ be a band module of quasi-length $1$.
Then $\tau_A(N) \cong N$.
Let $f_h \in \rad_A(M,N)$ be a standard homomorphism.
Since $N$ is a band module, we know that 
$h$ is of the form
$$
h = ((D_1,E_1,F_1),(D_2,E_2,F_2))
$$
with $l(D_2),l(F_2) \ge 1$.
Thus $f$ is described by one of the six cases in
Figures~\ref{fig:maintable1}, where the green curves in these pictures stand for $\widetilde{\gamma}_M$ and the red curves for
$\widetilde{\gamma}_N$.
(Actually we are then in 1st, 3rd or 5th case, where we count from
top to bottom.)
This implies 
$\Int(\widetilde{\gamma}_M,\widetilde{\gamma}_N) \not= 0$.
\end{proof}


\begin{Thm}\label{thm:main7b}
Let $M$ and $N$ be indecomposable $A$-modules.
If $M$ and $N$ are both band modules, then we assume that $M \not\cong N$.
Then the following are equivalent:
\begin{itemize}

\item[(i)]
$\Int(\gamma_M,\gamma_N) = 0$;

\item[(ii)]
$\Hom_A(M,\tau_A(N)) = 0$ and $\Hom_A(N,\tau_A(M)) = 0$.

\end{itemize}
\end{Thm}

\begin{proof}
(ii) $\implies$ (i):
This follows by combining 
Lemmas~\ref{lem:intersec1}, 
\ref{lem:intersec3} 
and \ref{lem:intersec4}. 

(i) $\implies$ (ii):
Assume that (ii) does not hold.
Without loss of generality let $\Hom_A(M,\tau_A(N)) \not= 0$.
If $N$ is a string module, then the result follows from 
Lemma~\ref{lem:intersec5}(i). 
Next, suppose $N = M(B,\lambda,q)$ is a band module.
The periodic curve $\widetilde{\gamma}_N$ and also the condition $\Hom_A(M,\tau_A(N)) \not= 0$ are independent of
$t$.
So we can assume that $q=1$.
By assumption we have $M \not\cong N$.
Thus $\rad_A(M,\tau_A(N)) \not= 0$.
Now the result follows from 
Lemma~\ref{lem:intersec5}(ii). 
\end{proof}

The following theorem corresponds to 
Theorem~\ref{thm:main7}.

\begin{Thm}\label{thm:main7c}
There is a bijection
$$
\eta_T\df \Lam\SM \to \decirr^\tau(A),
$$
which is natural in the sense that
$$
\eta_T(L) = \overline{\eta_T(\gamma_1,1)^{m_1} \oplus \cdots \oplus \eta_T(\gamma_t,1)^{m_t}}
$$
for each lamination $L = (\gamma,m)$ with $\gamma = \{ \gamma_1,\ldots,\gamma_t \}$ and $m(\gamma_i) = m_i$.
\end{Thm}

\begin{proof}
(a):
Let $M = M(C) \in \md(A)$ be a string module,
and let 
$$
Z_C := \overline{\cO_{(M,0)}} \subseteq \decrep(A,(\bd,0))
$$ 
where $\bd := \dimv(M)$.
By Theorem~\ref{thm:main7b} we have $\Int(\gamma_M,\gamma_M) = 0$ if and only if $\Hom_A(M,\tau_A(M)) = 0$ if and only if
$Z_C$  is a generically $\tau$-reduced decorated indecomposable irreducible component containing a dense orbit.

(b):
Next, let $M = M(B,\lambda,q) \in \md(A)$ be a band module, and
let 
$$
Z_{B,q} := \overline{\bigcup_{\lambda \in K^*} \cO_{(M(B,\lambda,q),0)}}
\subseteq \decrep(A,(\bd,0))
$$ 
where $\bd := \dimv(M)$.
If $q \ge 2$, then $\Int(\gamma_M,\gamma_M) \not= 0$.
Furthermore, 
$$
Z_{B,q} \subset \overline{Z_{B,1} \oplus \cdots \oplus Z_{B,1}}
$$
where $Z_{B,1} \oplus \cdots \oplus Z_{B,1}$ consists of
all decorated modules in $\decrep(A,(\bd,0))$ which are isomorphic to  
$(M(B,\lambda_1,1),0) \oplus \cdots
\oplus (M(B,\lambda_q,1),0)$ for some $(\lambda_1,\ldots,\lambda_q) \in (K^*)^q$.

Thus, we assume that $t=1$ and set $Z_B := Z_{B,1}$. 
Let $N = M(B,\mu,1)$ for some $\mu \in K^*$ with $\mu \not= \lambda$.
Note that $\gamma_M = \gamma_N$.
By Theorem~\ref{thm:main7b} we have $\Int(\gamma_M,\gamma_N) = 0$ if and only if $\Hom_A(M,N) = 0$ if and only if
$Z_B$  is a generically $\tau$-reduced decorated indecomposable irreducible component not containing a dense orbit.
Note here that $\tau_A(M) \cong M$ and $\tau_A(N) \cong N$ and that the condition $\Hom_A(M,N) = 0$ is equivalent to the condition
$\End_A(M) \cong K$.

(c):
The considerations in (a) and (b) show that there is a bijection
between the set of indecomposable components in
$\irr^\tau(A)$ and the set of laminations of the form
$L = (\{ \gamma_1 \},m)$ with $m(\gamma_1) = 1$ and
$\gamma_1 \notin T^\circ$.

(d):
Now the Theorem follows from 
Lemma~\ref{lem:negsimples} (which takes care of the negative simple decorated
modules) and 
Theorem~\ref{thm:main7b}
combined with the Decomposition Theorem~\ref{thm:decomp}.
\end{proof}

\subsection{Shear coordinates and $g$-vectors} \label{ssec:shearg}
Let $A = A_T$ as above.
As mentioned before,
a result by W. Thurston (see \cite[Theorem~12.3]{FT})
says that there is a bijection 
$\bs_T\df \Lam\SM \to \Z^n$ 
sending a lamination to its shear coordinate.
We briefly and informally recall the construction of $\bs_T$.

First, consider an arc
$$
\gamma = (a,\tau_{j_1},\ldots,\tau_{j_m},b) \in A\SM.
$$
Then
$$
\tau^{1/2}(\gamma) = 
(a',\tau_{j_{11}},\ldots,\tau_{j_{1t_a}},\tau_{j_1},\ldots,\tau_{j_m},\tau_{j_{m1}},\ldots,\tau_{j_{mt_b}},b'),
$$
where $a',b' \in \partial \bbS \setminus \bbM$, and
$(\tau_{j_{11}},\ldots,\tau_{j_{1t_a}})$ and 
$(\tau_{j_{m1}},\ldots,\tau_{j_{mt_b}})$ are possibly empty sequences
of curves in $T^\circ$ which are incident with $a$ and $b$, respectively.
Let $\tau_{a'}$ and $\tau_{b'}$ the boundary segments, which contain
$a'$ and $b'$, respectively.

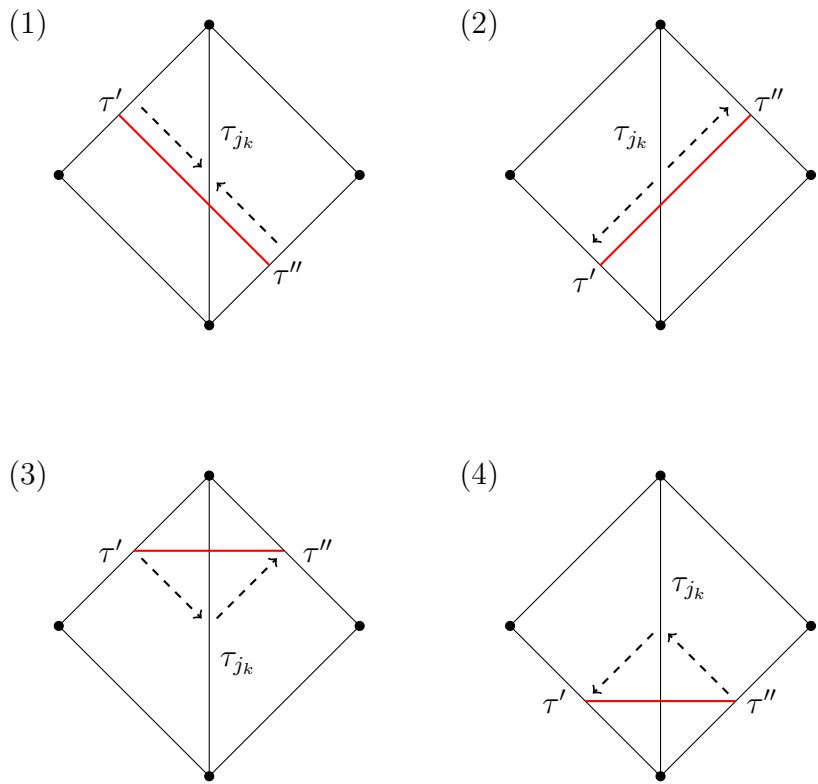
\begin{figure}[H]
\begin{tikzpicture}

\node [left] at (0,10) {$(1)$};
\draw [fill] (0,8) circle [radius=0.06];
\draw [fill] (2,10) circle [radius=0.06];
\draw [fill] (2,6) circle [radius=0.06];
\draw [fill] (4,8) circle [radius=0.06];
\draw[-](0,8) -- (2,10);
\draw[-](2,10) -- (4,8);
\draw[-](4,8) -- (2,6);
\draw[-](2,6) -- (0,8);
\draw[-](2,10) -- (2,6);
\node [right] at (2,8.5) {$\tau_{j_k}$};
\draw[thick,dashed,->] (1.1,8.9)--(1.9,8.1);
\draw[thick,dashed,->] (2.9,7.1)--(2.1,7.9);
\draw[thick,red,-] (0.8,8.8)--(2.8,6.8);
\node [left] at (1.0,9.0) {$\tau'$};
\node [right] at (2.7,6.7) {$\tau''$};


\node [left] at (6,10) {$(2)$};
\draw [fill] (6,8) circle [radius=0.06];
\draw [fill] (8,10) circle [radius=0.06];
\draw [fill] (8,6) circle [radius=0.06];
\draw [fill] (10,8) circle [radius=0.06];
\draw[-](6,8) -- (8,10);
\draw[-](8,10) -- (10,8);
\draw[-](10,8) -- (8,6);
\draw[-](8,6) -- (6,8);
\draw[-](8,10) -- (8,6);
\node [left] at (8,8.5) {$\tau_{j_k}$};
\draw[thick,dashed,->] (7.9,7.9)--(7.1,7.1);
\draw[thick,dashed,->] (8.1,8.1)--(8.9,8.9);
\draw[thick,red,-](7.2,6.8)--(9.2,8.8);
\node [left] at (7.3,6.6) {$\tau'$};
\node [right] at (9.1,9.0) {$\tau''$};


\node [left] at (0,4) {$(3)$};
\draw [fill] (0,2) circle [radius=0.06];
\draw [fill] (2,4) circle [radius=0.06];
\draw [fill] (2,0) circle [radius=0.06];
\draw [fill] (4,2) circle [radius=0.06];
\draw[-](0,2) -- (2,4);
\draw[-](2,4) -- (4,2);
\draw[-](4,2) -- (2,0);
\draw[-](2,0) -- (0,2);
\draw[-](2,4) -- (2,0);
\node [right] at (2,1.5) {$\tau_{j_k}$};
\draw[thick,dashed,->] (1.1,2.9)--(1.9,2.1);
\draw[dashed,thick,->] (2.1,2.1)--(2.9,2.9);
\draw[thick,red,-] (1.0,3.0)--(3.0,3.0);
\node [left] at (1.0,3.0) {$\tau'$};
\node [right] at (3.1,3.0) {$\tau''$};


\node [left] at (6,4) {$(4)$};
\draw [fill] (6,2) circle [radius=0.06];
\draw [fill] (8,4) circle [radius=0.06];
\draw [fill] (8,0) circle [radius=0.06];
\draw [fill] (10,2) circle [radius=0.06];
\draw[-](6,2) -- (8,4);
\draw[-](8,4) -- (10,2);
\draw[-](10,2) -- (8,0);
\draw[-](8,0) -- (6,2);
\draw[-](8,4) -- (8,0);
\node [right] at (8,2.5) {$\tau_{j_k}$};
\draw[thick,red,-] (7.0,1.0)--(9.0,1.0);
\node [left] at (6.9,1.0) {$\tau'$};
\node [right] at (9.0,1.0) {$\tau''$};
\draw[thick,dashed,->] (7.9,1.9)--(7.1,1.1);
\draw[thick,dashed,->] (8.9,1.1)--(8.1,1.9);

\end{tikzpicture}
\caption{Computing shear coordinates}\label{fig:shear}
\end{figure}

For each $1 \le k \le m$, we look at the triple 
$(\tau',\tau_{j_k},\tau'')$, where $\tau'$ and $\tau''$ are
the left and right neighbour, respectively, of $\tau_{j_k}$ in
the sequence 
$$
(\tau_{a'},\tau_{j_{11}},\ldots,\tau_{j_{1t_a}},\tau_{j_1},\ldots,\tau_{j_m},\tau_{j_{m1}},\ldots,\tau_{j_{mt_b}},\tau_{b'}).
$$
Then we are in one of the four cases displayed in Figure~\ref{fig:shear}, where the red line is a segment of the
curve $\tau^{1/2}(\gamma)$ and the dotted arrows indicate
possible arrows of $A$. (There is an arrow on the left if and only if
$\tau' \not= \tau_{a'}$, and there is an arrow on the right if and only if
$\tau'' \not= \tau_{b'}$.)

Next, consider a simple loop
$$
\gamma = (a,\tau_{j_1},\ldots,\tau_{j_m},\tau_{j_1},a) \in L\SM.
$$
For each $1 \le k \le m$, we look at the triple 
$$
(\tau',\tau_{j_k},\tau'') :=
\begin{cases}
(\tau_{j_{k-1}},\tau_{j_k},\tau_{j_{k+1}}) & \text{if $2 \le k \le m-1$},
\\
(\tau_{j_m},\tau_{j_1},\tau_{j_2}) & \text{if $k=1$},
\\
(\tau_{j_{m-1}},\tau_{j_m},\tau_{j_1}) & \text{if $k = m$}.
\end{cases}
$$

In both cases (i.e. $\gamma \in A\SM$ and $\gamma \in L\SM$),
the \emph{shear coordinate} of $\gamma$ (with respect to $T$) is defined as
$\bs_T(\gamma) := (s_1,\ldots,s_n)$,
where
$$
s_i := \sum_{k=1}^m \delta_{j_k,i} \delta_k
$$
for $1 \le i \le n$.
Here $\delta_{j_k,i}$ denotes the Kronecker delta and
$$
\delta_k :=
\begin{cases}
1 & \text{if $(\tau',\tau_{j_k},\tau'')$ looks as in case (1) of Figure~\ref{fig:shear}},
\\
-1 & \text{if we are in case (2)},
\\
0 & \text{if we are in cases (3) or (4)}.
\end{cases}
$$

Finally, let $L = (\gamma,m)$ be a lamination.
Then
$$
\bs_T(L) := \sum_{\gamma_i \in \gamma} 
m(\gamma_i) \bs_T(\gamma_i).
$$

Recall that by Plamondon \cite[Theorem~1.2]{P1}, there is a bijection 
$\bg_T\df \decirr(A)^\tau \to \Z^n$ 
sending a generically $\tau$-reduced decorated component 
to its $g$-vector.

The proof of the following result is a bit tedious but straightforward.
It follows essentially the ideas from
Labardini-Fragoso \cite[Theorem~10.0.5]{LF2}.
Note that \cite{LF2} deals with a dual situation and only
considers curves. 
The case of loops is however easier than the curve case and uses the same arguments.
Note also that \cite{LF2} uses a different (but equivalent) definition
of $g$-vectors.

\begin{Prop} \label{prp:eta}
With $A = A_T$ as above, the diagram
$$
\xymatrix{ 
\Lam\SM \ar[d]_{\eta_T} \ar[r]^<<<<{\bs_T} & \Z^n \ar@{=}[d]
\\
\decirr^\tau(A) \ar[r]^<<<<{\bg_T} & \Z^n
}
$$
commutes.
\end{Prop}


\section{Bangle functions and generic Caldero-Chapoton functions}
\label{sec:bangle}


We will assume throughout that our surface with marked points
$\SM$ is connected and has no punctures.
We fix a triangulation $T$ with
internal edges $T^\circ=(\tau_1, \tau_2,\ldots,\tau_n)$.

\subsection{Strings and Bands}\label{ssec:MSW-SB}
Recall from Section~\ref{ssec:curv-comb} that we identify each
curve $\gam\in A\SM \setminus T^\circ$
with a certain sequence $(a,\tau_{j_1},\ldots,\tau_{j_m},b)$, where
$a,b\in\bbM$ and the $\tau_{j_i}$ are the sequence of arcs of $T^\circ$ which are
crossed by $\gam$ in a minimal way, up to homotopy. Denote by $\Del_i$ the triangle of $T$,
which contains the arcs $\tau_{j_i}$ and $\tau_{j_{i+1}}$, and which contains
the segment $[\gam(t_i),\gam(t_{i+1})]$  of $\gam$ 
for $1 \le i \le m-1$.
This sequence can be coded into a (decorated) quiver $Q^T_\gam$ of type
$\sfA_m$ with vertices $\{1,2, \ldots, m\}$. Now, in $\Del_i$ there exists
an unique arrow $a_i$ of the quiver $Q_T$ (see Section~\ref{subsec:9.5}),
which goes either from $\tau_{j_i}$ to $\tau_{j_{i+1}}$, or from
$\tau_{j_{i+1}}$ to $\tau_{j_i}$.  In the first case we draw an arrow with
label $a_i$ from $i$ to $i+1$. In the second case, we draw an arrow with
the same label from $i+1$ to $i$. We call $Q^T_\gam$ the \emph{string}
of $\gam$  with respect to the triangulation $T$.

Analogously, we associate with a loop
$\gam=(a,\tau_{j_1},\ldots,\tau_{j_m},\tau_{j_1},a)\in L\SM$
a quiver of type $\tA_{m-1}$ with vertices $\{1,2,\ldots, m\}$.
The only difference is
that now we have an additional triangle  $\Del_m$, which  contains the edges
$\tau_{j_m}$, $\tau_{j_1}$, and the segment $[\gam(t_m),\gam(t_1)]$ of $\gam$. 
In this case $\Del_m$ determines the direction of the arrow between $a_m$
between $1$ and $m$.
We call in this case $Q^T_\gam$ the \emph{band} of $\gam$ with respect to $T$.

\subsection{MSW-functions}
In this section we will use the conventions and definitions
from \cite[Section~3]{MSW2} without further reference.

Musiker, Schiffler and Williams \cite{MSW2} assign to each homotopy
class $\gamma \in A\SM$ (resp. $\gam \in L\SM$)
a snake graph (resp. band graph) $G=G_{T,\gam}$.
We assume that in each tile $G_1, G_2,\ldots G_l$
of $G$, the diagonal goes from SE to NW,
and we always think that $G$ is drawn from SW to NE.

\begin{Rem}\label{Rem:10.1}
The  graph $G$ comes with a  distinguished good resp. perfect
matching $P_-$ which consists of the external edges of $G$
which are either vertical and belong to a negatively oriented tile,
or are horizontal and belong to a positively oriented tile.  
On the other hand, the tile $G_j$  and the position of its two neighbours
record how $\gam$ crosses the quadrilateral surrounding $\tau_{j_i}$ in the
neighbourhood of $\gam(t_i)$.
With these two observations it is an easy exercise to show that
\begin{equation}\label{eq:MSW-shear}
\frac{x(P_-)}{\cross(T,\gam)} =\ux^{\bs_T(\gam)},
\end{equation}
where $\bs_T(\gam)$ is the shear coordinate
vector (see Section~\ref{ssec:shearg})
of $\gam$ with respect to $T$.
\end{Rem}

\begin{Rem}
Following Derksen-Weyman-Zelevinsky \cite[p.~60]{DWZ1} 
each skew-symmetric matrix $B\in\Z^{n\times n}$ corresponds to a 2-acyclic quiver
$Q(B)$ with vertices
$\{1,2,\ldots, n\}$ and $b_{ij}$ arrows from $j$ to $i$ whenever $b_{ij}>0$.

In \cite[Definition~2.19]{MSW2} the (skew-symmetric) signed adjacency matrix
$B_T \in\Z^{n\times n}$ of a triangulation $T$ of $\SM$ is introduced.
With these conventions in place we have $Q_T=Q(-B_T)$ for our quiver
$Q_T$ from Section~\ref{subsec:9.5}.
The (coefficient-free) cluster algebra $\cA(B_T)$ associated with $B_T$ is
just $\cA_{\SM}$.
Let $\cA_\bullet(B_T)$ be the corresponding cluster algebra with principal coefficients.
\end{Rem}

\begin{Rem}
In \cite[Definitions~5.3 and 5.6]{MSW2} the authors associate to their graph
$G=G_{T,\gam}$ a poset structure $Q_G$ on the set $\{1,2,\ldots, m\}$
by describing its Hasse quiver.
We leave it as an exercise that our quiver $Q_\gam^T$ from
Section~\ref{ssec:MSW-SB}
is  opposite  to the Hasse quiver in \cite{MSW2}.
Thus, the poset ideals of $Q_G$ are
precisely the subsets $I$ of vertices of $Q_\gam^T$ which are closed under
predecessors. We call  such subsets \emph{order coideals} of $Q_\gam^T$.
\end{Rem}

In \cite[Definition~3.4]{MSW2} a Laurent polynomial
\begin{equation}
  X^T_\gam=\frac{1}{\cross(T,\gam)}\sum_{P} x(P)y(P) \in R :=
  \Z[x_i^\pm,y_i]_{i=1,2,\ldots, n} 
\end{equation} 
is defined, where the sum runs over the perfect resp.~good matchings of $G$.
We agree that $X^T_{\gam_i}=x_i$ for $\gam_i\in T^\circ$ and for
$L=(\xi, m)\in\Lam\SM$ one sets
\[
  X^T_L:= \prod_{\gam\in\xi} (X^T_\gam)^{m(\gam)}.
\]
The following result is  implicit in \cite[Sections~5 and 6]{MSW2}:
\begin{Lem}\label{Lem:MSW}
For each $\gam\in A\SM \cup L\SM$ holds
  \[
    X^T_\gam= \ux^{\bs_T(\gam)}\sum_{I\subset Q_\gam^T} \prod_{i\in I} \hy_{j_i},
  \]
where the summation runs over the order coideals $I$ of $Q_\gam^T$ and
\[
  \hy_j:= y_j\cdot \prod_{i=1}^n x_i^{b_{ij}}
  \in\cA_\bullet(B_T)
\]  
for $j=1,2,\ldots,n$.  
\end{Lem}

\begin{proof}
According to \cite[Theorem~5.7]{MSW2} the lattice $L(G)$ of good matchings of $G$
is in natural bijection with the distributive lattice of order coideals of
$Q_\gam^T$. More
precisely, to a good matching corresponds the coideal $I(P)$, which consists
of the labels of the tiles of $G$ which are enclosed by $P\ominus P_-$.

On the other hand, by \cite[Proposition~6.2]{MSW2} $x_\gam\in R$ is
homogeneous of degree
\[
  \bg(x_\gam)=\deg\left(\frac{x(P_-)}{\cross(T,\gam)}\right),
\]
if we agree that $\deg{x_j}=\be_j\in\Z^n$ and
$\deg{y_j} = -\sum_{i=1}^n b_{ij}\be_i \in \Z^n$. 

Thus in view of~\eqref{eq:MSW-shear} we have to show that
\begin{equation} \label{eq:MSW3}
  \frac{x(P)y(P)}{x(P_-)} =\prod_{i\in I(P)} \hy_{j_i}
 \qquad\text{ for all good matchings } P \text{ of }G.
\end{equation}
In order to show~\eqref{eq:MSW3} we proceed by induction on the Hasse diagram
of the distributive lattice 
$L(G)$ as in the proof of \cite[Theorem~5.1]{MSW2} at the end of \cite[Section~5]{MSW2}.
\end{proof}

\subsection{Dual CC-functions and MSW-functions}
\label{subsec:10.3}
We introduce the \emph{anti principal ice quiver}
$\tQ_T$, which is obtained
from $Q_T$ by adding an additional set of frozen vertices
$\{1',2',\ldots,n'\}$, and an additional arrow $p_i\df i'\ra i$ for
$i=1,2,\ldots, n$. 
The potential $W_T$ mentioned in Section~\ref{subsec:9.5}
can be
naturally viewed as a potential for $\tQ_T$ and it is not hard to see
that $(\tQ_T,W_T)$ is a non-degenerate QP with finite-dimensional Jacobian
algebra $\wtA_T = \cP_\C(\tQ_T,W_T)$.

\begin{Def} The \emph{dual Caldero-Chapoton function} with respect to $\wtA_T$
of a decorated representation $\cM=(M,V)$ of $A_T$ is the Laurent polynomial
  \[
  \CaCh'_{\wtA_T}(\cM):=
  \ux^{\bg_{\wtA_T}(\cM)}\sum_{\be\in\N^n} \chi(\Gr_{A_T}^\be(M))\cdot\uhy^\be\in R,
\]
where $\Gr_{A_T}^\be(M)$ is the quiver Grassmannian of factor modules with dimension vector $\be$  of the $A_T$-module $M$, and $\chi$ is
the  topological Euler characteristic.
\end{Def}

Note that for a decorated representation $\cM$ of $A_T$ we have in fact
$\bg_{\wtA_T}(\cM) = (\bg_{A_T},0,\ldots, 0)$.  This is so, since for each
$A_T$-module $M$ with minimal projective presentation  
\[
  P_1\ra P_0\ra M\ra 0,
\]
the same sequence can be taken as a minimal projective presentation of $M$
viewed as an $\wtA_T$-module, due to the shape of $\tQ_A$.

\begin{Rem}
  Obviously, the dual Caldero-Chapoton-function is the same as the usual
  Caldero-Chapoton-function for the corresponding dual module, more precisely 
\[
  \CaCh'_{\wtA_T}(\cM)=\CaCh_{\wtA_T^\op}(D\cM):=
  \ux^{\bg'_{\wtA_T^\op}(D\cM)}
  \sum_{\be\in\N^n} \chi(\Gr^{A_T^\op}_\be(DM))\cdot\uhy^\be\in R,
\]
where $D\cM=(DM,DV)$ is the $\C$-dual decorated $A_T^\op$-module, 
\[
\bg'_{\wtA_T^\op}(\cM)=\bg_{A_T^\op}(DM)+\dimv(DV)
\]
is the classical $g$-vector, calculated in terms of the minimal
injective copresentation $0\ra DM\ra DP_0\ra DP_1$, and
$\Gr^{A_T^\op}_\be(DM)$ is the quiver Grassmannian of $\be$-dimensional
$A_T^\op$-submodules of $DM$.

Thus we have in particular
\begin{equation}\label{eq:MSW-CC-mult}
  \CaCh'_{\wtA_T}(\cM_1\oplus\cM_2)=\CaCh'_{\wtA_T}(\cM_1)\CaCh'_{\wtA_T}(M_2)
\end{equation}  
for decorated representations $\cM_1$ and $\cM_2$.

Moreover the $\CaCh(\cM)$ for decorated reachable $E$-rigid $A_T$-modules
$\cM$ are precisely the cluster monomials for the cluster algebra
$\cA_\bullet(B_T) \subset R$ with principal coefficients, see for
example \cite{DWZ2}.
\end{Rem}

\begin{Rem} \label{MSW:RemCoeff}
For a curve $\gam \in A\SM \setminus T^\circ$ let 
$\cM_\gam := (M_\gam,0)$ be the corresponding decorated $A_T$-module.
For a primitive $\gam \in L\SM$ let 
$\cM_\gam := (M_{\gam,\la},0)$ for some $\la \in \C^*$.
Note that $M_{\gam,\la}$ is a band module of quasi-length $1$.
In these two cases, the
quiver $Q_\gam^T$ is the coefficient quiver of the string module
$M_\gam$ (resp. of the
band module $M_{\gam,\la}$).
Moreover, the order coideals of $Q_\gam^T$ can be identified with the
coordinate factor modules of $\cM_\gam$, see also \cite[Remark~5.8]{MSW2}.
Finally, for $\gam \in T^\circ$ let $\cM_\gam$ be the associated negative simple decorated $A_T$-module. 
\end{Rem}

\begin{Prop}\label{MSWProp-CC}
For  a curve or primitive loop $\gam\in A\SM \cup L\SM$
we have
\[
X^T_\gam = \CaCh'_{\wtA_T}(\cM_\gam).
\]
\end{Prop}

\begin{proof}
We  use Lemma~\ref{Lem:MSW} to compare both expressions.  
As a consequence of Proposition~\ref{prp:eta}, we get
$\bs_T(\gam) = \bg_{A_T}(\cM_\gam) = \bg_{\wtA_T}(\cM_\gam)$.
In view of Remark~\ref{MSW:RemCoeff} our claim follows now
from \cite[Theorem~1.2]{Hau}.  
\end{proof}

\subsection{Bangle functions are generic}
Recall that our set of laminations $\Lam\SM$ from Section~\ref{ssec:LamTri}
is  the same as the set of $\cC^\circ\SM$ of $\cC^\circ$-compatible
collection
of arcs and simple (= essential) loops in \cite[Def.~3.17]{MSW2}. 

Recall also that each irreducible component $Z\in\decirr^\tau(A_T)$ we can
consider the map
\[
  \CaCh'_Z\df Z\ra R,\quad\cM\mapsto\CaCh'_{\wtA_T}(\cM)
\]
as a constructible function, which indeed has a finite image. 
Thus there
exists an open dense subset $U \subseteq Z$ where $\CaCh'_Z$ takes a constant
value, say $\CaCh'_{\wtA_T}(Z)$.  
We define
$$
\widetilde{\cG}_T := 
\{ \CaCh'_{\wtA_T}(Z) \mid Z \in \decirr^\tau(A_T) \}
$$
and
$$
\widetilde{\cB}_T := 
\{ X_L^T \mid L \in \Lam\SM \}.
$$

With this definition we can state now
the main result of this section:

\begin{Thm}\label{thm:10.9}
For each lamination $L\in\Lam\SM$ we have
\[
    X^T_L=\CaCh'_{\wtA_T}(\eta_T(L)),
\]
where $\eta_T\df \Lam\SM \to \decirr^\tau(A_T)$ is the bijection 
from Theorem~\ref{thm:main7c}.
In particular,
we have 
$$
\widetilde{\cB}_T = \widetilde{\cG}_T.
$$
\end{Thm}

\begin{proof}
If an irreducible component $Z\in\decirr(A_T)$ decomposes as
  $Z=\overline{Z_1\oplus Z_2}$  then it follows from~\eqref{eq:MSW-CC-mult}
  and the above definition that
\[
  \CaCh'_{\wtA_T}(Z)=\CaCh'_{\wtA_T}(Z_1)\cdot\CaCh'_{\wtA_T}(Z_2).
\]  
Let $L=(\xi, m)\in\Lam\SM$ be a lamination,

In Theorem~\ref{thm:main7c} we assign  to $(\xi, m)$ a
generically $\tau$-reduced decorated irreducible component
\[
\eta_T(\gam,m) =
\overline{\bigoplus_{\gam\in\xi} \eta_T(\gam,1)^{m(\gam)}}\in\decirr^\tau(A_T).
\]
Since on the other hand, we have by definition
\[
X^T_{(\gam, m)}=\prod_{\gam\in\xi} (X^T_{\gam})^{m(\gam)},
\]
it is sufficient to prove
\[
X^T_\gam=\CaCh'_T(\eta_T(\gam))
\]
for $\gam$ an arc or a simple loop.  
This is trivial if $\gam \in T^\circ$, thus
we have to distinguish only two cases:

Case 1: $\gam$ is an arc which does not belong to $T$. 
In this case, the string
module $M_\gam$ is $\tau$-rigid and therefore
$$
\eta_T(\gam) = \overline{\cO_{\cM_\gam}},
$$
compare Theorem~\ref{thm:main7c}.
So our claim follows directly
from  Proposition~\ref{MSWProp-CC}.

Case 2: $\gam$ is a simple loop.  In this case $\eta_T(\gam)$ is the closure
of the union of a the orbits of a family of modules, namely
\[
\eta_T(\gam) = 
\overline{\bigcup_{\la\in\C^*} \cO_{(M_{\gam,\la},0)}}.
\]
In this case we have again by Proposition~\ref{MSWProp-CC}
$X^T_\gam = \CaCh'_{\wtA_T}((M_{\gam,\la},0))$ for all $\la\in\C^*$, and we
are done.
\end{proof}

By specializing the coefficients to $1$, the equality 
$\widetilde{\cB}_T = \widetilde{\cG}_T$ from Theorem~\ref{thm:10.9}
yields
$$
\cB_T = \cG_T. 
$$
Thus Theorem~\ref{thm:main8} is proved.


\section{An example}\label{sec:example}


Let $\SM$ be the sphere with three disks cut out,
and one marked point on
each boundary component.
In Figure~\ref{exfig1} we display 
a  triangulation $T$ of
$\SM$, where the arcs of $T$ are marked in green, together with a loop $\sigma$ in $\SM$.

\begin{figure}[H] 
\tikz{
\useasboundingbox (-7,-3.2) rectangle (7,3.2);
  \draw (0,0) circle[x radius=6, y radius=3];
  \filldraw[fill=lightgray] (-2.2,0) circle[radius=1.2];
  \filldraw[fill=lightgray] (2.2,0) circle[radius=1.2];
  \path (0,2.97) node(a){}
         (-2.2,0) ++(38:1.23) node(b){} 
         ( 2.2,0) ++(142:1.23) node(c){};
  \filldraw (a) circle[radius=.1]
            (b) circle[radius=.1]
            (c) circle[radius=.1];
  \draw[green]
 (a) .. controls (-11, 1.9) and (-2.9, -6.5) .. (c) node[pos=0.4,left]{1}    
 (a) .. controls (- 11.5, 1.6) and (2.1, -5.1) .. (b) node[pos=0.3,left]{2}
 (a) -- (b) node[pos=0.5,left]{3}
 (a) -- (c) node[pos=0.5,right]{4}
 (a) .. controls (10.3,1.8) and (-1.9,-5.8) .. (c) node[pos=0.3,right]{5}
 (b) -- (c) node[pos=0.5,below]{6};
\draw[red]
(0,0) .. controls (4.7,4.7) and (4.7,-4.7) .. (0,0) node[pos=.3,above]{$\si$}
(0,0) .. controls (-4.7,4.7) and (-4.7,-4.7) .. (0,0);
}
\caption{Triangulation \greenx{$T$} of $\SM$ and loop \redx{$\si$}}
\label{exfig1}
\end{figure}

\begin{figure}[ht]
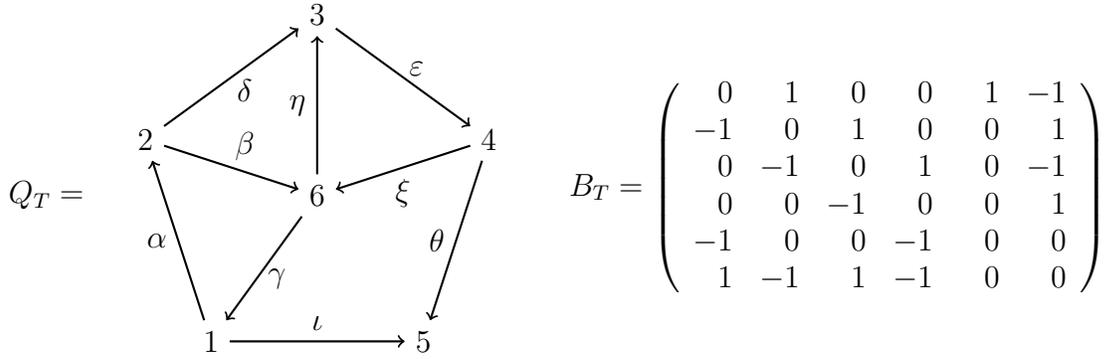
 
\tikz[baseline=0cm, scale=1.2]{
 \node(0) at (-3,0) {$Q_T=$}; 
 \node(1) at (234:2) {1};
 \node(2) at (162:2) {2};
 \node(3) at ( 90:2) {3};
 \node(4) at ( 18:2) {4};
 \node(5) at (306:2) {5};
 \node(6) at ( 0, 0) {6};
 \draw[->,thick] (1)--(2) node[pos=0.5,left]{$\alpha$}; 
 \draw[->,thick] (2)--(6) node[pos=0.6,above]{$\beta$};
 \draw[->,thick] (6)--(1) node[pos=0.6,right]{$\gamma$};
 \draw[->,thick] (2)--(3) node[pos=0.6,below]{$\delta$};
 \draw[->,thick] (6)--(3) node[pos=0.5,left]{$\eta$};
 \draw[->,thick] (3)--(4) node[pos=0.6,above]{$\eps$};
 \draw[->,thick] (4)--(6) node[pos=0.5,below]{$\xi$};  
 \draw[->,thick] (4)--(5) node[pos=0.5,left]{$\theta$};
 \draw[->,thick] (1)--(5) node[pos=0.5,above]{$\iota$}; 
}\qquad
$ B_T=\left(\begin{array}{rrrrrr}
  0 & 1 & 0 & 0 &\phantom{-}1 &-1\\
 -1 & 0 & 1 & 0 & 0 & 1\\
  0 &-1 & 0 & 1 & 0 &-1\\
  0 & 0 &-1 & 0 & 0 & 1\\
 -1 & 0 & 0 &-1 & 0 & 0\\
  1 &-1 & 1 &-1 & 0 & 0
\end{array}\right)$
\caption{Quiver $Q_T$ and signed adjacency matrix $B_T$ (MSW-convention)}
\label{exfig2}
\end{figure}

It is easy to read off the quiver $Q_T$ (following our convention) and the
signed adjacency matrix $B_T$ (following the convention of~\cite{MSW2}).
Both are shown in Figure~\ref{exfig2}.
Recall that with these convention in place we have $Q_T=Q(-B_T)$.

Musiker, Williams and Schiffler \cite{MSW2} associate to each loop $\si$
a band graph $G=G_{T,\si}$ with respect to a triangulation $T$. 
In our
example, we obtain the band graph $G$ displayed in
Figure~\ref{exfig3}.
Note that $G$ has $m=7$ tiles, corresponding
to the $7$ intersections of $\si$ with the edges of $T$.
The thick edges of $G$   correspond
to the distinguished good matching $P_-$.  
Note that the two extremal edges
have to be identified along the corresponding arrows.

\begin{figure}[H]    
\tikz[scale=1.3,baseline=1cm,>= triangle 45]{
\useasboundingbox (-.5,0) rectangle (4,3.2);
\draw[line width=.1cm] (0,1)--(0,2) node[pos=.5,left]{3}
                   (0,3)--(1,3) node[pos=.5,above]{5} 
                   (2,3)--(3,3) node[pos=.5,above]{6}
                   (4,3)--(5,3) node[pos=.5,above]{3}
                   (5,2)--(4,2) node[pos=.5,below]{6}
                   (3,2)--(2,2) node[pos=.5,below]{5} 
                   (1,2)--(1,1) node[pos=.5,right]{1};
 \draw[thin] (0,0)--(0,1) node[pos=.5, left]{2}
              (0,2)--(0,3) node[pos=.5, left]{6}
              (1,3)--(2,3) node[pos=.5, above]{4}
              (3,3)--(4,3) node[pos=.5, above]{2}
              (4,2)--(3,2) node[pos=.5,below]{4}
              (2,2)--(1,2) node[pos=.5,below]{1}
              (1,1)--(1,0) node[pos=.5,right]{6}
              (0,0)--(1,0) (5,3)--(5,2);   
  \draw[<-]  (5,2.6)--(5,2);               
  \draw[->]  (0,0)--(.6,0);                
  \draw[thin] (1,0)--(0,1) node[pos=.5, fill=white]{3}
              (0,1)--(1,1) node[pos=.5, fill=white]{4}
              (1,1)--(0,2) node[pos=.5, fill=white]{6}
              (0,2)--(1,2) node[pos=.5, fill=white]{2}
              (1,2)--(0,3) node[pos=.5, fill=white]{1}
              (1,2)--(1,3)
             (1,3)--(2,2) node[pos=.5, fill=white]{5} 
              (2,2)--(2,3)
              (2,3)--(3,2) node[pos=.5, fill=white]{4}
              (3,2)--(3,3) node[pos=.5, fill=white]{3}
              (3,3)--(4,2) node[pos=.5, fill=white]{6}
              (4,2)--(4,3) node[pos=.5, fill=white]{1}
              (4,3)--(5,2) node[pos=.5, fill=white]{2};
 \node at (0,0)[above right]{$+$};
 \node at (0,1)[above right]{$-$};
 \node at (0,2)[above right]{$+$};
 \node at (1,2)[above right]{$-$};
 \node at (2,2)[above right]{$+$};
 \node at (3,2)[above right]{$-$};
 \node at (4,2)[above right]{$+$};
}
$\displaystyle{\frac{x(P_-)}{\cross(T,\si)}
=\frac{x_3 x_1 x_5 x_6 x_5 x_3 x_6}{x_3 x_6 x_1 x_5 x_4 x_6 x_2}
=\frac{x_5 x_3}{x_4 x_2}=\ux^{\bs(\si)}}$
\caption{Band graph $G=G_{T,\si}$ with $P_-$ (thick edges) and
         $\ux^{\bs(\si)}$}
\label{exfig3}
\end{figure}

Recall from Section~\ref{subsec:10.3} that Musiker, Schiffler and Williams associate
to $G$ a Hasse quiver $Q_G$, which is opposite to our coefficient quiver
$Q_\si^T$  of the band module $M_{\si,\lambda}$ for $\lambda\in\C^*$, see Remarks~\ref{Rem:10.1} and \ref{MSW:RemCoeff}. 
We display the coefficient quiver
$Q_\si^T$ in Figure~\ref{exfig4}. Note that the two encircled vertices have to
be identified.

\begin{figure}[H]
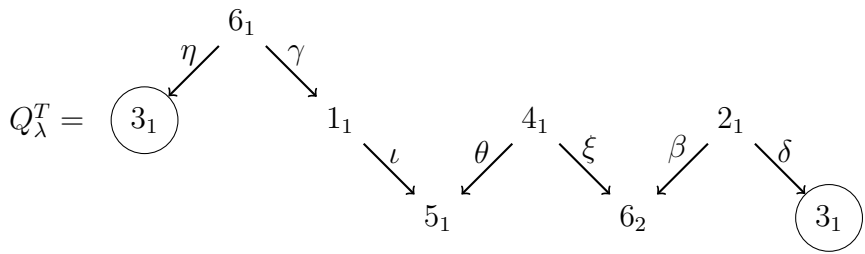
 
\tikz[baseline=0cm, scale=1.3]{
  \node(0) at (-1,0) {$Q_\lambda^T=$};
  \node(1) at (0,0)[circle,draw=black] {$3_1$};
  \node(2) at (1,1) {$6_1$};
  \node(3) at (2,0) {$1_1$};
  \node(4) at (3,-1){$5_1$};
  \node(5) at (4,0) {$4_1$};
  \node(6) at (5,-1){$6_2$};
  \node(7) at (6,0) {$2_1$};
  \node(8) at (7,-1)[circle,draw=black]{$3_1$};
\draw[->,thick] (2)--(1) node[pos=0.6,above]{$\eta$};  
\draw[->,thick] (2)--(3) node[pos=0.6,above]{$\gamma$};
\draw[->,thick] (3)--(4) node[pos=0.6,above]{$\iota$};
\draw[->,thick] (5)--(4) node[pos=0.6,above]{$\theta$};
\draw[->,thick] (5)--(6) node[pos=0.6,above]{$\xi$};
\draw[->,thick] (7)--(6) node[pos=0.6,above]{$\beta$};
\draw[->,thick] (7)--(8) node[pos=0.6,above]{$\delta$};  
}
\caption{Coefficient quiver $Q_\si^T$}
\label{exfig4}
\end{figure}

Thus the order coideals of $Q_\si^T$ (i.e. coordinate
factor modules of $M_{\si,\lambda}$) are in bijection with the good matchings
of $G$.  
More precisely, the tiles which are enclosed
by the symmetric difference $P\ominus P_-$ for a good matching $P$ are
identified with a basis of the corresponding coordinate factor module.
Finally we display in Figure~\ref{exfig5a},
three of the 27 good matchings of $G$.  In each case the edges of the matching
$P$ are  highlighted in orange, whilst the tiles which are enclosed
by $P\ominus P_-$ are highlighted in yellow. Moreover, we show in each case
the contribution of $P$ to $X_\si^T$. 

\begin{figure}[H]
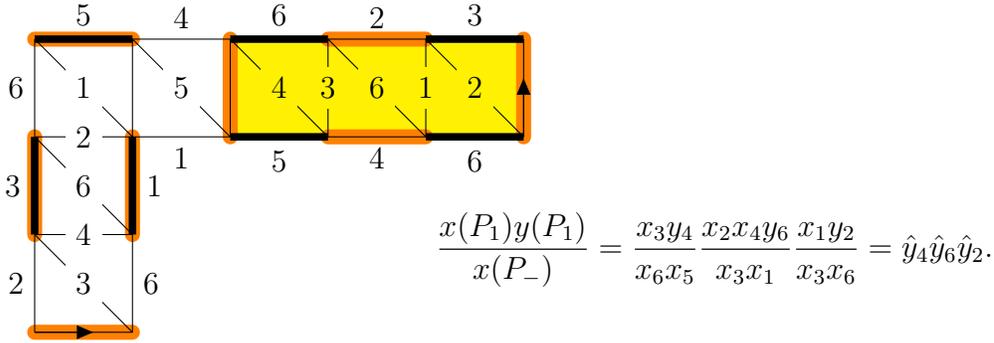

  \ContinuedFloat*
\tikz[scale=1.3, baseline=1cm,>= triangle 45]{
\useasboundingbox (-.5,0) rectangle (4,3.8);
\draw[fill, color=yellow] (2,2) rectangle (5,3);
\draw[line width=.2cm, color=orange, line cap=round]
(0,0)--(1,0) (0,1)--(0,2) (1,1)--(1,2) (0,3)--(1,3) (2,2)--(2,3)
(3,2)--(4,2) (3,3)--(4,3) (5,2)--(5,3);
\draw[line width=.1cm] (0,1)--(0,2) node[pos=.5,left]{3}
                   (0,3)--(1,3) node[pos=.5,above]{5}
                   (2,3)--(3,3) node[pos=.5,above]{6}
                   (4,3)--(5,3) node[pos=.5,above]{3}
                   (5,2)--(4,2) node[pos=.5,below]{6}
                   (3,2)--(2,2) node[pos=.5,below]{5}
                   (1,2)--(1,1) node[pos=.5,right]{1};
  \draw[thin] (0,0)--(0,1) node[pos=.5, left]{2}
              (0,2)--(0,3) node[pos=.5, left]{6}
              (1,3)--(2,3) node[pos=.5, above]{4}
              (3,3)--(4,3) node[pos=.5, above]{2}
              (4,2)--(3,2) node[pos=.5,below]{4}
              (2,2)--(1,2) node[pos=.5,below]{1}
              (1,1)--(1,0) node[pos=.5,right]{6}
(0,0)--(1,0) (5,3)--(5,2);   
  \draw[<-]  (5,2.6)--(5,2);               
  \draw[->]  (0,0)--(.6,0);                
  \draw[thin] (1,0)--(0,1) node[pos=.5, fill=white]{3}
              (0,1)--(1,1) node[pos=.5, fill=white]{4}
              (1,1)--(0,2) node[pos=.5, fill=white]{6}
              (0,2)--(1,2) node[pos=.5, fill=white]{2}
              (1,2)--(0,3) node[pos=.5, fill=white]{1}
              (1,2)--(1,3)
              (1,3)--(2,2) node[pos=.5, fill=white]{5}
              (2,2)--(2,3)
              (2,3)--(3,2) node[pos=.5, fill=yellow]{4}
              (3,2)--(3,3) node[pos=.5, fill=yellow]{3}
              (3,3)--(4,2) node[pos=.5, fill=yellow]{6}
              (4,2)--(4,3) node[pos=.5, fill=yellow]{1}
              (4,3)--(5,2) node[pos=.5, fill=yellow]{2};
}
${\displaystyle
  \frac{x(P_1)y(P_1)}{x(P_-)}
  =\frac{x_3 y_4}{x_6 x_5} \frac{x_2 x_4 y_6}{x_3 x_1} \frac{x_1 y_2}{x_3 x_6}
  = \hy_4 \hy_6 \hy_2}$.\hfill\phantom{X}
\caption{ Good matching $P_1$ and its contribution to $X_\si^T$}
\label{exfig5a}
\end{figure}

\begin{figure}[H]
  \ContinuedFloat
\tikz[scale=1.3, baseline=1cm,>= triangle 45]{
\useasboundingbox (-.5,0) rectangle (4,3.8);
\draw[fill, color=yellow] (2,2) rectangle (3,3);
\draw[fill, color=yellow] (4,2) rectangle (5,3);
\draw[line width=.2cm, color=orange, line cap=round]
(0,0)--(1,0) (0,1)--(0,2) (1,1)--(1,2) (0,3)--(1,3) (2,2)--(2,3)
(3,2)--(3,3) (4,2)--(4,3) (5,2)--(5,3);
\draw[line width=.1cm] (0,1)--(0,2) node[pos=.5,left]{3}
                   (0,3)--(1,3) node[pos=.5,above]{5}
                   (2,3)--(3,3) node[pos=.5,above]{6}
                   (4,3)--(5,3) node[pos=.5,above]{3}
                   (5,2)--(4,2) node[pos=.5,below]{6}
                   (3,2)--(2,2) node[pos=.5,below]{5}
                   (1,2)--(1,1) node[pos=.5,right]{1};
  \draw[thin] (0,0)--(0,1) node[pos=.5, left]{2}
              (0,2)--(0,3) node[pos=.5, left]{6}
              (1,3)--(2,3) node[pos=.5, above]{4}
              (3,3)--(4,3) node[pos=.5, above]{2}
              (4,2)--(3,2) node[pos=.5,below]{4}
              (2,2)--(1,2) node[pos=.5,below]{1}
              (1,1)--(1,0) node[pos=.5,right]{6}
(0,0)--(1,0) (5,3)--(5,2);   
  \draw[<-]  (5,2.6)--(5,2);               
  \draw[->]  (0,0)--(.6,0);                
  \draw[thin] (1,0)--(0,1) node[pos=.5, fill=white]{3}
              (0,1)--(1,1) node[pos=.5, fill=white]{4}
              (1,1)--(0,2) node[pos=.5, fill=white]{6}
              (0,2)--(1,2) node[pos=.5, fill=white]{2}
              (1,2)--(0,3) node[pos=.5, fill=white]{1}
              (1,2)--(1,3)
              (1,3)--(2,2) node[pos=.5, fill=white]{5}
              (2,2)--(2,3)
              (2,3)--(3,2) node[pos=.5, fill=yellow]{4}
              (3,2)--(3,3) node[pos=.5, inner xsep=0pt, fill=orange]{3}
              (3,3)--(4,2) node[pos=.5, fill=white]{6}
              (4,2)--(4,3) node[pos=.5, inner xsep=0pt, fill=orange]{1}
              (4,3)--(5,2) node[pos=.5, fill=yellow]{2};
}
${\displaystyle
  \frac{x(P_2)y(P_2)}{P_-}=\frac{x_3 y_4}{x_5 x_6} \frac{x_1 y_2}{x_3 x_6}
  =\hy_4 \hy_2}$.\hfill\phantom{X}
\caption{Good matching $P_2$ and its contribution to $X_\si^T$}
\end{figure}

\begin{figure}[H]
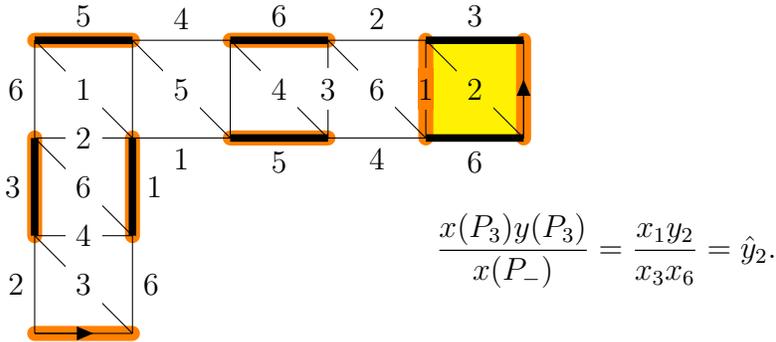

    \ContinuedFloat
\tikz[scale=1.3, baseline=1cm,>= triangle 45]{
\useasboundingbox (-.5,0) rectangle (4,3.8);
\draw[fill, color=yellow] (4,2) rectangle (5,3);
\draw[line width=.2cm, color=orange, line cap=round]
(0,0)--(1,0) (0,1)--(0,2) (1,1)--(1,2) (0,3)--(1,3) (2,2)--(3,2)
(2,3)--(3,3) (4,2)--(4,3) (5,2)--(5,3);
\draw[line width=.1cm] (0,1)--(0,2) node[pos=.5,left]{3}
                   (0,3)--(1,3) node[pos=.5,above]{5}
                   (2,3)--(3,3) node[pos=.5,above]{6}
                   (4,3)--(5,3) node[pos=.5,above]{3}
                   (5,2)--(4,2) node[pos=.5,below]{6}
                   (3,2)--(2,2) node[pos=.5,below]{5}
                   (1,2)--(1,1) node[pos=.5,right]{1};
  \draw[thin] (0,0)--(0,1) node[pos=.5, left]{2}
              (0,2)--(0,3) node[pos=.5, left]{6}
              (1,3)--(2,3) node[pos=.5, above]{4}
              (3,3)--(4,3) node[pos=.5, above]{2}
              (4,2)--(3,2) node[pos=.5,below]{4}
              (2,2)--(1,2) node[pos=.5,below]{1}
              (1,1)--(1,0) node[pos=.5,right]{6}
(0,0)--(1,0) (5,3)--(5,2);   
  \draw[<-]  (5,2.6)--(5,2);               
  \draw[->]  (0,0)--(.6,0);                
  \draw[thin] (1,0)--(0,1) node[pos=.5, fill=white]{3}
              (0,1)--(1,1) node[pos=.5, fill=white]{4}
              (1,1)--(0,2) node[pos=.5, fill=white]{6}
              (0,2)--(1,2) node[pos=.5, fill=white]{2}
              (1,2)--(0,3) node[pos=.5, fill=white]{1}
              (1,2)--(1,3)
              (1,3)--(2,2) node[pos=.5, fill=white]{5}
              (2,2)--(2,3)
              (2,3)--(3,2) node[pos=.5, fill=white]{4}
              (3,2)--(3,3) node[pos=.5, fill=white]{3}
              (3,3)--(4,2) node[pos=.5, fill=white]{6}
              (4,2)--(4,3) node[pos=.5, inner xsep=0pt, fill=orange]{1} 
              (4,3)--(5,2) node[pos=.5, fill=yellow]{2}; 
            }
$\displaystyle \frac{x(P_3) y(P_3)}{x(P_-)}=\frac{x_1 y_2}{x_3 x_6}=\hy_2$.
\hfill\phantom{X}
\caption{Good matching $P_3$ and its contribution to $X_\si^T$}
\end{figure}

The relation between perfect matchings and coordinate submodules of string modules has been also studied in a more general setup by
Canakci and Schroll \cite{CS}.

\bigskip
{\parindent0cm \bf Acknowledgements.}\,
The first named author acknowledges partial support from 
CONACyT grant no. 239255.
The second author gratefully acknowledges the support he received from a C\'atedra Marcos Moshinsky and the grants 
CONACyT-238754 and PAPIIT-IN112519.
The third author thanks the Sonderforschungsbereich/Transregio 
SFB 45 for support.
He was also funded by the Deutsche Forschungsgemeinschaft
(DFG, German Research Foundation) under Germany's Excellence Strategy - GZ 2047/1, Projekt-ID 390685813.
He thanks the Instituto de Matem\'atics  of the UNAM (Mexico City)
for two weeks of hospitality in September 2019, where part of this work
was done.
The first and second author thank the Mathematical Institute of the University of Bonn for two weeks of hospitality in January 2020.
We thank Ryan Kinser for helpful comments on an earlier version of this article.
We are indepted to the anonymous referee for suggesting improvements of the proofs of Proposition~\ref{prop:orbitclosure} and Lemma~\ref{lem:simplesummands}.



\begin{thebibliography}{999}


\bibitem[AIR]{AIR}
T. Adachi, O. Iyama, I. Reiten,
\emph{$\tau$-tilting theory}. 
Compos. Math. 150 (2014), no. 3, 415--452.

\bibitem[A]{A}
D. Allegretti,
\emph{Categorified canonical bases and framed BPS states}.
Preprint (2018), 51 pp., 
arXiv:1806.10394v2

\bibitem[AB]{AB}
C. Amiot, T. Br\"ustle,
\emph{Derived equivalences between skew-gentle algebras using orbifolds}.
Preprint (2019), 46 pp., arXiv:1912.04367

\bibitem[ABCP]{ABCP}
I. Assem, T. Br\"ustle, G. Charbonneau-Jodoin, 
P.-G. Plamondon, 
\emph{Gentle algebras arising from surface triangulations}. 
Algebra Number Theory 4 (2010), no. 2, 201--229. 

\bibitem[AS]{AS}
I. Assem, A. Skowro\'nski,
\emph{Iterated tilted algebras of type $\widetilde{A}_n$}. 
Math. Z. 195 (1987), no. 2, 269--290.

\bibitem[ASS]{ASS}
I. Assem, D. Simson, A. Skowro\'nski,
\emph{Elements of the representation theory of associative algebras}, 
Vol. 1. Techniques of representation theory. London Mathematical Society Student Texts, 65. Cambridge University Press, Cambridge, 2006. x+458 pp. 

\bibitem[ARS]{ARS}
M. Auslander, I. Reiten, S. Smal\o,
\emph{Representation theory of Artin algebras}. 
Corrected reprint of the 1995 original. Cambridge Studies in Advanced Mathematics, 36. Cambridge University Press, Cambridge, 1997. 
xiv+425 pp.

\bibitem[BS]{BS}
K. Baur, S. Schroll,
\emph{Higher extensions for gentle algebras}.
Preprint (2019), 18 pp., arXiv:1906.05257

\bibitem[B]{B}
T. Br\"ustle,
\emph{Kit algebras}.
J. Algebra 240 (2001), no. 1, 1--24. 

\bibitem[BZ]{BZ}
T. Br\"ustle, J. Zhang,
\emph{On the cluster category of a marked surface without punctures}.
Algebra Number Theory 5 (2011), no. 4, 529--566.

\bibitem[BR]{BR}
M.C.R. Butler, C.M. Ringel,
\emph{Auslander-Reiten sequences with few middle terms and applications to string algebras}. 
Comm. Algebra 15 (1987), no. 1-2, 145--179. 

\bibitem[CS]{CS}
I. Canakci, S. Schroll,
\emph{Lattice bijections for string modules, snake graphs and the weak Bruhat order}.
Preprint (2018), 17 pp., arXiv:1811.06064

\bibitem[C]{C}
A. Carroll,
\emph{Generic modules for string algebras}.
J. Algebra 437 (2015), 177--201. 

\bibitem[CC]{CC}
A. Carroll, C. Chindris, 
\emph{On the invariant theory for acyclic gentle algebras}. 
Trans. Amer. Math. Soc. 367 (2015), no. 5, 3481--3508.

\bibitem[CCKW]{CCKW}
A. Carroll, C. Chindris, R. Kinser, J. Weyman,
\emph{Moduli Spaces of Representations of Special Biserial Algebras}. Int. Math. Res. Not. IMRN 2020, no. 2, 403--421.

\bibitem[CW]{CW}
A. Carroll, J. Weyman,
\emph{Semi-invariants for gentle algebras}. 
Noncommutative birational geometry, representations and combinatorics, 111--136, 
Contemp. Math., 592, Amer. Math. Soc., Providence, RI, 2013. 

\bibitem[CLFS]{CLFS}
G. Cerulli Irelli, D. Labardini-Fragoso, J. Schr\"oer,
\emph{Caldero-Chapoton algebras}.
Trans. Amer. Math. Soc. 367 (2015), no. 4, 2787--2822. 

\bibitem[CB1]{CB1}
W. Crawley-Boevey,
\emph{Maps between representations of zero-relation algebras}. 
J. Algebra 126 (1989), no. 2, 259--263.

\bibitem[CB2]{CB2}
W. Crawley-Boevey,
\emph{On tame algebras and bocses}. 
Proc. London Math. Soc. (3) 56 (1988), 451--483.

\bibitem[CBS]{CBS}
W. Crawley-Boevey, J. Schr\"oer,
\emph{Irreducible components of varieties of modules}.
 J. Reine Angew. Math. 553 (2002), 201--220.

\bibitem[DS]{DS}
C. De Concini, E. Strickland,
\emph{On the variety of complexes}.
Adv. in Math. 41 (1981), no. 1, 57--77. 

\bibitem[DWZ1]{DWZ1}
H. Derksen, J. Weyman, A. Zelevinsky,
\emph{Quivers with potentials and their representations. I. Mutations}. 
Selecta Math. (N.S.) 14 (2008), no. 1, 59--119. 

\bibitem[DWZ2]{DWZ2}
H. Derksen, J. Weyman, A. Zelevinsky,
\emph{Quivers with potentials and their representations II: 
applications to cluster algebras}. 
J. Amer. Math. Soc. 23 (2010), no. 3, 749--790.

\bibitem[FG]{FG}
V. Fock, A. Goncharov, 
\emph{Dual Teichm\"uller and lamination spaces}. 
Handbook of Teichm\"uller theory. Vol. I, 647--684, 
IRMA Lect. Math. Theor. Phys., 11, Eur. Math. Soc., Z\"urich, 2007. 

\bibitem[FST]{FST}
S. Fomin, M. Shapiro, D. Thurston,
\emph{Cluster algebras and triangulated surfaces. I. Cluster complexes}. 
Acta Math. 201 (2008), no. 1, 83--146.

\bibitem[FT]{FT}
S. Fomin, D. Thurston, 
\emph{Cluster algebras and triangulated surfaces Part II: Lambda lengths}. 
Mem. Amer. Math. Soc. 255 (2018), no. 1223, v+97 pp. 

\bibitem[FZ]{FZ}
S. Fomin, A. Zelevinsky, 
\emph{Cluster algebras. I. Foundations}. 
J. Amer. Math. Soc. 15 (2002), no. 2, 497--529.

\bibitem[Ga]{Ga}
P. Gabriel, 
\emph{Finite representation type is open}.
Proceedings of the International Conference on Rep- resentations of Algebras (Carleton Univ., Ottawa, Ont., 1974), Paper No. 10, 23pp. Carleton Math. Lecture Notes, No. 9, Carleton Univ., Ottawa, Ont., 1974.

\bibitem[GLFS]{GLFS}
C. Gei{\ss}, D. Labardini-Fragoso, J. Schr\"oer,
\emph{Generic bases of surface cluster algebras with geometric coefficients}, 36 pp., Preprint in preparation (2020).

\bibitem[GLS]{GLS}
C. Gei{\ss}, B. Leclerc, J. Schr\"oer,
\emph{Generic bases for cluster algebras and the chamber Ansatz}. 
J. Amer. Math. Soc. 25 (2012), no. 1, 21--76.

\bibitem[G]{G}
C. Gei{\ss},
\emph{Geometric methods in representation theory of finite-dimensional algebras}. 
Representation theory of algebras and related topics (Mexico City, 1994), 53--63, 
CMS Conf. Proc., 19, Amer. Math. Soc., Providence, RI, 1996. 

\bibitem[GP]{GP}
C. Gei{\ss},  J.A. de la Pe{\~n}a,
\emph{On the deformation theory of finite-dimensional algebras}.
Manuscripta Math. 88 (1995), no. 2, 191--208. 

\bibitem[Go]{Go}
N. Gonciulea,
\emph{Singular Loci of Varieties of Complexes. II}.
J. Algebra 235 (2001), no. 2, 547--558.

\bibitem[HKK]{HKK}
F. Haiden, L. Katzarkov, M. Kontsevich,
\emph{Flat surfaces and stability structures}.
Publ. Math. Inst. Hautes \'Etudes Sci. 126 (2017), 247--318.

\bibitem[H]{H}
R. Hartshorne,
\emph{Algebraic geometry}. 
Graduate Texts in Mathematics, No. 52. Springer-Verlag, 
New York-Heidelberg, 1977. xvi+496 pp. 

\bibitem[Hau]{Hau}  
N. Haupt,
\emph{Euler characteristics of quiver Grassmannians and Ringel-Hall algebras of string algebras}.  
Algebr. Represent. Theory 15 (2012), no. 4, 755--793. 

\bibitem[K]{K}
H. Krause, 
\emph{Maps between tree and band modules}.
J. Algebra 137 (1991), no. 1, 186--194.

\bibitem[LF1]{LF1}
D. Labardini-Fragoso, 
\emph{Quivers with potentials associated to triangulated surfaces}.
Proc. London Math. Soc, 98 (2009), 797--839. 

\bibitem[LF2]{LF2}
D. Labardini-Fragoso,
\emph{Quivers with potentials associated with triangulations of Riemann surfaces}. 
Ph.D. thesis, Northeastern University (2010).

\bibitem[L]{L}
V. Lakshmibai,
\emph{Singular loci of varieties of complexes}. 
Commutative algebra, homological algebra and representation theory (Catania/Genoa/Rome, 1998). 
J. Pure Appl. Algebra 152 (2000), no. 1-3, 217--230.

\bibitem[LP]{LP}
Y. Lekili, A. Polishchuk,
\emph{Derived equivalences of gentle algebras via Fukaya categories}.
Preprint (2018), 35 pp., arXiv:1801.06370

\bibitem[Mu1]{Mu1}
G. Muller, 
\emph{Locally acyclic cluster algebras}.
Adv. Math. 233 (2013) 207--247.

\bibitem[Mu2]{Mu2}
G. Muller, 
\emph{$\cA=\cU$ for locally acyclic cluster algebras}.
SIGMA Symmetry Integrability Geom. Methods Appl. 10 (2014), Paper 094, 8 pp. 

\bibitem[MSW1]{MSW1}
G. Musiker, R. Schiffler, L. Williams, 
\emph{Positivity for cluster algebras from surfaces}. 
Adv. Math. 227 (2011), 2241--2308.

\bibitem[MSW2]{MSW2}
G. Musiker, R. Schiffler, L. Williams, 
\emph{Bases for cluster algebras from surfaces}.
Compos. Math. 149 (2013), no. 2, 217--263. 

\bibitem[OPS]{OPS}
S. Opper, P.-G. Plamondon, S. Schroll,
\emph{A geometric model for the derived category of gentle algebras}.
Preprint (2018), 41 pp., arXiv:1801.09659v5

\bibitem[P1]{P1}
P.-G. Plamondon,
\emph{Generic bases for cluster algebras from the cluster category}.
Int. Math. Res. Not. IMRN 2013, no. 10, 2368--2420. 

\bibitem[P2]{P2}
P.-G. Plamondon,
\emph{$\tau$-tilting finite gentle algebras are representation-finite}.
Pacific J. Math. (to appear), 
Preprint (2018), 7 pp. arXiv:1809.06313

\bibitem[Q]{Q}
F. Qin,
\emph{Bases for upper cluster algebras and tropical points}.
Preprint (2019), 45 pp., 
arXiv:1902.09507

\bibitem[R]{R}
C.M. Ringel,
\emph{Tame algebras and integral quadratic forms}. 
Lecture Notes in Mathematics, 1099. Springer-Verlag, Berlin, 1984. xiii+376 pp.

\bibitem[Sch]{Sch}
J. Schr\"oer, 
\emph{Modules without self-extensions over gentle algebras}. 
J. Algebra 216 (1999), no. 1, 178--189.

\bibitem[Sh1]{Sh1}
I. Shafarevich, 
\emph{Basic algebraic geometry. 1. Varieties in projective space}. 
Third edition. Translated from the 2007 third Russian edition. Springer, Heidelberg, 2013. xviii+310 pp. 

\bibitem[Sh2]{Sh2}
I. Shafarevich, 
\emph{Basic algebraic geometry. 2. 
Schemes and complex manifolds}. 
Third edition. Translated from the 2007 third Russian edition by Miles Reid. Springer, Heidelberg, 2013. xiv+262 pp. 

\bibitem[St]{St}
E. Strickland, 
\emph{On the conormal bundle of the determinantal variety}. 
J. Algebra 75 (1982), no. 2, 523--537. 

\bibitem[V]{V}
D. Voigt, 
\emph{Induzierte Darstellungen in der Theorie der endlichen, algebraischen Gruppen}. 
Lecture Notes in Mathematics, Vol. 592. Springer-Verlag, 
Berlin-New York, 1977. iv+413 pp.

\bibitem[WW]{WW}
B. Wald, J. Waschb\"usch, 
\emph{Tame biserial algebras},
J. Algebra 95 (1985), no. 2, 480--500. 

\bibitem[Z]{Z}
G. Zwara,
\emph{Degenerations for modules over representation-finite 
algebras}.
Proc. Amer. Math. Soc. 127 (1999), 1313--1322.

\end{thebibliography}
\end{document}